\documentclass[12pt]{amsart}

\usepackage{amssymb}
\usepackage{mathrsfs}
\usepackage{xypic}
\usepackage{graphicx}
\usepackage{psfrag}
  
\theoremstyle{plain}
\newtheorem{theorem}{Theorem}[subsection]
\newtheorem{lemma}[theorem]{Lemma}
\newtheorem{proposition}[theorem]{Proposition}
\newtheorem{corollary}[theorem]{Corollary}

\newtheorem{question}[theorem]{Question}

\theoremstyle{definition}
\newtheorem{definition}[theorem]{Definition}

\theoremstyle{remark}
\newtheorem{remark}[theorem]{Remark}

\newtheorem{example}[theorem]{Example}
\newtheorem{notation}[theorem]{Notation}

\newtheorem{observations}[theorem]{Observations}

\renewenvironment{proof}[1]{\vspace*{.1in}\noindent{\bf
Proof{#1}. \/}}{\qed\vspace{3ex}} 
		     
\newcommand{\bbP}{{\mathbb P}}
\newcommand{\bbQ}{{\mathbb Q}}
\newcommand{\bbN}{{\mathbb N}}
\newcommand{\bbR}{{\mathbb R}}
\newcommand{\bbZ}{{\mathbb Z}}
\newcommand{\bbH}{{\mathbb H}}
\newcommand{\bbE}{{\mathbb E}}

\newcommand{\bA}{{\mathbf A}}
\newcommand{\bB}{{\mathbf B}}
\newcommand{\bI}{{\mathbf I}}
\newcommand{\bH}{{\mathbf H}}
\newcommand{\bF}{{\mathbf F}}
\newcommand{\bD}{{\mathbf D}}
\newcommand{\bE}{{\mathbf E}}

\newcommand{\be}{{\mathbf e}}

\newcommand{\bu}{{\mathbf u}}
\newcommand{\bs}{{\mathbf s}}

\newcommand{\cA}{{\mathcal A}}
\newcommand{\cS}{{\mathcal S}}
\newcommand{\cF}{{\mathcal F}}

\newcommand{\cM}{{\mathcal M}}
\newcommand{\cN}{{\mathcal N}}
\newcommand{\cP}{{\mathcal P}}
\newcommand{\cT}{{\mathcal T}}
\renewcommand{\cR}{{\mathcal R}}

\newcommand{\cc}{{Z}}
\newcommand{\ccc}{{X}}
\newcommand{\cp}{{x}}

\newcommand{\sK}{{\mathscr K}}
\newcommand{\sC}{{\mathscr C}}
\newcommand{\sS}{{\mathscr S}}
\newcommand{\sD}{{\mathscr D}}

\newcommand{\card}{\operatorname{Card}}
\newcommand{\Lk}{\operatorname{Lk}}
\newcommand{\clstar}{\overline{\operatorname{St}}}
\renewcommand{\star}{\operatorname{St}}

\newcommand{\Span}{\operatorname{Span}}
\newcommand{\Aut}{\operatorname{Aut}}
\newcommand{\Stab}{\operatorname{Stab}}
\newcommand{\Bd}{\partial}
\newcommand{\supp}{\operatorname{Supp}}
\newcommand{\Id}{\operatorname{Id}}

\renewcommand{\hat}[1]{\widehat{#1}}
\renewcommand{\bar}[1]{\overline{#1}}
\renewcommand{\tilde}[1]{\widetilde{#1}}

\catcode`\@=11
\def\@secnumfont{\bfseries}
\def\section{\@startsection{section}{1}%
  \z@{.7\linespacing\@plus\linespacing}{.5\linespacing}%
  {\normalfont\centering\bfseries}}
\def\subsection{\@startsection{subsection}{2}%
  \z@{.5\linespacing\@plus.7\linespacing}{-.5em}%
  {\normalfont\bfseries}}
\catcode`\@12

\title{Fundamental Groups of Blow-ups}
\author{M.~Davis, T.~Januszkiewicz, and R.~Scott}
\thanks{First author partially supported by NSF grant DMS9803374.\\
\indent Second author partially supported by KBN grants P03A 023 14 and 5 P03A 035 20.} 
\address{The Ohio State University}
\email{mdavis@math.ohio-state.edu}
\address{Wroc{\l}aw University and IM PAN}
\email{tjan@math.uni.wroc.pl}
\address{Santa Clara University} 
\email{rscott@math.scu.edu}

\begin{document}

\begin{abstract}
Many examples of nonpositively curved closed manifolds arise as
blow-ups of projective hyperplane arrangements.  If the hyperplane
arrangement is associated to a finite reflection group $W$, and the
blow-up locus is $W$-invariant, then the resulting manifold will admit
a cell decomposition whose maximal cells are all combinatorially
isomorphic to a given convex polytope $P$.  In other words, $M$ admits
a tiling with tile $P$.  The universal covers of such examples yield  
tilings of $\bbR^n$ whose symmetry groups are generated by involutions
but are not, in general, reflection groups.  We begin a study of these
``mock reflection groups'', and develop a theory of tilings that
includes the examples coming from blow-ups and that generalizes the
corresponding theory of reflection tilings.  We apply our general
theory to classify the examples coming from blow-ups in the case where
the tile $P$ is either the permutohedron or the associahedron.   
\end{abstract}

\maketitle

\section{Introduction}\label{s:intro}

Suppose $M^n$ is a connected closed manifold equipped with a cubical
cell structure.  (In other words, $M^n$ is homeomorphic to a regular
cell complex in which each $k$-dimensional cell is combinatorially
isomorphic to a $k$-dimensional cube.)  It turns out that there is a
rich class of examples of such manifolds satisfying the following
three properties.
\begin{enumerate}
\item[(1)] There is a group $G$ of symmetries of the cellulation such that
the action of $G$ on the vertex set is simply transitive and such that
the stabilizer of each edge is cyclic of order $2$.
\item[(2)] In the dual cell structure on $M^n$ each top-dimensional cell is
combinatorially isomorphic to some given simple convex polytope, for
example, to a permutohedron or an associahedron.  (Such a
top-dimensional dual cell will be called a ``tile''.)
\item[(3)] The natural piecewise Euclidean metric on $M^n$ (in which each
combinatorial cube is isometric to a regular cube in Euclidean space)
is nonpositively curved.
\end{enumerate}
It follows from (2) and (3) that the universal cover $\tilde{M}^n$ is
homeomorphic to $\bbR^n$.  The cubical cell structure on $M^n$ lifts
to a cellulation of $\tilde{M}^n$ as does the dual cell structure.
Although these two cell structures on $M^n$ (the cubical one and its dual)
carry exactly the same combinatorial information, they correspond to two
distinct geometric pictures.  Throughout this paper we shall go back and
forth between these two pictures.  For example, property (1), that $G$ acts
simply transitively on the vertex set of the cubical cellulation, 
means that $G$
acts simply transitively on the set of $n$-dimensional dual cells.  So,
$\tilde{M}^n$ is ``tiled'' by isomorphic copies of such an $n$-dimensional
dual cell.

Let $A$ denote the group of all lifts of the $G$-action to
$\tilde{M}^n$.  Fix a vertex $x$ of the cubical structure on $\tilde{M}^n$.
By property (1)
each edge containing $x$ is flipped by a unique involution in $A$.
Since $\tilde{M}^n$ is connected, these involutions generate $A$ and the
$1$-skeleton of $\tilde{M}^n$ is the Cayley graph of $A$ with respect to
this set of generators.
Since $\tilde{M}^n$ is simply connected, a presentation for $A$ can be
derived by examining the $2$-cells that contain $x$ and the $2$-skeleton
of $\tilde{M}^n$ is the Cayley $2$-complex of this presentation.  (This is
explained in Section~\ref{ss:group-A} and \ref{s:tilings}.)
Furthermore, the fundamental group of $M^n$ is naturally identified
with the kernel of the epimorphism $A\rightarrow G$ induced by the
projection $\tilde{M}^n\rightarrow M^n$.  One of the purposes of this paper is
to initiate the study of such symmetry groups $A$.

This paper has two major thrusts:
\begin{itemize}
\item to describe a large class of examples of the above type (in Sections
\ref{s:intro} - \ref{s:dual}, \ref{s:permuto} and
\ref{s:associahedral-tilings}), and 
\item to develop a general theory of tilings and their symmetry groups
(in Sections~\ref{s:tilings} and \ref{s:linearity}).
\end{itemize}
We first give a rough description of the examples. 
The first examples are fairly standard and arise from actions of
right-angled reflection groups on manifolds.  (In this setting 
$\tilde{M}^n$ is the manifold, and $A$ is the reflection group.)  The
other examples that we discuss arise by performing an equivariant
blow-up procedure to (not necessarily right-angled) reflection group
actions and lifting to the universal cover.  In this case,
$\tilde{M}^n$ is the universal cover, and $A$ is the group of lifts of
the reflection group action.  An important guiding principle
underlying this paper is that the group actions in the blow-up setting
are tantalizingly similar to, but different from, reflection group
actions.  

Our reflection-type examples can be constructed as in \cite{D} or
\cite{DM}.  Given a simple polytope $P^n$ that is a candidate for the 
fundamental tile, let $W$ be the right-angled Coxeter group with one
generator for each codimension-one face and one relation for each
codimension-two face.  Let $\tilde{M}^n$ be the result of applying the
reflection group construction to $P^n$ and $W$, and let $\Gamma$ be a
torsion-free, normal subgroup of $W$.  Then we get examples of the
above type with $M^n=\tilde{M}^n/\Gamma$, $G=W/\Gamma$, and $A=W$.
Again, we note that these reflection type examples are {\em not} the
ones of primary interest in this paper. 

Our primary examples are manifolds that are constructed by blowing up
certain subspaces of projective hyperplane arrangements in
$\bbR\bbP^n$.  The theory of such blow-ups was developed in
\cite{DJS}.  Given a hyperplane arrangement in $\bbR^{n+1}$, there is
an associated $(n+1)$-dimensional convex polytope $Z$ called a
``zonotope''.  An equivalent formulation of the blowing-up procedure
is described in \cite{DJS}: one ``blows-up'' certain cells of
$\partial Z/a$ (where $a$ denotes the antipodal map).  In this
generality, the resulting cubical cell complex might not admit a
suitable symmetry group $G$ satisfying property (1).  The condition
needed is that the original zonotope $Z$ admit a group of symmetries
that is simply transitive on the vertex set of $Z$.  The most obvious
zonotopes with this property are the so-called ``Coxeter cells''.
So, this paper is a continuation and specialization of
\cite{DJS} to the case of hyperplane arrangements associated to finite
reflection groups.  (N.B. a
{\em Coxeter cell} is a zonotope corresponding to a hyperplane
arrangement associated to a finite reflection group $W$ on
$\bbR^{n+1}$. A {\em Coxeter cell complex} is a regular cell complex in
which each cell is isomorphic to a Coxeter cell.  For example, since an
$(n+1)$-cube is the Coxeter cell associated to $(\bbZ_2)^{n+1}$, any cubical
complex is a Coxeter cell complex.)

In \cite{DJS} we also discussed a generalization of the blow-up
procedure to zonotopal cell complexes.  Again, in order for property
(1) to hold we need to require that the zonotopal cell complex admit a
group of automorphisms that acts simply transitively on its vertex
set.  Examples of zonotopal cell complexes with this property are
provided by Coxeter groups.  Associated to any Coxeter system $(W,S)$,
there is a Coxeter cell complex $\Sigma(W,S)$ such that $W$ acts
simply transitively on its vertex set.  (Here $W$ might be infinite.)
Thus, we also want to apply our blowing up procedures to the complexes
$\Sigma(W,S)$.

As data for such a blowing up procedure it is necessary to
specify the set of cells which are to be blown up.  There are two extreme
cases, the ``minimal blow-up'' and the ``maximal blow-up''.  In the case of a
minimal blow-up, this set of cells is the collection of all cells 
that cannot be decomposed as a nontrivial product.  In the case of a
maximal blow-up, it is the set of all cells.

Next we describe the motivating example for this paper (which was also
one of the motivating examples for \cite{DJS}).  Consider the action
of the symmetric group $S_{n+2}$ as a reflection group on $\bbR^{n+1}$.
The associated hyperplane arrangement is called the ``braid
arrangement''.  Let $M^n$ denote the minimal blow-up (as in \cite{DP}
or \cite{DJS}) of the corresponding arrangement in $\bbR\bbP^n$.  The
interesting feature of $M^n$ lies in the result of Kapranov \cite{Ka1,Ka2}
that $M^n$ can be identified with $\overline{\cM}_{0,n+3}(\bbR)$, the
real points of the Grothendieck-Knudsen moduli space of stable
$(n+3)$-pointed curves of genus $0$ (which, in turn, coincides with
the Chow quotient $(\bbR\bbP^1)^{n+3}//PGL(2,\bbR)$).  Kapranov also
showed that each tile of the dual cellulation of $M^n$ was a copy of
Stasheff's polytope, the $n$-dimensional associahedron, $K^n$.

As explained in \cite{Lee} or in
Section~\ref{s:associahedral-tilings}, below, the
set of codimension-one faces of $K^n$ can be identified with the set
of proper subintervals of $[1,n+1]$ with integer endpoints.  Moreover,
given two such subintervals $T$ and $T'$, the corresponding faces
intersect if and only if either (i) the distance between $T$ and $T'$ (as
subsets of $[1,n+1]$) is at least $2$ or (ii) $T'\subset T$ (or
$T\subset T'$).

In the case at hand, where $M^n=\overline{\cM}_{0,n+3}(\bbR)$, the
symmetry group $G$ is $S_{n+2}$.  The group $A$ has one involutory
generator $\alpha_T$ for each proper subinterval $T$ of $[1,n+1]$.  The
codimension-two faces of $K^n$ impose additional relations of two
types: (i) if the distance between $T$ and $T'$ is at least $2$ then
$(\alpha_{T}\alpha_{T'})^2=1$ and (ii) if $T'\subset T$, then
$\alpha_T\alpha_{T'}\alpha_T=\alpha_{T''}$ (where $T''$ denotes the
image of $T'$ under the order-reversing involution of $T$).  The
epimorphism $A\rightarrow S_{n+2}$ sends $\alpha_T$ to the
order-reversing involution in the subgroup of $S_{n+2}$ corresponding
to $T$.  Looking at relation (ii), it is clear that if the interval
$T$ is not a single point, then $\alpha_T$ will not act as a
reflection on $\tilde{M}^n$.  We call it a ``mock reflection'' and $A$
a ``mock reflection group''.

Similarly, given any finite Coxeter group $W$, one can take the minimal
blow-up of the associated projective hyperplane arrangement to obtain
a manifold $M^n$ with a cubical cell structure.  When the Coxeter
diagram of $W$ is an interval, the tiles will again be associahedra.

Other examples arise by taking the maximal blow-up of
an arrangement associated to a finite reflection group $W$.  In any such
example each tile is a permutohedron.  (In the case where $W=(\bbZ_2)^{n+1}$
these examples occur in nature as real toric varieties associated to
flag manifolds.)

In Sections~\ref{s:permuto} and \ref{s:associahedral-tilings} we prove some
classification results for the universal covers of the permutohedral and
associahedral tilings which arise from blow-ups.
In Section~\ref{s:permuto}, we show that the universal covers of all
such permutohedral tilings yield the same tiling of $\bbR^n$;
moreover, the various symmetry groups $A$ that arise in this fashion
are commensurable with each other (and with the right-angled
reflection group associated to the permutohedron).  By way of
contrast, in Section~\ref{s:associahedral-tilings}, we show that the various
associahedral tilings of $\bbR^n$ tend not to be isomorphic with each
other.  The reason for this dichotomy lies in the fact that the
associahedron is much less symmetric than is the permutohedron.
It turns out, however, that in dimensions $\leq 3$ all of the
symmetry groups arising from these permutohedral and associahedral
tilings are quasi-isometric to each other (Theorem~\ref{thm:2d-assoc-class} and
Theorem~\ref{thm:q-iso}).

Section~\ref{s:tilings}, the longest section of the paper,
concerns the general theory of tilings.  The results in this section are
of a somewhat different nature than in the rest of the paper.
We develop the theory in a context which is considerably
more general than is indicated by the above examples.  All of our
previous requirements are either weakened or dropped as explained
below in statements (a) through (e).
\begin{enumerate}
\item[(a)] The space is not required to be a manifold.
\item[(b)] The cell structure on the space need not be cubical;
however, each cell is required to be a Coxeter cell.  (This is to
accomodate the ``partial blow-ups'' of \cite{DJS} and also to
include arbitrary reflection type tilings in our general theory).
\item[(c)] In view of (a), there may no longer be a well-defined dual
cell structure; however, there are still ``dual cones'' and ``tiles''
(that is, cones dual to vertices).
\item[(d)] The requirement that there exist a group $G$ that acts
simply transitively on the vertex set is replaced
by the requirement that the cell complex
$\ccc$ admit a ``framing'' (which amounts to specifying isomorphisms
between the links of any two vertices of $\ccc$).
\item[(e)] The requirement of nonpositive curvature is dropped.
\end{enumerate}
Thus,  in Section~\ref{s:tilings} we shall largely abandon the notion
of a fundamental tile (since, in view of (c), it need not be a cell).
By a ``tiling'' we will simply mean a Coxeter cell
complex $\ccc$ in which the links of any two vertices are isomorphic. 
The tiling
is ``symmetric'' if $\ccc$ admits a group action which is simply transitive on
its vertex set.  If $\ccc$ is symmetric and simply connected, then one can
read off a presentation for its symmetry group $A$
(cf. Section~\ref{ss:homogeneous-framings}) as before. 

A key ingredient in our analysis of framed tilings and their symmetry groups
is the notion of a ``gluing isomorphism.''   This is is an isomorphism
between two ``codimension-one faces'' of a fundamental tile.  It determines how
two adjacent tiles are glued together.  In practice (e.g. when the
symmetry group is generated by involutions), these gluing isomorphisms
will always be involutions.  For simplicity, let us assume this.  For example,
in a reflection
type tiling, each gluing involution is the identity map.  For blow-ups
of $\Sigma(W,S)$, the gluing involutions are determined by the
elements of longest length in various finite special subgroups of $W$.  In
\ref{ss:glu-iso} and \ref{ss:iso-reflection}, we give
necessary and sufficient conditions on the sets of gluing involutions
for the universal covers of two tilings to be isomorphic.  (This result is
then used in Sections~\ref{s:permuto} and \ref{s:associahedral-tilings} to
classify certain permutohedral and associahedral tilings.)
In Theorem~\ref{thm:gluings-to-group}, we describe necessary and
sufficient conditions under which a given fundamental tile and set
of gluing involutions can be realized by a symmetric tiling.

In Sections~\ref{s:tilings} and \ref{s:linearity} we also prove, under very mild hypotheses,
some general results about the
symmetry group $A$ of a symmetric and simply connected $\ccc$.
Since these hypotheses are satisfied in our examples,
these results apply to the
universal cover of a blow-up.  First of all, even without the
nonpositive curvature requirement, $\ccc$ can always be ``completed''
to a CAT($0$)-complex
$\hat{\ccc}$ by adding a finite number of $A$-orbits of cells.  (This
is proved in Sections~\ref{ss:hat-X} and \ref{ss:cat-0}.)  Thus, each such
$A$ is a ``CAT($0$) group''
(cf. Theorem~\ref{thm:CAT(0)-groups-from-framings}).  In particular, if $\pi$
is any torsion-free subgroup of finite index in $A$, then $\hat{\ccc}/\pi$
is a finite $K(\pi ,1)$-complex.   Secondly, $A$ has a
linear representation analogous to the canonical represention for a
Coxeter group (Section~\ref{s:linearity}).  Whether this
representation is always faithful, however, remains an interesting
open question. 

\section{Some cell complexes associated to Coxeter groups}\label{s:Cox}

\subsection{Coxeter systems}

We first recall some standard facts about Coxeter groups; we refer the
reader to \cite{Bo}, \cite{Bro}, or \cite{H} for details.  

\begin{definition}\label{def:coxeter}
Let $S$ be a finite set.  A {\em Coxeter matrix} on $S$ is a symmetric
$S\times S$ matrix $M=(m(s,s'))$ with entries in $\bbN\cup\{\infty\}$
such that $m(s,s')=1$ if $s=s'$ and $m(s,s')\geq 2$ if $s\neq s'$.
Given a Coxeter matrix $M$, the corresponding {\em Coxeter group} is the
group $W$ defined by the presentation  
\[W=\langle S \mid (ss')^{m(s,s')}=1\; \mbox{for all}\; s,s'\in S\rangle.\]
The pair $(W,S)$ is a {\em Coxeter system}.  
\end{definition}

A Coxeter system $(W,S)$ determines a {\em length function}
$l:W\rightarrow\bbZ_{\geq 0}$.  Given $w\in W$, the length $l(w)$ is defined to
be the minimal $n$ such that $w=s_1s_2\cdots s_n$ and $s_i\in S$.  

To any Coxeter system $(W,S)$ we associate a labeled graph
$\Gamma(W,S)$, called the {\em Coxeter diagram}, as follows.  The vertex
set is $S$ and two vertices $s,s'$ determine an edge if and only if
$m(s,s')>2$.  In this case, the edge joining $s$ and $s'$ is labeled
$m(s,s')$.     

Any subset $T\subset S$ generates a subgroup $W_T\subset W$, 
which is itself a Coxeter group with Coxeter system $(W_T,T)$.  
$W_T$ is called a {\em special subgroup} of $W$.  The
Coxeter diagram $\Gamma(W_T,T)$ is the induced labeled subgraph of
$\Gamma(W,S)$ with vertex set $T$.  A subset $T\subset S$ is {\em
spherical} if $W_T$ is finite.  Any finite special subgroup $W_T$ has
a unique element of longest length, denoted by $w_T$.  Moreover, $w_T$
is an involution and $w_TTw_T=T$ (Exercise 22, page 43 in \cite{Bo}).   

If $W_1$ and $W_2$ are Coxeter groups, then so is $W=W_1\times W_2$.
Let $(W_1,S_1)$ and $(W_2,S_2)$ be Coxeter systems for $W_1$ and
$W_2$, respectively, and let $S=S_1\times\{1\}\cup\{1\}\times S_2$.
Then $(W,S)$ is a Coxeter system for $W$, and the Coxeter diagram
$\Gamma(W,S)$ is the disjoint union of the labeled graphs $\Gamma(W_1,S_1)$
and $\Gamma(W_2,S_2)$.

\subsection{Coxeter cells}

Let $W$ be a finite Coxeter group.  Then $W$ can be represented as a
group generated by orthogonal reflections on a finite dimensional
Euclidean space $V$.  The reflection hyperplanes of this representation
separate $V$ into simplicial cones, called {\em  chambers}, and $W$
acts transitively on the set of chambers.  In fact, the representation
of $W$ can be chosen such that for any Coxeter system $(W,S)$, the
generators in $S$ correspond to the reflections through the supporting
hyperplanes of a fixed chamber $C$.  We call this fixed chamber the {\em
fundamental chamber}, and we call the reflections in $S$ the {\em
simple reflections}.  

\begin{definition}\label{def:cox-cell}
Let $x$ be a point in $C$ that is unit distance from each of the
supporting hyperplanes, and let $\cc=\cc(W,S)$ be the convex hull of
the orbit $W x$. The polytope $\cc$ is called a {\em
(normalized) Coxeter cell of type $W$}.  The intersection of $\cc$
with $C$ (or with any translate of $C$ by an element of $W$) is called
a {\em Coxeter block of type $W$}.
\end{definition}

\begin{figure}[ht]
\begin{center}
\psfrag{1}{$1$}
\psfrag{D}{$B$}
\psfrag{Coxeter Diagram}{\hspace{-.5in}$\scriptsize\begin{array}{cc}\mbox{Coxeter cell and
Coxeter block}\\ \mbox{associated to the diagram:}\end{array}$}  
\psfrag{3}{\scriptsize\hspace{-.05in}$3$}
\includegraphics[scale = .7]{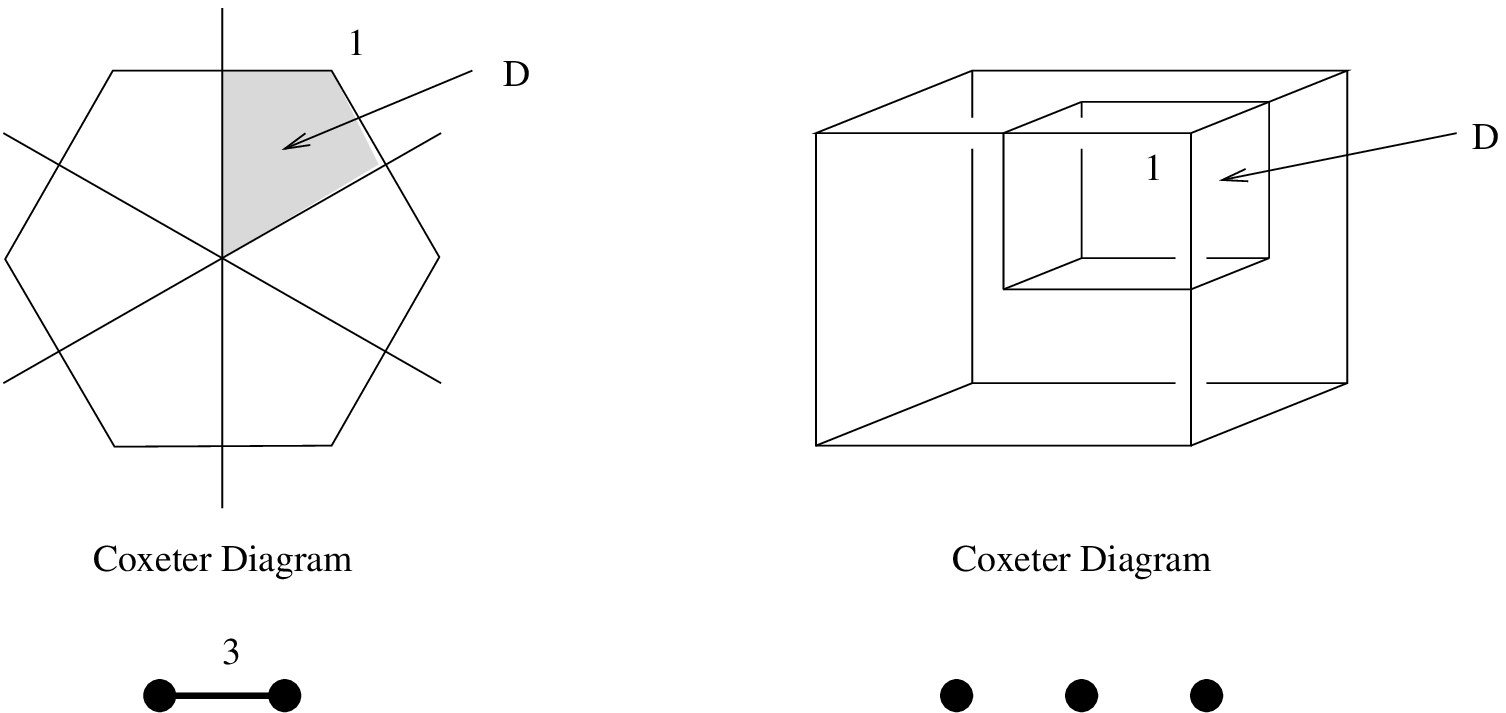}
\caption{\label{fig:cox1}}
\end{center}
\end{figure}

The group $W$ acts isometrically on $\cc$, and the Coxeter block
$B=\cc\cap C$ is a fundamental domain (Figure~\ref{fig:cox1}).  $B$ is
combinatorially 
equivalent to a cube of dimension $\card S$.  Since $W$ acts freely
and transitively on the vertices of $\cc$, we can identify the
vertices of $\cc$ with the elements of $W$ (once we identify $x$
with $1$).  Each vertex is contained in a unique Coxeter block of type
$W$. 

The Coxeter block $B$ has two types of codimension-one faces.  One
type is an intersection of $B$ with a codimension-one face of $C$.
These are the {\em mirrors} of $B$.  To descibe the others, we first
describe the faces of a Coxeter cell.  Let $W_T$ be a special subgroup
and let $\cc_T$ be a normalized Coxeter cell of type $W_T$.  Then the
inclusion $W_T\rightarrow W$ induces an isometry from $\cc_T$ onto a
face of $\cc$ (in fact, every face of $\cc$ is of the form $w\cc_T$ for some
$T\subset S$ and $w\in W$).  The remaining codimension-one faces of the 
Coxeter block $B$ can now be identified with the Coxeter blocks
associated to the codimension-one faces $\cc_T\subset \cc$ (i.e., where
$T$ has cardinality one less than $S$).  

If $\cc_1$ and $\cc_2$ are (normalized) Coxeter cells of types
$W_1$ and $W_2$, respectively, then the product
$\cc_1\times \cc_2$ is a (normalized) Coxeter cell of type
$W_1\times W_2$.  In particular, the $n$-cube $[-1,1]^n$ is a Coxeter
cell; the Coxeter block containing the vertex $(1,1,\ldots,1)$ is
$[0,1]^n$ (Figure~\ref{fig:cox1}). 

\subsection{Coxeter cell complexes}\label{ss:cox-cell-complex}

A locally finite, regular cell complex $\ccc$ is a {\em Coxeter cell
complex} if all of its cells are Coxeter cells.  Since any
combinatorial isomorphism between two (normalized) Coxeter cells $\cc$
and $\cc'$ is induced by an isometry, any Coxeter cell complex $\ccc$ has a
canonical piecewise Euclidean metric, and any isomorphism of Coxeter
cell complexes is induced by a unique isometry.    

Let $\ccc$ be a Coxeter cell complex $\ccc$, and let $\ccc^{(i)}$ denote the
set of $i$-cells in $\ccc$.  Let $\cp\in \ccc^{(0)}$ be a vertex of $\ccc$.
The {\em link of $\cp$}, denoted $\Lk(\cp,\ccc)$, is the
(piecewise spherical) simplicial complex consisting of all points in $\ccc$
that are unit distance from $\cp$.  The {\em closed star of $\cp$},
denoted $\clstar(\cp,\ccc)$, is the subcomplex of $\ccc$ consisting of
all cells that contain $\cp$ and all of their faces.   

\begin{definition}\label{def:tiling}
A Coxeter cell complex $\ccc$ is a {\em tiling} if for any
two vertices $\cp$ and $\cp'$, there exists a combinatorial isomorphism 
$\clstar(\cp,\ccc)\rightarrow\clstar(\cp',\ccc)$ taking $\cp$ to $\cp'$. 
\end{definition}

Let $(W,S)$ be a (possibly infinite) Coxeter system and let $\cS$ be 
the set of spherical subsets of $S$.  (It is clear that
$\cS_{>\emptyset}$ is an abstract simplicial complex with vertex set
$S$.)  We can then define a Coxeter cell complex $\Sigma=\Sigma(W,S)$
that generalizes the Coxeter cell of a finite Coxeter group.  The
vertex set of $\Sigma$ is $W$.  We take a   
Coxeter cell of type $W_T$ for each coset $wW_T$ where $T\in\cS$,
and identify its vertices with the elements of $wW_T$.  We then
identify two faces of two Coxeter cells if they have the same vertex
set.  If $Z$ is the cell in $\Sigma$ corresponding to the coset
$wW_T$, we define its {\em type} to be the subset $T$.  (The type of a
cell is well-defined since $wW_T=w'W_{T'}$ implies $T=T'$.)

It is clear that $\Sigma$ is a Coxeter cell complex and that $W$ acts
via combinatorial automorphisms.  Since $W$ acts simply transitively
on the vertices, $\Sigma$ is a tiling.  The cells of $\Sigma$ that
contain the vertex $1$ are in bijection with the set
$\cS_{>\emptyset}$; thus, $\Lk(1,\Sigma)$ can be identified
with $\cS_{>\emptyset}$.  More generally, we have the following:  

\begin{definition} \label{def:reflection-tiling}
Let $L$ be a subcomplex of $\cS_{>\emptyset}$ with the same
$1$-skeleton.  Then the {\em reflection tiling} of type $(W,S,L)$,
denoted $\Sigma(W,S,L)$, is the subcomplex of $\Sigma(W,S)$ consisting
of cells corresponding to the cosets $wW_T$ where $T$
is either the empty set or the vertex set of a simplex in $L$.  A
reflection tiling is {\em complete} if $L=\cS_{>\emptyset}$. 
\end{definition} 

\begin{example}
If $L$ is the $0$-skeleton of $\cS_{>\emptyset}$, then $\Sigma(W,S,L)$
is the Cayley graph of $(W,S)$.  Similarly, if $L$ is the $1$-skeleton
of $\cS_{>\emptyset}$, then $\Sigma(W,S,L)$ is the Cayley $2$-complex
associated to the standard presentation of $W$.
\end{example}

\begin{example}
Suppose $W$ is finite.  Then $\Sigma(W,S)$ is the zonotope $Z(W,S)$.
If $L$ is the boundary of the simplex $\cS_{>\emptyset}$, then
$\Sigma(W,S,L)=\partial Z$.
\end{example}

\subsection{Coxeter tiles and the local geometry of Coxeter cell
complexes}\label{ss:local-geometry}

Let $\ccc$ be a Coxeter cell complex.  For any vertex $\cp\in \ccc^{(0)}$, the
{\em dual cone} at $\cp$, denoted $D(\cp,\ccc)$ is the union of all Coxeter
blocks in $\ccc$ that contain the vertex $\cp$.  It can be identified with
a subcomplex of the barycentric subdivision of $\ccc$ and therefore, has
a natural piecewise Euclidean metric.

\begin{definition}\label{def:coxeter-tile}
Let $\Sigma$ be the reflection tiling of type $(W,S,L)$.  Then the
{\em Coxeter tile} of type $(W,S,L)$, denoted $D(W,S,L)$, is the
dual cone $D(1,\Sigma)$.  If $L=\cS_{>\emptyset}$, we shall denote the
corresponding Coxeter tile $D(W,S)$.  An {\em isomorphism}
$D(W,S,L)\rightarrow D(W',S',L')$ of Coxeter tiles is a bijection
$S\rightarrow S'$ that (1) extends to a simplicial isomorphism
$L\rightarrow L'$, and (2) extends to a group isomorphism
$W\rightarrow W'$.   
\end{definition}    

It is clear that the geometric realization of an isomorphism of
Coxeter tiles is an isometry.

\begin{example}
Let $(W,S)$ be the (infinite) Coxeter group with generating set
$S=\{a,b,c\}$ and Coxeter diagram as in Figure~\ref{fig:cox2}.  Then
the complete reflection tiling $\Sigma(W,S)$ is an infinite
$2$-dimensional cell complex whose $2$-cells are regular hexagons and
squares.  The Coxeter tile $D(W,S)$ is a union of two Coxeter blocks,
one from the hexagon and one from the square.  See Figure~\ref{fig:cox2}.  

\begin{figure}[ht]
\begin{center}
\psfrag{D}{$D$}
\psfrag{CT}{\hspace{-1.2in}$\scriptsize\begin{array}{cc}
\mbox{Complete reflection tiling (and Coxeter tile)}
\\ \mbox{associated to the (infinite) Coxeter
group:}\end{array}$}  
\psfrag{IT}{\hspace{-1.1in}$\scriptsize\begin{array}{cc}
\mbox{Reflection tiling (and Coxeter tile)}
\\ \mbox{associated to the subcomplex}\\
\mbox{$L=\{\{a\},\{b\},\{c\}\}$}\end{array}$}   
\psfrag{3}{\scriptsize\hspace{-.01in}$3$}
\psfrag{I}{\scriptsize\hspace{-.03in}$\infty$}
\psfrag{a}{\scriptsize\hspace{-.03in}$a$}
\psfrag{b}{\scriptsize\hspace{-.03in}$b$}
\psfrag{c}{\scriptsize\hspace{-.03in}$c$}
\includegraphics[scale = .7]{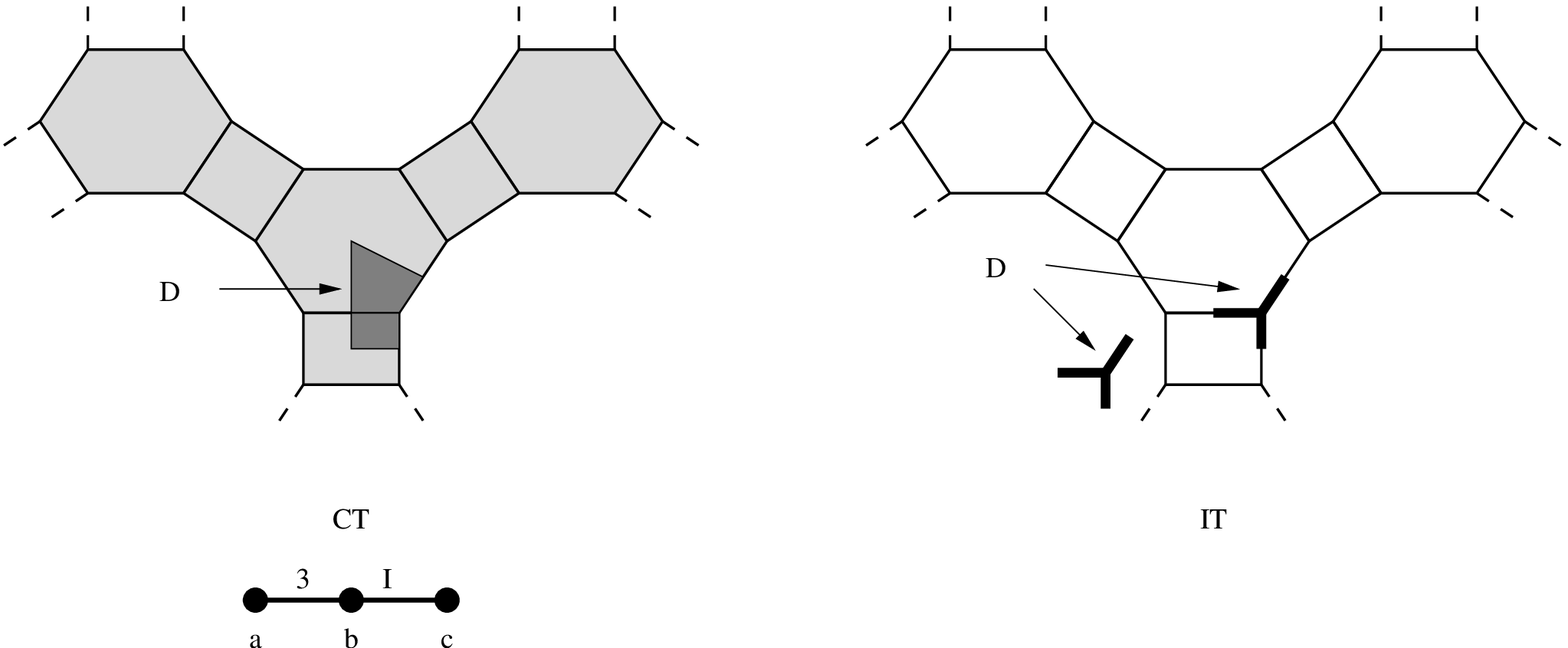}
\caption{\label{fig:cox2}}
\end{center}
\end{figure}

If we take $L$ to be the $0$-skeleton of $\cS_{>\emptyset}$ obtained
by removing the spherical subsets $\{a,b\}$ and $\{a,c\}$, we obtain
the second reflection tiling $\Sigma(W,S,L)$ shown in Figure~\ref{fig:cox2}.
The associated Coxeter tile $D(W,S,L)$ is the cone on three vertices.
(Note that in this example, $\Sigma(W,S,L)$ is the Cayley 
graph of $(W,S)$, and $\Sigma(W,S)$ is the Cayley $2$-complex.)
\end{example}

\begin{example}
Suppose $W$ is finite.  Then $\mathcal{S}_{>\emptyset}$ is a simplex
and $D(W,S)$ is the Coxeter block of type $(W,S)$, as in
Definition~\ref{def:cox-cell}.  If $L=\cS_{>\emptyset}-\{S\}$ is 
the boundary of the simplex $\cS_{>\emptyset}$, then $D(W,S,L)$ is 
called the {\em fundamental simplex}.  It is the intersection of
$\partial \cc$ with the fundamental chamber $C$.  The image of
$\partial Z \cap C$ in projective space is a simplex.  So, we think of
the cell complex $\Delta = D(W,S,L)$ as being a subdivision of the
simplex into Coxeter blocks.\end{example} 

\begin{remark}
The Coxeter tile $D(W,S,L)$ is a fundamental domain for the
$W$-action on the reflection tiling $\Sigma(W,S,L)$.  
\end{remark}

If $\ccc$ is any Coxeter cell complex and $\cp$ is a vertex of $\ccc$, then
the dual cone $D(\cp,\ccc)$ is a Coxeter tile.  To determine its type, we
let $L_{\cp}=\Lk(\cp,\ccc)$, we let $V_{\cp}$ be the set
of vertices in $L_{\cp}$ (equivalently, the set of edges containing
$\cp$), and we define a Coxeter matrix $M_{\cp}$ on $V_{\cp}$ as follows:
\[ m_{\cp}(v,v')=\left\{\begin{array}{ll}
1 & \mbox{if $v=v'$}\\
m & \mbox{if $v$ and $v'$ correspond to distinct edges of a $2m$-gon
in $\ccc$}\\
\infty & \mbox{otherwise}.\end{array}\right.\]
Letting $(W_{\cp},V_{\cp})$ be the corresponding Coxeter system, it is
clear that $L_{\cp}$ consists of spherical subsets of
$V_{\cp}$.  It follows that $D(\cp,\ccc)$ can be identified
with the Coxeter tile of type $(W_{\cp},V_{\cp},L_{\cp})$.   

\begin{remark} 
If $\ccc$ is a tiling, then there exists a Coxeter tile
$D=D(W,S,L)$ and, for each $\cp\in \ccc^{(0)}$, an isomorphism
$\phi_{\cp}:D\rightarrow D(\cp,\ccc)$.  Since $D$ may have nontrivial
automorphisms, the isomorphism $\phi_{\cp}$ need not be unique. 
\end{remark} 

\section{The blow-up $\Sigma_{\#}$ of a reflection tiling $\Sigma$}\label{s:blow-ups}

\subsection{The blow-up of a Coxeter cell}\label{ss:blow-up-cell}
Let $W$ be a finite Coxeter group, $(W,S)$ be a Coxeter system, and
let $\cc$ ($=\cc(W,S)$) be the corresponding Coxeter cell.  Then the
boundary $\Bd \cc$ is a Coxeter cell complex homeomorphic to
$S^{n-1}$.  It is the reflection tiling $\Sigma(W,S,L)$ where $L$ is
the boundary complex of the simplex $\cS_{>\emptyset}$.  A fundamental
domain for the $W$-action on $\Bd \cc$ is the fundamental simplex
$\Delta=\Bd \cc\cap C$ ($=D(W,S,L)$).  A {\em mirror} of 
$\Delta$ is its intersection with a codimension-one face of $C$.  The
Coxeter block $B=\cc\cap C$ is homeomorphic to the cone on $\Delta$.

The antipodal map $a:\cc\rightarrow \cc$ defined
by $x\mapsto -x$, restricts to an isometric involution on
$\partial \cc$.  This involution on $\Bd \cc$ freely permutes the 
cells; hence the quotient $\bbP(\cc)=\Bd \cc/a$ is a Coxeter cell
complex, homeomorphic to $\bbR\bbP^{n-1}$.

\begin{definition}\label{def:blow-up}
The {\em blow-up of $ \cc$ at its center}, denoted
$ \cc_{\#}$, is the quotient of $\Bd \cc\times[-1,1]$ by
the involution $\hat{a}$ defined by $\hat{a}(x,t)=(-x,-t)$, i.e.,  
\[ \cc_{\#}=\Bd \cc\times_{\hat{a}}[-1,1].\]
\end{definition}

The blow-up $ \cc_{\#}$ is an interval bundle over $\bbP( \cc)$.
Let $\pi: \cc_{\#}\rightarrow\bbP( \cc)$ be the projection.  For
any face $F\subset\Bd \cc$, $-F\cap F=\emptyset$; hence, $\pi(F)$ is
embedded in $\bbP( \cc)$.  Moreover, $\pi^{-1}(\pi(F))$ can be
identified with $F\times[-1,1]$ (where $F\times\{1\}$ corresponds to
$F$ and $F\times\{-1\}$ to $-F$).  Since $F$ is a Coxeter cell, so is
$F\times[-1,1]$.  Thus, $ \cc_{\#}$ is naturally a Coxeter cell
complex.  Since the antipodal map $a$ commutes with the $W$-action,
there is an induced $W$-action on $ \cc_{\#}$.

Similarly, $\pi(\Delta)$ is embedded in $\bbP( \cc)$ and
$\pi^{-1}(\pi(\Delta))\cong\Delta\times[-1,1]$.  Since $w_S$ maps
$\Delta$ to $-\Delta$, the stabilizer of $\Delta\times[-1,1]$ is the
cyclic group of order two generated by $w_S$.  Thus,
$\Delta\times[0,1]$ is a fundamental domain for $W$ on $ \cc_{\#}$.
Put 
\[B_{\#}=\Delta\times[0,1].\]
We think of $B_{\#}$ as being obtained from $B$ by ``truncating'' its
vertex (the cone point) and introducing a new codimension-one face
corresponding to $\Delta\times\{0\}$.  A {\em mirror} of $B_{\#}$ is
either a codimension-one face of the form $\Delta_{\{s\}}\times[0,1]$,
where $\Delta_{\{s\}}$ is a mirror of $\Delta$, or the codimension-one
face $\Delta\times\{0\}$.

The space $ \cc_{\#}$ is tiled by the translates of $B_{\#}$ by
elements of $W$. The adjacent tile to $B_{\#}$ across the mirror
$\Delta\times\{0\}$ is $w_SB_{\#}$.  If we identify $w_SB_{\#}$ with
$\Delta\times[0,1]$ via $w_S$, then $B_{\#}\cup w_SB_{\#}$
($\cong\Delta\times[-1,1]$) is homeomorphic to two copies of
$\Delta\times[0,1]$ glued along $\Delta\times\{0\}$.  The gluing map
$j_S:\Delta\rightarrow\Delta$ is given by $j_S=a\circ w_S$.

Let $1$ denote the vertex of $ \cc$ that is contained in the
interior of $C$.  We note that $\Delta$ is the union of all Coxeter
blocks in $\Bd \cc$ which contain the vertex $1$ and that $B$ is the
union of all Coxeter blocks in $ \cc$ that contain $1$.  Similarly,
$B_{\#}$ is the union of all Coxeter blocks in $ \cc_{\#}$ that
contain $1$.

\begin{remark}\label{rem:a-is-longest-word}
Suppose $(W,S)$ is a finite irreducible Coxeter group.  The element of
longest length, $w_S$, is equal to the antipodal map, $a$, in the
following cases (\cite{Bo}, Appendix I-IX, pp. 250-275): $\bA_1$ (the
cyclic group of order $2$), $\bI_2(p)$ with $p$ even (the dihedral group
of order $2p$), $\bB_n$ (the hyperoctahedral groups), $\bD_n$ with $n$
even, $\bH_3$, $\bH_4$, $\bF_4$, $\bE_7$, and $\bE_8$.  Hence, in all these
cases, $j_S$ is the identity map.

In the remaining cases $w_S$ is not the antipodal map, and conjugating
by it induces a nontrivial diagram automorphism of $\Gamma(W,S)$.
These cases are: $\bA_n$ with $n>1$, $\bD_n$ with $n$ odd, $\bI_2(p)$ with
$p$ odd, and $\bE_6$.  In each of these cases, the Coxeter diagram
admits a unique non-trivial automorphism which, in fact, coincides
with the diagram automorphism induced by conjugation by $w_S$. 

If $W$ is finite and reducible, then $w_S$ is the product of the
elements of longest length in each factor and, therefore, is the
antipodal map if and only if it is antipodal in each factor.
\end{remark}

\subsection{The blow-up of a reflection tiling
$\Sigma(W,S,L)$}\label{ss:blow-up-functor}  

Our goal in this section is to describe a functorial
generalization of the blow-up of a Coxeter cell that
(1) allows for iterated blow-ups of faces, (2) makes sense more
generally for the complexes $\Sigma(W,S,L)$, and (3) preserves the
$W$-action.  Given a complex $\Sigma(W,S,L)$ and a suitable
collection of cells to be blown-up, this functor will produce for
every subcomplex $K\subset \Sigma$ a cell complex $K_{\#}$ called the
``blow-up of $K$''.  Our primary interest in this paper is the
topology of the blow-up $\Sigma_{\#}=\Sigma(W,S,L)_{\#}$.  

\begin{definition}
Suppose $L$ is a simplicial complex with vertex set $V$.  We let
$\cP(L)$ denote the poset consisting of those subsets $T$ of $V$ such
that either $T=\emptyset$ or $T$ is the vertex set of a simplex in
$L$.  Thus, we can identify $L$ with $\cP(L)_{>\emptyset}$. 
\end{definition}

Let $\Sigma$ be the reflection tiling $\Sigma(W,S,L)$, and let $\cP$
be the poset $\cP(L)$.  Then the cells of $\Sigma$ are indexed by the
set $W\cP=\{wW_T\;|\; T\in\cP\}$ (see
Definition~\ref{def:reflection-tiling}).  Any collection of cells  
defining one of our blow-ups $\Sigma_{\#}$ will satisfy the following
condition:  if a face $ \cc'$ of a cell $ \cc$ is blown-up, then
either   
\begin{enumerate}
\item[(1)] $ \cc$ is blown up also, or 
\item[(2)] $ \cc$ is of the form $ \cc'\times \cc''$ and the
blow-up functor is applied to the two factors independently.  
\end{enumerate}
In addition, if the $W$-action on $\Sigma$ is to induce an action on the
blow-up $\Sigma_{\#}$, the collection of cells we blow-up should be
$W$-invariant.  Since the cells of $\Sigma$ are indexed by the poset
$W\cP$, the blow-up collection will be indexed by a subset of the
form $W\cR$ for some subset $\cR\subset \cP$.  With this in mind, we
let $\cR$ be any subset of $\cP$, and define a category whose
objects are all subcomplexes of $\Sigma$, but whose morphisms depend on
$\cR$.  

A {\em subcomplex} of $\Sigma$ is a Coxeter cell complex $K$ together with an
injective cellular map $K\rightarrow \Sigma$.  Via this map, we will often
identify $K$ with its image in $\Sigma$.  In particular, any cell $ \cc$
of a subcomplex $K$ has a well-defined type $T$ where $T\in\cP$
(cf. section~\ref{ss:cox-cell-complex}).  An {\em $\cR$-morphism}
$f:K_1\rightarrow K_2$ between two subcomplexes of $\Sigma$ is an
injective cellular map satisfying the condition: for every cell $
\cc\subset K_1$, the type of $\cc$ is in $\cR$ if and only if the type
of $f( \cc)$ is in $\cR$. 

There are primarily two kinds of $\cR$-morphisms relevant to our
construction.  Since the $W$-action on $\Sigma$ preserves the cells of a
given type, the automorphism $w:\Sigma\rightarrow \Sigma$ is an $\cR$-morphism
for every $w\in W$.  Other $\cR$-morphisms are provided by antipodal
maps on cells in $\Sigma$ (though in general only some of these are
$\cR$-morphisms).  To make sense of (1) above, i.e., to iterate the
blow-up procedure, we need the antipodal map on $ \cc$ to be a
morphism in our category.  

\begin{lemma}\label{lem:antipodal}
Given $T\in\cP$, let $j_T:\cP_{\leq T}\rightarrow\cP_{\leq
T}$ be the involution defined by $T'\mapsto w_T T'w_T^{-1}$ where $w_T$ is
the longest element in $W_T$.  Let $ \cc_T\subset \Sigma$ be the cell
with vertices $W_T$, and let $a_T:\Bd \cc_T\rightarrow\Bd \cc_T$
be the antipodal map $x\mapsto -x$.  Then $a_T$ is an $\cR$-morphism
if and only if    
\[j_T(\cR_{\leq T})=\cR_{\leq T}.\]
\end{lemma}

\begin{proof}{}
The involution $w_T$ is an $\cR$-morphism and the composition $a\circ
w_T$ maps the face $Z_{T'}$ onto the face $Z_{j_TT'}$.   
\end{proof}

To make sense of (2), we recall the basic facts about product
decompositions of Coxeter cells.  Let $\Gamma$ denote the Coxeter
diagram $\Gamma(W,S)$, and for every $T\subset S$, let $\Gamma_T$ denote
the subdiagram $\Gamma(W_T,T)$.  Two spherical sets $T_1$ and $T_2$ in
$\cP$ are {\em completely disjoint} if they are disjoint and no edge
of $\Gamma_{T_1\cup T_2}$ connects a vertex of $\Gamma_{T_1}$ to a
vertex of $\Gamma_{T_2}$.  Given $T\in\cP$, a collection
$\{T_1,\ldots,T_k\}\subset\cP$ is called a {\em decomposition of
$T$} if $T=T_1\cup\cdots\cup T_k$ and the $T_i$'s are pairwise completely
disjoint. If $\{T_1,\ldots,T_k\}$ is a decomposition of $T$, then
there is a corresponding direct product decomposition of
groups   
\[W_T=W_{T_1}\times\cdots\times W_{T_k}\]
and a corresponding direct product decomposition of cells
\[ \cc_T= \cc_{T_1}\times\cdots\times \cc_{T_k}.\]
If our blow-up collection $\cR$ is to satisfy (1) and (2) above, then
whenever $ \cc_T$ is not blown-up and
$ \cc_{T_1},\ldots, \cc_{T_k}$ are the maximal faces of 
$ \cc_T$ that are blown up, there should be some
other subset $T_0\subset T$ and a decomposition
\[ \cc_T= \cc_{T_0}\times \cc_{T_1}\times\cdots\times \cc_{T_k}.\] 
We make this precise as follows.   

Given $T\in\cP$, we define $\cR T$, the {\em
$\cR$-maximals in $T$}, to be the set of maximal elements of $\cR_{\leq 
T}$, and we define $T_0$, the {\em $\cR$-fixed part of $T$}, to be the
complement of the union of the $\cR$-maximals in $T$.  To satisfy (1) and
(2), $\{T_0\}\cup\cR T$, must be a decomposition of $T$.  Thus, if
$\cR T=\{T_1,\ldots,T_k\}$, then the $T_i$'s must be pairwise
completely disjoint.  Combining this requirement with the condition of 
Lemma~\ref{lem:antipodal}, we make the following definition. 

\begin{definition}\label{def:admissible}
A collection $\cR\subset\cP$ is {\em admissible} if it satisfies the
following two conditions:
\begin{enumerate}
\item For every $T\in\cR$, $j_T(\cR_{\leq T})=\cR_{\leq T}$.
\item For every $T\in\cP$, the set $\{T_0\}\cup\cR T$ is a
decomposition of $T$.
\end{enumerate}
The collection $\cR$ is {\em fully admissible} if, in addition,
whenever $s,t\in S$ and $2<m(s,t)<\infty$, then $\{s,t\}\in\cR$.  
The collection $\{T_0\}\cup\cR T$ is called the {\em $\cR$-decomposition
of $T$}, and the elements $T_1,\ldots,T_k$ of $\cR T$ are the
{\em blow-up factors of the decomposition}.
\end{definition}

Let $\sK=\sK_{\cR}(W,S,L)$ be the 
category whose objects are the subcomplexes of $\Sigma$ and whose morphisms
are $\cR$-morphisms.  Let $\sC$ be the category whose objects are
Coxeter cell complexes and whose morphisms are isomorphisms onto
subcomplexes.   

\begin{proposition}\label{prop:blow-up}
Let $\cR$ be admissible.  Then there exists a unique functor
$\sK\rightarrow\sC$, denoted by $K\mapsto K_{\#}$ (and $f\mapsto
f_{\#}$), satisfying the two properties:
\begin{enumerate}
\item If $T\in\cR$, then $( \cc_T)_{\#}$ is the quotient of
$\left(\partial \cc_T\right)_{\#}\times[-1,1]$ by the involution
$\hat{a}:(x,t)\mapsto(a_{\#}(x),-t)$, where 
$a:\partial \cc_T\rightarrow\partial \cc_T$ is the antipodal 
map.  That is,  
\[( \cc_T)_{\#}=\left(\partial \cc_T\right)_{\#}\times_{\hat{a}}[-1,1].\]
\item If $T\in\cP$ has $\cR$-decomposition
$\{T_0,T_1,\ldots,T_k\}$, then 
\[\left( \cc_{T}\right)_{\#}=
 \cc_{T_0}\times\left( \cc_{T_1}\right)_{\#}  
\times\cdots\times\left( \cc_{T_k}\right)_{\#}.\]
\end{enumerate}
\end{proposition}

\begin{proof}{}
This functor is a special case of the blow-up of ``partially mined
zonotopal cell complexes'' as defined in \cite{DJS}.  If $\cR$ is
admissible, then the pair $(\Sigma,W\cR)$ is a partially mined zonotopal
cell complex (cf., \cite{DJS}), and for any subcomplex $K\subset \Sigma$
one obtains a partially mined subcomplex $(K,\cM)$.  Any $\cR$-morphism
$f:K_1\rightarrow K_2$ between subcomplexes of $\Sigma$ correspond to a
morphism $f:(K_1,\cM_1)\rightarrow(K_2,\cM_2)$ of partially mined
zonotopal cell complexes.  The functor $K\mapsto K_{\#}$ (and
$f\mapsto f_{\#}$) is then precisely the blow-up functor
$(K,\cM)\mapsto K_{\#\cM}$ (and $f\mapsto f_{\#}$) defined in
\cite{DJS}.  (If $\cR$ is fully admissible, then we can omit the term
``partially'' from the above discussion.)

Properties (1) and (2) are evident from the construction in
\cite{DJS}, and uniqueness follows from the fact that the functor is
defined inductively by these properties (and the $W$-action).
\end{proof} 
 
\begin{definition}
Let $\Sigma$ be a reflection tiling $\Sigma(W,S,L)$, and let $\cR$ be
an admissible subset of $\cP(L)$.  Then the Coxeter cell complex
$\Sigma_{\#}$ obtained by applying the blow-up functor to $\Sigma$
is called the {\em $\cR$-blow-up of $\Sigma$} (or just the {\em
blow-up of $\Sigma$} if there is no ambiguity).  
\end{definition}

\begin{remark}
If $\cR$ is fully admissible, then each Coxeter cell in $\Sigma_{\#}$ is a
cube.
\end{remark}

\begin{remark} 
If $T\in\cR$, then since $( \cc_T)_{\#}$ is an interval bundle over
$(\Bd \cc_T)_{\#}/a_{\#}$, the fundamental group of
$( \cc_T)_{\#}$ is isomorphic to that of
$(\Bd \cc_T)_{\#}/a_{\#}$. 
\end{remark}

\begin{remark}\label{rem:bu-compatibility}
Blow-up functors are compatible in the following sense.  Suppose
$\Sigma$ is the complex $\Sigma(W,S,L)$, and $\cR$ is an admissible subset
of $\cP(L)$.  Let $\Sigma'$ be a complex of the form
$\Sigma(W',S',L')$ where $S'\subset S$, $W'\subset W$, and $L'\subset
L$.  Then $\Sigma'$ can be naturally identified with a subcomplex of
$\Sigma$.  Moreover, the set $\cR'=\cR\cap\cP(L')$ is an admissible subset of 
$\cP(L')$, and the $\cR'$-blow-up of $\Sigma'$ coincides with the
$\cR$-blow-up of $\Sigma'$ (where $\Sigma'$ is viewed as a subcomplex
of $\Sigma$).
\end{remark}

\begin{remark}
If $\cR$ is admissible, then so is the collection $\cR'$, obtained by
removing from $\cR$ all elements of the form $\{s\}$ where $s\in S$.
Moreover, it follows from properties (1) and (2) of the blow-up
functor that for any singleton set $\{s\}\in\cP$, the blow-up
$( \cc_{\{s\}})_{\#}$ is canonically isomorphic to $ \cc_{\{s\}}$
regardless of whether $\{s\}$ is in $\cR$ or not.  Thus, for any
subcomplex $K\subset \Sigma$, the $\cR$-blow-up of $K$ and the
$\cR'$-blow-up of $K$ are the same.  For this reason, we shall assume
for the remainder of the article that all admissible sets $\cR$
contain no singleton subsets in $\cP$. 
\end{remark}

\begin{definition}\label{def:max-bu-set}
$\cR$ is the {\em maximal blow-up set} if it contains every spherical
set in $\cP$ of cardinality at least two.
\end{definition}

\begin{definition}\label{def:min-bu-set}
$\cR$ is the {\em minimal blow-up set} if it consists of the spherical
sets $T$ in $\cP$ such that $\card(T)\geq 2$ and $\Gamma_T$ is
connected. 
\end{definition}

\subsection{The link of a vertex in $\Sigma_{\#}$}\label{ss:nested-sets}
Let $L_{\#}$ denote the link of the vertex $1$ in an $\cR$-blow-up
$\Sigma_{\#}$.  In this subsection we shall decribe the poset of
simplices in $L_{\#}$.  First we introduce a subset $S_{\#}$ of $\cP$
which will be used to index the vertices of the link.  Let   
\[S_{\#}=\cR\cup\{\{s\}\;|\; s\in S\}.\]

Next, we want to define a poset $\cN$ of subsets of $S_{\#}$.  It plays
the same role for $\Sigma_{\#}$ as does $\cP$ for $\Sigma$.  

Let $\cT$ be a subset of $S_{\#}$.  It is partially ordered by inclusion.
The {\em support of $\cT$}, denoted $\supp\cT$ is the union of all
elements of $\cT$.  The {\em $\cR$-fixed part} of $\cT$, denoted $\cT_0$ is
the set of all singletons $\{s\}\in\cT$ where
$s\not\in\supp\cT\cap\cR$.  (If $\cR$ contains no singleton sets, then
$\cT_0$ consists of the elements $\{s\}$ that are maximal in $\cT$.)

\begin{definition}\label{def:R-nested}
Given a subset $\cT\subset S_{\#}$, let $T=\supp\cT$, $T_0=\supp\cT_0$, and
let $T_1,\ldots,T_k$ be the maximal elements of $\cT\cap\cR$.  We give
an inductive definition of what it means for $\cT$ to be
{\em $\cR$-nested}; the induction is on $\card\cT$.  If
$\cT=\emptyset$, then it is $\cR$-nested.  If $\card{\cT}>0$, 
then it is $\cR$-nested if and only if the following two conditions
hold:
\begin{enumerate}
\item[(1)] $T\in\cP$ and $\{T_0,T_1,\ldots,T_k\}$ is the
$\cR$-decomposition of $T$,
\item[(2)] $\cT_{<T_i}$ is $\cR$-nested (as defined by induction) for
all $1\leq i\leq k$. 
\end{enumerate}
\end{definition}

It is clear that for each $T\in S_{\#}$, $\{T\}$ is $\cR$-nested.  The next
lemma (which is easily verified) describes which two-element subsets
of $S_{\#}$ are $\cR$-nested.

\begin{lemma}\label{lem:R-nested-pairs}
Let $T$ and $T'$ be distinct elements of $S_{\#}$.  Then $\{T,T'\}$ is
$\cR$-nested if and only if one of the following three cases holds:
\begin{enumerate}
\setlength{\itemindent}{.4in}
\item[Case 1:] $T=\{s\}$ and $T'=\{s'\}$ where $s,s'\in S$,
$m(s,s')<\infty$, and $\{s,s'\}\not\in\cR$.
\item[Case 2:] $\{T,T'\}$ is the $\cR$-decomposition of $T\cup T'$.  
\item[Case 3:] $T'\subset T$ (or $T\subset T'$).
\end{enumerate}
\end{lemma}

Let $\cN$ be the set of all $\cR$-nested subsets of $S_{\#}$, partially
ordered by inclusion.  It is not difficult to see that any subset of
an $\cR$-nested set is $\cR$-nested.  In other words,
$\cN_{>\emptyset}$ is an abstract simplicial complex with vertex
set $S_{\#}$.  A subset of $S_{\#}$ spans a simplex of
$\cN_{>\emptyset}$ if and only if it is $\cR$-nested.  The edges of
$\cN_{>\emptyset}$ are described explicitly in
Lemma~\ref{lem:R-nested-pairs}.  It is proved in Section 3.3 of
\cite{DJS} that $\cN_{>\emptyset}$ is a certain simplicial
subdivision of $\cP_{>\emptyset}$ and that it can be identified with
the link $L_{\#}$ of the vertex $1$ in $\Sigma_{\#}$. 
 
\subsection{Nonpositive curvature} \label{ss:flag-complex}  
Suppose that $(W,S)$ is a Coxeter system, that $L$ is a subcomplex of
$\cS_{>\emptyset}$, that $\Sigma=\Sigma(W,S,L)$ (as in
Definition~\ref{def:reflection-tiling}), and that $\cR$ is an
admissible subset of $\cP(L)$.  Let $\Sigma_{\#}$ denote the
$\cR$-blow-up of $\Sigma$.  In this section, we want to determine when
the natural piecewise Euclidean metric on $\Sigma_{\#}$ is
nonpositively curved.  Throughout this section, we shall assume, for
simplicity, {\em that $\cR$ is fully admissible}
(cf. Definition~\ref{def:admissible}), so that $\Sigma_{\#}$ is a
cubical cell complex.  Let $L_{\#}$ denote the link of the vertex $1$
in $\Sigma_{\#}$, as described in the previous section.

In \cite{G} Gromov showed that a cubical cell complex is nonpositively
curved if and only if the link of each of its vertices is a flag
complex.  (Recall that a simplicial complex $K$ is a {\em flag
complex} if any finite, nonempty, collection of vertices in $K$ that
are pairwise connected by edges spans a simplex in $K$.)  Therefore,
we have the following lemma.

\begin{lemma}
$\Sigma_{\#}$ is nonpositively curved if and only if $L_{\#}$ is a
flag complex.
\end{lemma}

\subsection{Condition (F)} 
Consider the following condition on a spherical set $T$ in $\cS_{>\emptyset}$.
\begin{enumerate}
\item[(F)] Suppose there exists a decomposition $\{T_1,\ldots,T_k\}$
of $T$ such that 
\begin{enumerate}
\item[(1)] $T_i\in S_{\#}$ for $1\leq i\leq k$, 
\item[(2)] $T_i\cup T_j\in\cP-\cR$ for $1\leq i<j\leq k$, and 
\item[(3)] $k\geq 3$.  
\end{enumerate}
Then $T\in\cP-\cR$. 
\end{enumerate}

In \cite{DJS}, we analyzed when links of vertices in blow-ups were
flag complexes.  Here we are interested in the case where each such
link is isomorphic to $L_{\#}$.  The analysis in \cite{DJS} yields the
following.

\begin{theorem}
$\Sigma_{\#}$ is nonpositively curved if and only if Condition (F)
holds for each $T\in\cS_{>\emptyset}$.
\end{theorem}

\begin{proof}{}
We want to see that Condition (F) is equivalent to the condition that
$L_{\#}$ is a flag complex.  To check this we must consider
collections $\{T_1,\ldots,T_k\}$ where each $T_i$ is a vertex of
$L_{\#}$, and each pair $\{T_i,T_j\}$ spans an edge of $L_{\#}$ and
then show that $\{T_1,\ldots,T_k\}$ is $\cR$-nested.  To verify this,
we only need consider the case where the $T_i$'s are pairwise
incomparable, and this is precisely the case covered by Condition
(F). 
\end{proof}

\begin{corollary}
Suppose $\Sigma=\Sigma(W,S,L)$ and that $\cR$ is the maximal blow-up
set (cf., Definition~\ref{def:max-bu-set}).  Then $\Sigma_{\#}$ is
nonpositively curved.
\end{corollary}

\begin{proof}{}
Condition (F) holds for any $T\in\cS_{>\emptyset}$ (since condition (2) will
never be satisfied).
\end{proof}
  
\begin{corollary}
Suppose that $\Sigma=\Sigma(W,S,L)$ and that $\cR$ is the minimal
blow-up set (cf., Definition~\ref{def:min-bu-set}).  Then
$\Sigma_{\#}$ is nonpositively curved if and only if the following
condition holds: given any three completely disjoint spherical sets
$T_1,T_2,T_3$ in $\cP(L)_{>\emptyset}$ with $\Gamma_{T_i}$ connected
for $i=1,2,3$ and with $T_i\cup T_j\in\cP(L)$ for
$\{i,j\}\subset\{1,2,3\}$, then $T_1\cup T_2\cup T_3\in\cP(L)$.
\end{corollary} 

\begin{corollary}
Suppose $\Sigma=\Sigma(W,S)$.  Then the minimal blow-up of $\Sigma$ is
nonpositively curved.
\end{corollary}

\section{The dual tiling of $\Sigma_{\#}$}\label{s:dual}

\subsection{The role of $D_{\#}$}
The quick and clean discription of the blow-up functor is the one
given above in \ref{ss:blow-up-functor}, in terms of blowing up
certain collections 
of Coxeter cells.  However, this description leaves obscure an
important aspect of the construction, namely, the definition of the
fundamental domain (or ``fundamental tile'') $D_{\#}$ for the
$W$-action on $\Sigma_{\#}$.  In fact, many of the geometric ideas about
these blow-ups are best explained in terms of $D_{\#}$.  The goal of
this section is to give the definition of $D_{\#}$ and discuss some
of its properties.  The case where $W$ is finite is particularly
simple since $D_{\#}$ can be thought of as a convex polytope.  So we
shall deal with this case first.

\subsection{The case where $W$ is finite}\label{ss:finite-W-tiles}
When $W$ is finite, $D$ is the Coxeter block of type $(W,S)$.  It is a
fundamental domain for $W$ on $ \cc$.  The {\em tile} 
$D_{\#}$ is obtained from $D$ by truncating those faces that
correspond to elements of $\cR$.  Thus, $D_{\#}$ is combinatorially
equivalent to a convex polytope.  The mirrors of $D_{\#}$ either
correspond to the original mirrors of $D$ (these are indexed by $S$)
or to the new codimension-one faces introduced in the truncation
process (these are indexed by $\cR$).

Similarly, the simplex $\Delta$ is a fundamental domain for $W$ on
$\Bd \cc$.  By Remark~\ref{rem:bu-compatibility}, the $\cR$-blow-up of
the $\Bd \cc$ (viewed as a subcomplex of $\cc$) coincides with the
$\cR'$-blow-up (where $\cR'=\cR-\{S\}$) of $\Bd \cc$ (viewed as a
reflection tiling in its own right).  The {\em fundamental tile}
$\Delta_{\#}$ is obtained by truncating the faces of $\Delta$
corresponding to elements of $\cR'$.  (Here we are thinking of
$\Delta$ simply as a convex simplex, that is, we are temporarily
ignoring its subdivision into Coxeter blocks.)  Again, $\Delta_{\#}$
is a convex polytope, each codimension-one face is a mirror, and the
mirrors are indexed either by the elements of $S$ or by elements of
$\cR'$.   

\begin{example} 
The {\em maximal blow-up} of $ \cc$ means the case where
$\cR=\cS_{>\emptyset}$.  In this case, $\Delta_{\#}$ is a
permutohedron and $D_{\#}=\Delta_{\#}\times[0,1]$.
\end{example}

\begin{example}\label{ex:minbu-assoc}
By the {\em minimal blow-up} of $ \cc$, we mean the case where
$T\in\cR$ if and only if the Coxeter diagram $\Gamma(W_T,T)$ is
connected.  If the Coxeter diagram is a straight line segment
(i.e., if $\Gamma(W,S)$ is of type $\bA_n$, $\bB_n$, $\bI_2(p)$,
$\bH_3$, $\bH_4$, or $\bF_4$), then it turns out that $\Delta_{\#}$ is
an {\em associahedron} (as defined by \cite{Sta}).  (See
Section~\ref{s:associahedral-tilings}.)  
\end{example}

The next lemma gives an alternative description of $D_{\#}$ in terms
of Coxeter blocks.  Its proof follows from the discussion in
\ref{ss:blow-up-cell}. 

\begin{lemma}\label{lem:D-sharp}
$D_{\#}$ is the union of all Coxeter blocks in $ \cc_{\#}$ that
contain the vertex $1$.
\end{lemma}

It follows from this lemma that $D_{\#}$ is a fundamental domain for
the $W$-action on $ \cc_{\#}$ in the following sense:  the $W$-orbit
of any point $x$ in $ \cc_{\#}$ intersects $D_{\#}$ in at least one
point; moreover, this intersection is exactly one point if $x$ does
not belong to a mirror of any Coxeter block.

\subsection{The general case}
We now allow $W$ to be infinite and consider its action on the complex
$\Sigma_{\#}$.  As before, we choose a vertex of $\Sigma_{\#}$ and denote it by
$1$.  Lemma~\ref{lem:D-sharp} suggests the following.

\begin{definition}
The {\em fundamental tile} $D_{\#}$ for $W$ on $\Sigma_{\#}$ is the union
of all Coxeter blocks in $\Sigma_{\#}$ that contain the vertex $1$.
\end{definition}

As before, we see that $D_{\#}$ is a fundamental domain for the
$W$-action on $\Sigma_{\#}$.  We will decribe the mirrors of $D_{\#}$
in the next subsection.

\subsection{Mirrors}\label{ss:mirrors}
Suppose that $L$ is a subcomplex of $\mathcal{S}_{>\emptyset}$.  
We begin by recalling some facts about the poset $\cP=\cP(L)$.  We
note that $\cP_{>\emptyset}$ is an abstract simplicial complex, in
other words, the vertex set of $\cP_{>\emptyset}$ is $S$, and a subset
$T\subset S$ spans a simplex if and only if $T\in\cP_{>\emptyset}$.
The poset $\cP$ plays two roles in the description of $\Sigma$.  First
of all 
$\cP_{>\emptyset}$ can be identified with the link of a vertex in the
cellulation of $\Sigma$ by Coxeter cells.  Secondly, if the
fundamental chamber $D$ is defined to be the union of all Coxeter
blocks in $\Sigma$ that contain the vertex $1$, then $D$ can be
identified with the geometric realization of $\cP$. 
(If $\cP$ is a poset of subsets of some set $V$, if $\emptyset\in\cP$,
and if $\cP_{>\emptyset}$ is an abstract simplicial complex, then for
any $T\in\cP$, the geometric realization of $\cP_{\leq T}$ is
isomorphic to a standard simplicial subdivision of a cube of dimension
$\card T$.  It follows that the geometric realization of $\cP$ can be
identified with the union of all such cubes.)  Having made this
identification, the mirror of $D$ corresponding to $s\in S$ is then
identified with the geometric realization of $\cP_{\geq\{s\}}$.

We now apply the same construction to the blow-up $\Sigma_{\#}$.  In
this case, the link of the vertex $1$ is associated to the abstract
simplicial complex $\cN_{>\emptyset}$ where $\cN$ is the poset of
$\cR$-nested subsets of $\cP$.  It follows that there is
one Coxeter block $B_{\cT}$ containing $1$ for each $\cR$-nested set
$\cT$. (The dimension of $B_{\cT}$ is $\card\cT$.)  In other words,
$D_{\#}$ can be identified with the geometric realization of $\cN$.  

\begin{definition}
For each $T\in S_{\#}$, the {\em mirror of $D_{\#}$ corresponding to $T$},
denoted $D_{\#T}$, is the geometric realization of $\cN_{\geq\{T\}}$.
\end{definition}

Each Coxeter block in $\Sigma_{\#}$ has a well-defined {\em type} -- it is an
element of $\cN$.  (Translate the Coxeter block by an element of $W$
so that it contains the vertex $1$.)  Each $1$-dimensional Coxeter
cell also has a well-defined {\em type} -- it is an element of $S_{\#}$.
(In Section~\ref{s:tilings} the function that assigns to each oriented edge of 
$\Sigma_{\#}$ the type of its underlying 1-cell will be called a 
``framing''.)  In general, however, there is no consistent way to assign
a ``type'' to the Coxeter cells in $\Sigma_{\#}$ of dimension $\geq
2$. (See Case 3 of Figure~\ref{fig:cells-2-dim}.)

\subsection{The involution $j_T$.}\label{ss:blow-up-gluing}
Let $T\in S_{\#}$.  Next we shall describe the self-homeomorphism of
$D_{\#T}$ which is needed to glue together the tiles $D_{\#T}$ and
$w_TD_{\#T}$.  To this end we want to define an appropriate extension
of the automorphism $j_T:\cR_{\leq T}\rightarrow\cR_{\leq T}$, defined
in Lemma~\ref{lem:antipodal}, to an automorphism of $\cN_{\geq\{T\}}$, the
closed star of the vertex $T$ in $\cN_{>\emptyset}$.  This
automorphism will also be denoted by $j_T$.

For any $T'\in\cP_{>\emptyset}$, let $a_{T'}$ denote the antipodal
map on the Coxeter cell $ \cc(W_{T'},T')$.  If $T'\in S_{\#}$, then it
corresponds to a Coxeter $1$-cell connecting the vertex $1$ to
$a_{T'}(1)=w_{T'}(1)$.  The gluing map $j_T$ is induced by the action
of $a_Tw_T$ on those vertices $T'$ that lie in $\cN_{\geq\{T\}}$.  By
Lemma ~\ref{lem:R-nested-pairs} there are three cases to consider.
First we consider Case 3, when $T'\subset T$.  The face
$\cc(W_{T'},T')$ of $ \cc(W_T,T)$ is mapped by $a_Tw_T$ to the
face $ \cc(W_{T''},T'')$, where $T''=w_TT'w_T$.  Hence, $a_Tw_T$
maps the vertex $a_{T'}(1)$ to $a_{T''}(1)$.  So in this case, the
correct definition of $j_T$ is as in Lemma~\ref{lem:antipodal}: $j_T(T')=T''$.
On the other hand, if $T\subset T'$, then $a_{T'}$ commutes with $a_T$
and with $w_T$.  Hence $a_Tw_T$ maps $a_{T'}(1)$ to
$a_{T'}a_Tw_T(1)=a_{T'}(1)$.  Thus, in this case, we see that the
appropriate definition is: $j_T(T')=T'$.  A similar argument shows
that in Cases 1 and 2, the appropriate definition is $j_T(T')=T'$.  Thus
$j_T$ is defined on any vertex $T'$ of $\cN_{\geq\{T\}}$ by 
\[j_T(T')=\left\{\begin{array}{cl}
w_TT'w_T &\mbox{if $T'\subset T$}\\ 
T' &\mbox{if $T'\not\subset T$.}\end{array}\right.\]    
It is then not difficult to check that for any $\cR$-nested set
$\cT\in\cN_{\geq\{T\}}$, $j_T(\cT)$ is also $\cR$-nested.  That is to
say, $j_T$ induces an automorphism of $\cN_{\geq\{T\}}$.  The
geometric realization of $j_T$ is a self-homeomorphism of the mirror
$D_{\#T}$, which we shall continue to denote by $j_T$.

\begin{remark}
In general, $j_T$ will not extend to an automorphism of $\cN$.
\end{remark}

\begin{remark}
The subspace $D_{\#}\cup w_TD_{\#}$ of $\Sigma_{\#}$ is homeomorphic to two
copies of $D_{\#}$ glued together along the mirror $D_{\#T}$ via the
homeomorphism $j_T:D_{\#T}\rightarrow D_{\#T}$.
\end{remark}

\begin{definition}
The mirror $D_{\#T}$ is a {\em reflecting mirror} if $j_T$ is the
identity map on $D_{\#T}$.  If $j_{T}$ is not the identity, then
$D_{\#T}$ is a {\em mock reflecting mirror}.  
\end{definition}

\subsection{Local pictures around codimension-two
corners}\label{ss:codim-2-pictures}
By Lemma~\ref{lem:R-nested-pairs}, the intersection of distinct
mirrors $D_{\#T}$ and $D_{\#T'}$ is nonempty if and only if one of the
three cases in the lemma holds.  In Case 1 and Case 2 the picture of
the tiles around $D_{\#T}\cap D_{\#T'}$ is the usual two-dimensional
picture for reflection groups.  In Case 3 it is slightly different.
All three cases are depicted in Figure~\ref{fig:tiling-codim-2} below.
In the figure we have labeled the tile $wD_{\#}$ by the corresponding
element $w$ and the mirror $wD_{\#T}$ by the corresponding element
$T\in S_{\#}$.

\begin{figure}[ht]
\begin{center}
\psfrag{1}{\scriptsize $1$}
\psfrag{s}{\scriptsize $s$}
\psfrag{s'}{\scriptsize $s'$}
\psfrag{ss'}{\scriptsize $ss'$}
\psfrag{ss's=s'ss'}{\scriptsize $ss's=s'ss'$}
\psfrag{s's}{\scriptsize $s's$}
\psfrag{{s}}{\scriptsize $\{s\}$}
\psfrag{{s'}}{\scriptsize \hspace{-.1in}$\{s'\}$}
\psfrag{C1}{\hspace{-.6in}{\bf Case 1} ($m(s,s')=3$)}
\psfrag{T}{\scriptsize $T$}
\psfrag{T'}{\scriptsize $T'$}
\psfrag{T''}{\scriptsize $T''$}
\psfrag{w_T}{\scriptsize $w_T$}
\psfrag{w_T'}{\scriptsize $w_{T'}$}
\psfrag{w_Tw_T'=w_T'w_T}{\scriptsize $w_Tw_{T'}=w_{T'}w_T$}
\psfrag{C2}{\hspace{-.2in}{\bf Case 2}}
\psfrag{w_Tw_T''=w_T'w_T}{\scriptsize $w_Tw_{T''}=w_{T'}w_T$}
\psfrag{C3}{\hspace{-.8in}{\bf Case 3} ($T'\subset T$, $T''=j_T(T')$)} 
\includegraphics[scale = .4]{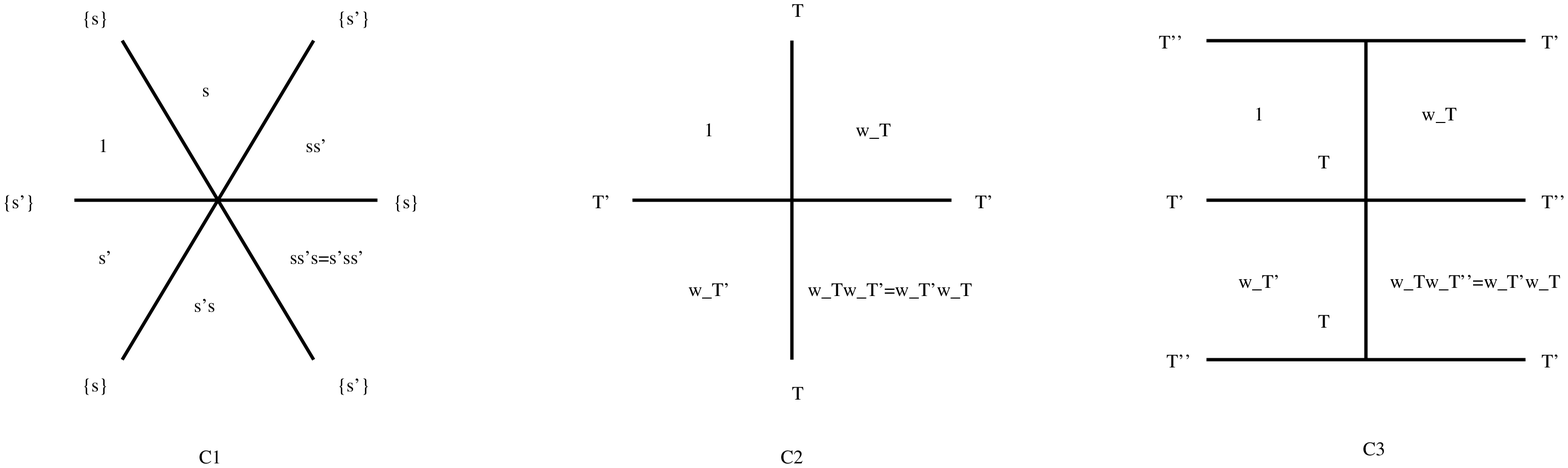}
\caption{\label{fig:tiling-codim-2}}
\end{center}
\end{figure}

\begin{figure}[ht]
\begin{center}
\psfrag{1}{\scriptsize $1$}
\psfrag{s}{\scriptsize $s$}
\psfrag{s'}{\scriptsize $s'$}
\psfrag{ss'}{\scriptsize $ss'$}
\psfrag{ss's=s'ss'}{\scriptsize $ss's=s'ss'$}
\psfrag{s's}{\scriptsize $s's$}
\psfrag{{s}}{\scriptsize $\{s\}$}
\psfrag{{s'}}{\scriptsize $\{s'\}$}
\psfrag{C1}{\hspace{-.2in}{\bf Case 1}}
\psfrag{T}{\scriptsize $T$}
\psfrag{T'}{\scriptsize $T'$}
\psfrag{T''}{\scriptsize $T''$}
\psfrag{w_T}{\scriptsize $w_T$}
\psfrag{w_T'}{\scriptsize $w_{T'}$}
\psfrag{w_Tw_T'=w_T'w_T}{\scriptsize $w_Tw_{T'}=w_{T'}w_T$}
\psfrag{C2}{\hspace{-.3in}{\bf Case 2}}
\psfrag{w_Tw_T''=w_T'w_T}{\scriptsize $w_Tw_{T''}=w_{T'}w_T$}
\psfrag{C3}{\hspace{-.2in}{\bf Case 3}}
\includegraphics[scale = .55]{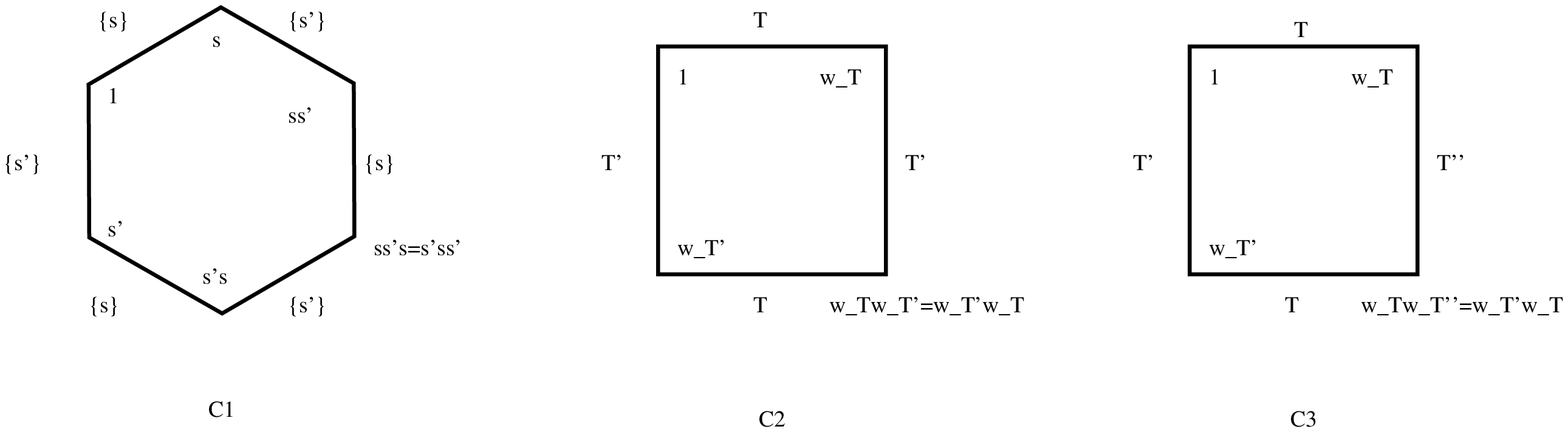}
\end{center}
\caption{\label{fig:cells-2-dim}}
\end{figure}

There are similar pictures for the corresponding Coxeter $2$-cells
containing the vertex $1$ as in Figure~\ref{fig:cells-2-dim}, below.
Here we have labeled the vertices by group elements and the Coxeter
$1$-cells by their corresponding type.

Since $D_{\#}$ is the dual cone in a Coxeter cell complex, it is a
Coxeter tile as defined in \ref{def:coxeter-tile}.  Given the local
pictures around $2$-dimensional cells, we can now determine the type 
of this dual cone (cf., Section~\ref{ss:local-geometry}).  We define a Coxeter matrix
$M_{\#}$($=(m_{\#}(T,T'))$) on $S_{\#}$ according to the three cases
of Lemma~\ref{lem:R-nested-pairs}:
\[m_{\#}(T,T')=\left\{\begin{array}{ll}
1 & \mbox{if $T=T'$}\\
m(s,s') & \mbox{in Case 1}\\
2 & \mbox{in Cases 2 and 3}\\
\infty &\mbox{if $\{T,T'\}\not\in\cN$.}\end{array}\right.\]
Let $(W_{\#},S_{\#})$ be the Coxeter system determined by $M_{\#}$.
With respect to this Coxeter system, we have the notion of a spherical
subset of $S_{\#}$  (namely, $\cT\subset S_{\#}$ is spherical if and
only if the special subgroup $(W_{\#})_{\cT}$ is finite).  Let
$\cS(S_{\#})$ denote the set of spherical subsets of $S_{\#}$.  It is
not hard to see that the simplicial complex $\cN_{>\emptyset}$ is a
subcomplex of $\cS(S_{\#})_{>\emptyset}$ (i.e., that every
$\cR$-nested set $\cT$ is spherical).  Denoting this subcomplex
$L_{\#}$, we then have the following.

\begin{proposition}
The fundamental tile $D_{\#}$ is isomorphic to the Coxeter tile\\
$D(W_{\#},S_{\#},L_{\#})$ (as defined in \ref{def:coxeter-tile}).
\end{proposition}

\subsection{The universal cover of $\Sigma_{\#}$ and the group
$A$}\label{ss:group-A}
Let $p:\tilde{\Sigma}_{\#}\rightarrow \Sigma_{\#}$ denote the projection map of
the universal cover $\tilde{\Sigma}_{\#}$ onto $\Sigma_{\#}$.  Let $A$
denote the group of all lifts of the $W$-action on $\Sigma_{\#}$ to
$\tilde{\Sigma}_{\#}$.  The projection map
$p:\tilde{\Sigma}_{\#}\rightarrow \Sigma_{\#}$ induces a surjective
homomorphism $\phi:A\rightarrow W$, and $\pi_1(\Sigma_{\#})$ is
naturally identified with the kernel of $\phi$. The goal of this
subsection is to give a presentation for $A$.  The method is
essentially due to Poincar\'{e}.    

The cellulation of $\Sigma_{\#}$ by Coxeter cells lifts to a
cellulation of $\tilde{\Sigma}_{\#}$.  Thus, $\tilde{\Sigma}_{\#}$ is
a Coxeter cell complex.  Since $D_{\#}$ is contractible, $p$ maps each
component of $p^{-1}(D_{\#})$ homeomorphically onto $D_{\#}$.  Hence, 
$\tilde{\Sigma}_{\#}$ is also tiled by copies of $D_{\#}$.  First
choose a component of $p^{-1}(D_{\#})$ and denote it by
$\tilde{D}_{\#}$.  Let $\tilde{1}$ denote the point of
$\tilde{D}_{\#}$ that lies above the vertex $1$ in $D_{\#}$.
($\tilde{D}_{\#}$ should be thought of as the Dirichlet domain
centered at $\tilde{1}$ for the $A$-action on $\tilde{\Sigma}_{\#}$.)
The set of elements of $A$ that map $\tilde{D}_{\#}$ to an adjacent
chamber across a mirror is a set of generators for $A$.  A set of
relations can be read off by considering the local pictures around the
intersections of two mirrors. 

There is a dual and equivalent method to this which is easier to make
mathematically precise.  Consider the $2$-skeleton of $\tilde{\Sigma}_{\#}$
in its Coxeter cell structure.  Since $W$ acts simply transitively on
the vertex set of $\Sigma_{\#}$, $A$ acts simply transitively on the vertex
set of $\tilde{\Sigma}_{\#}$.  Since we have chosen a distinguished vertex
$\tilde{1}$ in $\tilde{\Sigma}_{\#}$, each vertex of $\tilde{\Sigma}_{\#}$ is
labeled by an element of $A$.  Each $1$-cell in $\tilde{\Sigma}_{\#}$ is
then labeled by its type (an element of $S_{\#}$).  The set of edges in the
star of the vertex $\tilde{1}$ then gives a set of generators for $A$
while the set of $2$-cells in this star gives a set of relations.

We now give the details of this method.  For each $T\in S_{\#}$, let
$\alpha_T$ denote the unique lift of $w_T$ that maps $\tilde{D}_{\#}$
to the adjacent tile across $\tilde{D}_{\#T}$.  (Here is a more
precise definition of $\alpha_T$.  Let $x_T$ be the midpoint of the
edge in $\Sigma_{\#}$ connecting the vertex $1$ to $w_T(1)$.  Then $w_T$
fixes $x_T$.  Let $\tilde{x}_T$ denote the lift of $x_T$ in
$\tilde{D}_{\#}$.  Then $\alpha_T$ is defined to be the unique lift of
$w_T$ that fixes $\tilde{x}_T$.)  Since $(\alpha_T)^2$ also fixes
$\tilde{x}_T$ and covers the identity map on $\Sigma_{\#}$, $(\alpha_T)^2$
must be the identity on $\tilde{\Sigma}_{\#}$, i.e., $\alpha_T$ is an
involution.

\begin{definition} 
The involution $\alpha_T$ is called a {\em reflection} if $j_T=\Id$
and a {\em mock reflection} otherwise.
\end{definition}

\begin{theorem}\label{thm:mock-presentation}
The set $\cA=\{\alpha_T\}_{T\in S_{\#}}$ is a set of generators for $A$.
Moreover, the following relations give a presentation for $A$:
\begin{enumerate}
\item[(0)] $(\alpha_T)^2=1$ for all $T\in S_{\#}$.
\item[(1)] $(\alpha_{\{s\}}\alpha_{\{s'\}})^{m(s,s')}=1$ whenever
$m(s,s')<\infty$ and $\{s,s'\}\not\in\cR$.
\item[(2)] $(\alpha_T\alpha_{T'})^2=1$ whenever $T\cup T'\in\cS$ and
$\{T,T'\}$ is the $\cR$-decomposition of $T\cup T'$.
\item[(3)] $\alpha_T\alpha_{T'}=\alpha_{T''}\alpha_T$ whenever
$T'\subset T$ (where $T''=w_TT'w_T$).
\end{enumerate}
\end{theorem}

\begin{proof}{}
For each $T\in S_{\#}$ there is a Coxeter $1$-cell connecting $\tilde{1}$
to $\alpha_T(\tilde{1})$.  Since the $1$-skeleton of $\tilde{\Sigma}_{\#}$
is connected, $\cA$ is a set of generators for $A$.

As we have already observed, $\alpha_T$ is an involution, so the
relations in (0) hold in $A$.  Examining the local pictures around the
intersection of two mirrors in Figure~\ref{fig:tiling-codim-2}, we see
that the relations of type (1), (2), and (3) hold in $A$.  In terms of
the $2$-skeleton of $\tilde{\Sigma}_{\#}$, the relations of the form
(1), (2), and (3) correspond to the $2$-dimensional Coxeter cells that
contain the vertex $\tilde{1}$.  Hence, the Cayley $2$-complex of the
group defined by this presentation is a covering space of the
$2$-skeleton.  Since the $2$-skeleton is simply-connected, this
covering must be trivial; hence, this presentation is a presentation
for $A$.
\end{proof}

\begin{remark}  
The homomorphism $\phi:A\rightarrow W$ is defined on
the generating set $\cA$ by $\phi(\alpha_T)=w_T$.
\end{remark}

\section{Tilings of Coxeter cell complexes}\label{s:tilings}

\subsection{Framings}
A framing is an additional structure on a Coxeter cell complex
$\ccc$.  It is used to rigidify the situation and to provide a notion
of ``parallel transport'' along curves.  The existence of a framing
implies that $X$ is tiled (cf., Definition~\ref{def:tiling}).  In the
converse direction, if $\ccc$ admits a group of automorphisms that acts
simply transitively on its vertex set $\ccc^{(0)}$, then there is an
associated framing (cf., Example~\ref{ex:mainframe}, below).

Recall that an {\em orientation} for an edge of $\ccc$ is an ordering
of its two 
endpoints.  Let $OE(\ccc)$ denote the set of oriented edges in $\ccc$.
For any $e\in OE(\ccc)$, let $i(e)$ denote its first endpoint (its
``initial vertex'') and $t(e)$ its second (``terminal vertex'').
Also, $\bar{e}$ denotes the same edge with the opposite orientation.
An {\em edge path} in $\ccc$ is a sequence $\be=(e_1,\ldots,e_n)$
of oriented edges such that $i(e_{k+1})=t(e_k)$, for $1\leq k<n$.  Its
{\em initial vertex} $i(\be)$ is defined to be $i(e_1)$ and its {\em
terminal vertex} $t(\be)$ is $t(e_n)$.

\begin{definition}\label{def:framing-system}
A {\em framing system} is a triple $(V,M,L)$ where 
\begin{enumerate}
\item[(i)] $V$ is a finite set,
\item[(ii)] $M$ ($=m(v,v')$) is a Coxeter matrix on $V$, and 
\item[(iii)] $L$ is a subcomplex of $\cS(V)_{>\emptyset}$ containing
its $1$-skeleton.
\end{enumerate}
\end{definition}

There is an obvious notion of an {\em isomorphism} between two framing
systems $(V,M,L)$ and $(V',M',L')$  (namely, it is a bijection
$\phi:V\rightarrow V'$ that induces a simplicial 
isomorphism $L\rightarrow L'$ and pulls back $M'$ to $M$).

Associated to a framing system $(V,M,L)$ there is the reflection
tiling $\Sigma(W(M),V,L)$ defined in
\ref{def:reflection-tiling}. 

\begin{example}
Associated to any vertex $\cp\in\ccc^{(0)}$, we have three sets
$E_{\cp}$, $O_{\cp}$, and $I_{\cp}$ called, respectively, the {\em
unoriented edges}, the {\em outward pointing edges}, and the {\em
inward pointing edges} at $\cp$.  Their definitions are:
\[\begin{array}{l}
E_{\cp}=\{\mbox{unoriented edges with one endpoint $\cp$}\}\\
O_{\cp}=\{e\in OE\;|\; i(e)=\cp\}\\
I_{\cp}=\{e\in OE\;|\; t(e)=\cp\}.
\end{array}\]
We note that these three sets are canonically isomorphic.  A Coxeter
matrix $M_{\cp}$ ($=m_{\cp}(e,e')$) on $E_{\cp}$ (or on $O_{\cp}$ or
$I_{\cp}$) is defined as follows: if $e=e'$, then $m_{\cp}(e,e')=1$;
if $e$ and $e'$ are distinct edges of a $2$-cell in $\ccc$ that is a
$2m$-gon, then $m_{\cp}(e,e')=m$; otherwise, $m_{\cp}(e,e')=\infty$.
Let $W(M_{\cp})$ denote the Coxeter group corresponding to $M_{\cp}$
with fundamental set of generators $E_{\cp}$, and let $\cS(E_{\cp})$
be the poset of spherical subsets of $E_{\cp}$.  The link, $L_{\cp}$,
of $\cp$ in $\ccc$ is then a subcomplex of $\cS(E_{\cp})$.  Thus, to
each $\cp\in\ccc^{(0)}$ we have associated three canonically
isomorphic framing systems: $(E_{\cp},M_{\cp},L_{\cp})$,
$(O_{\cp},M_{\cp},L_{\cp})$, and $(I_{\cp},M_{\cp},L_{\cp})$.
\end{example} 
 
\begin{definition}
Suppose $\ccc$ is a Coxeter cell complex, that $(V,M,L)$ is a framing
system and that $\nabla:OE(\ccc)\rightarrow V$ is a function.  For
each $\cp\in\ccc^{(0)}$, denote the restriction of $\nabla$ to
$O_{\cp}$ and $I_{\cp}$ by $\nabla_{\cp}:O_{\cp}\rightarrow V$ and
$\bar{\nabla}_{\cp}:I_{\cp}\rightarrow V$, respectively.  Then
$\nabla$ is called a {\em framing} if for each $\cp\in\ccc^{(0)}$,
both $\nabla_{\cp}$ and $\bar{\nabla}_{\cp}$ are isomorphisms of
framing systems.
\end{definition}

\begin{remark}\label{rem:iota}
The canonical isomorphism $c_{\cp}:O_{\cp}\rightarrow I_{\cp}$ is
defined by $c_{\cp}(e)=\bar{e}$.  If $\nabla$ is a framing on $\ccc$,
then for each $\cp\in\ccc^{(0)}$, we get an involution of framing
systems $\iota_{\cp}:V\rightarrow V$ defined by the condition that the
following diagram commutes
\[\xymatrix{O_{\cp}\ar[d]^{\nabla_{\cp}}\ar[r]^{c_{\cp}} &
I_{\cp}\ar[d]^{\bar{\nabla}_{\cp}}\\
V\ar[r]^{\iota_{\cp}}& V}\]
In all cases of interest in this paper, the involution $\iota_{\cp}$
will be the identity map for each $\cp\in\ccc^{(0)}$.
\end{remark}

\begin{example}\label{ex:blow-up-framing}
Suppose that $(W,S)$ is a Coxeter system, that $\cR$ is an admissible
subset of $\cS$, and that $S_{\#}$ is the subset of $\cS$ defined in
\ref{ss:nested-sets}.  Let $M_{\#}$ be the Coxeter
matrix on $S_{\#}$ defined in \ref{ss:codim-2-pictures}, and let
$L_{\#}$ be the simplicial complex $\cN_{>\emptyset}$.  Then
$(S_{\#},M_{\#},L_{\#})$ is a framing system.  If $\Sigma_{\#}$
denotes the $\cR$-blowup of $\Sigma(W,S)$, then there is a natural
framing on $\Sigma_{\#}$ with framing system $(S_{\#},M_{\#},L_{\#})$
which associates to each edge its type (as defined in
\ref{ss:mirrors}).  
\end{example}

\begin{example}\label{ex:mainframe}
Suppose $A$ is a group of automorphisms of a Coxeter cell complex
$\ccc$ that acts simply transitively on $\ccc^{(0)}$.  Set
$V=O_{\cp}$, $M=M_{\cp}$, and $L=L_{\cp}$.  A
framing $\nabla:OE(\ccc)\rightarrow V$ is defined as follows.  For any
$e\in OE(\ccc)$, let $a$ be the unique element of $A$ that takes
$i(e)$ to $\cp$.  Then $\nabla(e)=a(e)$.
\end{example}

\begin{example}
Suppose $\nabla$ is a framing of $\ccc$ with framing system
$(V,M,L)$.  Let $(V',M',L')$ be another framing system and
$\phi:V\rightarrow V'$ an isomorphism.  Then there is an induced
framing $\nabla^{\phi}$ of $\ccc$ with system $(V',M',L')$ defined by
$\nabla^{\phi}(e)=\phi(\nabla(e))$.
\end{example}

\begin{example}
Suppose that $p:\ccc\rightarrow\ccc'$ is a covering projection of
Coxeter cell complexes and that $\nabla'$ is a framing on $\ccc'$.  We
define a framing $p^*\nabla'$ on $\ccc$ (with the same framing system
as $\nabla'$) by the formula $(p^*\nabla')(e)=\nabla'(p(e))$.
\end{example}

\begin{definition}\label{def:framed-map}
Suppose that $(\ccc,\nabla)$ and $(\ccc',\nabla')$ are framed Coxeter
cell complexes with the same framing systems.  A {\em framed map} from
$(\ccc,\nabla)$ to $(\ccc',\nabla')$ is a covering projection
$p:\ccc\rightarrow\ccc'$ such that $\nabla=p^*\nabla'$.
\end{definition}

\begin{remark}\label{rem:framed-map-uniqueness}
Suppose $p_1$ and $p_2$ are two framed maps from $(\ccc',\nabla')$ to
$(\ccc,\nabla)$ and that for some vertex $\cp'$ of $\ccc'$,
$p_1(\cp')=p_2(\cp')$.  Since the maps preserve the framing, this
implies that $p_1$ and $p_2$ also agree on any vertex adjacent to
$\cp'$.  (Two vertices are {\em adjacent} if they are connected by an
edge.)  Hence, if $\ccc'$ is connected, we must have $p_1=p_2$.  In
other words, provided $\ccc'$ is connected, if two framed maps agree
at a single vertex, then they are equal.
\end{remark}

Henceforth, we shall assume that all Coxeter cell complexes are
connected.  

As a special case of Definition~\ref{def:framed-map}, note that any
automorphism $f:\ccc\rightarrow\ccc$ is a covering projection.  It is
{\em framed} if $f^*\nabla=\nabla$.  Let $\Aut(\ccc,\nabla)$ denote
the group of framed automorphisms of $(\ccc,\nabla)$.  By
Remark~\ref{rem:framed-map-uniqueness}, $\Aut(\ccc,\nabla)$ acts freely
on $\ccc^{(0)}$.  

Suppose that $(\ccc,\nabla)$ and $(\ccc',\nabla')$ are framed with
framing systems $(V,M,L)$ and $(V',M',L')$, respectively, and that
$p:\ccc\rightarrow\ccc'$ is a covering projection (not necessarily
framed).  Then for each vertex $\cp\in\ccc^{(0)}$ there is an {\em
induced isomorphism} $p_{\cp}:V\rightarrow V'$ of framing systems
defined by the condition that the following diagram commutes:
\[\xymatrix{O_{\cp}\ar[d]^{\nabla_{\cp}}\ar[r]^{p|_{O_{\cp}}} &
O_{p(\cp)}\ar[d]^{\nabla'_{p(\cp)}}\\
V\ar[r]^{p_{\cp}}& V'}\]
Thus, after changing $(V',M',L')$ by an isomorphism of framing
systems, we may assume that $(V',M',L')=(V,M,L)$ and that at a given
vertex $\cp$, $p_{\cp}$ is the identity. 

The framing $\nabla$ is {\em symmetric} if $\Aut(\ccc,\nabla)$ acts
transitively on $\ccc^{(0)}$.  Thus, up to isomorphism of framing
systems, every symmetric framing arises from the construction in
Example~\ref{ex:mainframe}.  

\subsection{The action of $F_V$ on the vertex set}
Suppose $(\ccc,\nabla)$ is a framed Coxeter cell complex with framing
system $(V,M,L)$.  Let $F_V$ denote the free group on $V$.  Given a
vertex $\cp\in\ccc^{(0)}$ and an element $v\in V$, let $e\in O_{\cp}$
and $e^*\in I_{\cp}$ be the oriented edges defined by $\nabla_{\cp}(e)=v$
and $\bar{\nabla}_{\cp}(e^*)=v$, respectively.  Define $\cp\cdot v=t(e)$
and $\cp\cdot v^{-1}=i(e^*)$.  Clearly, $(\cp\cdot v)\cdot
v^{-1}=\cp=(\cp\cdot v^{-1})\cdot v$; hence, these formulas define a
right $F_V$-action on $\ccc^{(0)}$.  The action is transitive since
$\ccc$ is connected. 

If $f:(\ccc,\nabla)\rightarrow(\ccc',\nabla')$ is a framed map then
the restriction of $f$ to the vertex set, $f^{(0)}:\ccc^{(0)}\rightarrow
{\ccc'}^{(0)}$, is obviously $F_V$-equivariant.  Conversely, we have the
following. 

\begin{lemma}\label{lem:frame-fv-maps}
Suppose $(\ccc,\nabla)$ and $(\ccc',\nabla')$ are two framed Coxeter
cell complexes with the same framing system $(V,M,L)$.  Let
$\theta:{\ccc}^{(0)}\rightarrow{\ccc'}^{(0)}$ be an equivariant map of
$F_V$-sets.  Then there exists a unique framed map
$f:(\ccc,\nabla)\rightarrow(\ccc',\nabla')$ with $f^{(0)}=\theta$.
\end{lemma}

\begin{proof}{}
The map $\theta$ defines $f$ on the $0$-skeleton of $\ccc$.  If
$e\in OE(\ccc)$, then there is a unique oriented edge $e'$ in $\ccc'$
from $\theta(i(e))$ to $\theta(t(e))$, namely, the edge $e'$ with
$i(e')=\theta(i(e))$ and with $\nabla'(e')=\nabla(e)$.  (By
$F_V$-equivariance,
$t(e')=i(e')\cdot\nabla'(e')=\theta(i(e))\cdot\nabla(e)=\theta(t(e))$.)
We extend $f$ to the $1$-skeleton by mapping $e$ isometrically to
$e'$.  A similar argument for higher dimensional cells shows that
$\theta$ extends to a map $f$, which is clearly a covering projection.
Uniqueness follows from Remark~\ref{rem:framed-map-uniqueness}.
\end{proof}

\begin{corollary}
Suppose that $(\ccc,\nabla)$ and $(\ccc',\nabla')$ are two framed
Coxeter cell complexes with the same framing system.  Then the
correspondence $f\mapsto f^{(0)}$ is a bijection from the set of
framed maps from $(\ccc,\nabla)$ to $(\ccc',\nabla')$ with the set of
$F_V$-equivariant maps from ${\ccc}^{(0)}$ to ${\ccc'}^{(0)}$.  In
particular, $(\ccc,\nabla)$ is framed isomorphic to $(\ccc',\nabla')$
if and only if their vertex sets are isomorphic as $F_V$-sets.
\end{corollary}

\subsection{Gluing isomorphisms}\label{ss:glu-iso}
Suppose $\cc$ is a Coxeter cell and that $\cp$ is a vertex of $\cc$.
Let $O_{\cp}(\cc)$ denote the set of oriented edges of $\cc$ with
initial vertex $\cp$.  Let $W_{\cc}$ denote the Coxeter group
associated to $\cc$.  ($W_{\cc}$ is a well-defined group of symmetries
of $\cc$ generated by reflections across hyperplanes.)  Suppose $\cp$
and $\cp'$ are two vertices of $\cc$.  Then there is a unique element
$w\in W_{\cc}$ such that $w\cp=\cp'$.  The element $w$ gives us a
canonical bijection, which we denote by $w_*:O_{\cp}(\cc)\rightarrow
O_{\cp'}(\cc)$.

Now suppose that $(\ccc,\nabla)$ is a framed Coxeter cell complex,
that $\cc$ is a Coxeter cell in $\ccc$, and that $\cp$ is a vertex of
$\cc$.  Set $\nabla_{\cp}(\cc)=\nabla_{\cp}(O_{\cp}(\cc))$.  It is a
spherical subset of $V$.  If $\cp'$ is another vertex of $\cc$, then
we have a canonical isomorphism
$j_{(\cp,\cp',\cc)}:\nabla_{\cp}(\cc)\rightarrow\nabla_{\cp'}(\cc)$
defined by $j_{(\cp,\cp',\cc)}=\nabla_{\cp}\circ
w_*\circ(\nabla_{\cp'})^{-1}$, where $w\in W_{\cc}$ is such that
$\cp'=w\cp$.

For each $v\in V$ let $\star(v)$ denote the star of $v$ in the
simplicial complex $L$.  We shall now define, for each $e\in OE(\ccc)$,
an isomorphism $j_e:\star(\nabla(e))\rightarrow\star(\nabla(\bar{e}))$.  
Let $\cp_0=i(e)$ and $\cp_1=t(e)$ be the endpoints of $e$.  For any
simplex $\sigma$ in $\star(\nabla(e))$ let $\cc_{\sigma}$ be the
corresponding Coxeter cell (i.e., $\cc_{\sigma}$ contains $e$, and
$\nabla_{\cp_0}(\cc_{\sigma})$ is the vertex set of $\sigma$).  Let
$\sigma'$ be the simplex in $\star(\nabla(\bar{e}))$ with vertex set
$\nabla_{\cp_1}(\cc_{\sigma})$.  Then $j_e$ is defined to be the
simplicial map that takes $\sigma$ to $\sigma'$ via
$j_{(\cp_0,\cp_1,\cc_{\sigma})}$.  We note that the maps $j_e$ and
$j_{\bar{e}}$ are inverses of one another; hence, $j_e$ is an
isomorphism, called the {\em gluing isomorphism} associated to $e$.  

\begin{remark}
The definition of $j_e$ depends only on the values of $\nabla$ on the
edges of the $2$-cells that contain $e$.
\end{remark}

\begin{remark}\label{rem:gluing-composition}
Suppose that $\be=(e_1,\ldots,e_n)$ is an edge path in a Coxeter cell
$\cc$ in $\ccc$ with $i(\be)=\cp$ and $t(\be)=\cp'$.  Then
$j_{e_n}\circ\cdots\circ j_{e_1}$ maps $\nabla_{\cp}(\cc)$ to
$\nabla_{\cp'}(\cc)$ and $j_{e_n}\circ\cdots\circ
j_{e_1}|_{\nabla_{\cp}(\cc)}=j_{(\cp,\cp',\cc)}$.  That is to say, the
composition of gluing isomorphisms corresponding to an edge path in
$\cc$ depends only on the endpoints of the path.  In particular, suppose
that $F$ is a $2$-dimensional face of $\cc$ and that
$\be=(e_1,\ldots,e_{2m})$ is an oriented edge path around the boundary
of $F$.  Then $j_{e_{2m}}\circ\cdots\circ j_{e_1}$ is the identity
map on $\nabla_{\cp}(\cc)$.  This last condition
can be rephrased as follows: if $\sigma$ is the $1$-cell in $L$ with
vertex set $\nabla_{\cp}(F)$ and if $\star(\sigma)$ denotes its star
in $L$, then $j_{e_{2m}}\circ\cdots\circ j_{e_1}$ maps $\star(\sigma)$
to itself and $j_{e_{2m}}\circ\cdots\circ
j_{e_1}|_{\star(\sigma)}=\Id_{\star(\sigma)}$.  

Suppose $(e_1,\ldots,e_n)$ is an edge path in $\cc$ and that
$\cp_k=t(e_k)$.  Then $\cp_k=w_k\cp$ for some unique element $w_k\in
W_{\cc}$.  Since $\nabla_{\cp}(\cc)$ is a fundamental set of
generators for $W_{\cc}$ we get a corresponding word $(v_1,\ldots,
v_n)$ in $\nabla_{\cp}(\cc)$ defined by $w_k=v_k\cdots v_1$.  How do
we express the $\nabla(e_k)$ in terms of the $v_k$?  The answer is
simple: $\nabla(e_1)=v_1$ and for $1\leq k<n$, 
$\nabla(e_{k+1})=j_{e_{k}}\circ\cdots\circ j_{e_1}(v_{k+1})$.
\end{remark}

\begin{remark}
If $\ccc$ is a Coxeter cell complex of reflection type (i.e.,
isomorphic to a reflection tiling) with the obvious symmetric framing,
then each gluing automorphism $j_e$ is the identity map on
$\star(\nabla(e))$. 
\end{remark}

The following lemma gives a necessary condition for two framable
Coxeter cell complexes to be isomorphic (not necessarily by a framed
isomorphism).  

\begin{lemma}\label{lem:iso-nec-cond}
Suppose that $(\ccc,\nabla)$ and $(\ccc',\nabla')$ are two framed
Coxeter cell complexes with the same framing system $(V,M,L)$.  Let
$f:\ccc\rightarrow\ccc'$ be a covering projection that induces the map $\phi$
at a given vertex $\cp$ (i.e., $f_{\cp}=\phi:V\rightarrow V$). 
For each $e\in O_{\cp}$, let
$j_e:\star(\nabla(e))\rightarrow\star(\nabla(\bar{e}))$ and 
$j_{f(e)}:\star(\nabla'(f(e)))\rightarrow\star(\nabla'(\overline{f(e)}))$
be the corresponding gluing isomorphisms of $\ccc$ and $\ccc'$,
respectively.  Set
$J_{e}=j_{f(e)}\circ\phi\circ(j_e)^{-1}:\star(\nabla(\bar{e}))\rightarrow
\star(\nabla'(\overline{f(e)}))$.  Then $J_e$ extends to an
automorphism $\tilde{J}_e:V\rightarrow V$ of framing systems.
\end{lemma}

\begin{proof}{}
Let $\bar{\cp}$ denote $t(e)$.  Then the
following diagram of isomorphisms clearly commutes:
\[\xymatrix{\star(\nabla(e))\ar[d]^{j_{e}}\ar[r]^{f_{\cp}=\phi} &
\star(\nabla'(f(e)))\ar[d]^{j_{f(e)}}\\
\star(\nabla(\bar{e}))\ar[r]^{f_{\bar{\cp}}}&
\star(\nabla'(\overline{f(e)}))}\] 
Hence, $\tilde{J}_e=f_{\bar{\cp}}$ is the desired extension of
$J_e$. 
\end{proof} 

\subsection{Homogeneous framings}\label{ss:homogeneous-framings}
\begin{definition}
A framing $\nabla$ on $\ccc$ is {\em homogeneous} if for each $e\in
OE(\ccc)$, the gluing isomorphism
$j_e:\star(\nabla(e))\rightarrow\star(\nabla(\bar{e}))$ depends only on the
value of $\nabla(e)$ (in other words, whenever $\nabla(e)=\nabla(e')$,
we have $\nabla(\bar{e})=\nabla(\bar{e}')$ and $j_e=j_{e'}$).
\end{definition}

For homogeneous framings we shall often write $j_{\nabla(e)}$ instead
of $j_e$.

\begin{remark}
If $\ccc$ is homogeneously framed, then the involution
$\iota_{\cp}:V\rightarrow V$ defined in Remark~\ref{rem:iota} is
independent of $\cp$.  We shall denote this involution by
$v\mapsto\bar{v}$.
\end{remark}

Suppose that $(\ccc,\nabla)$ is homogeneously framed.  We shall now
explain the consequences of homogeneity for the $F_V$-action on
$\ccc^{(0)}$.  Fix a vertex $\cp\in\ccc^{(0)}$.  For each $1$-simplex
$\sigma$ in $L$, let $\cc_{\sigma}$ denote the corresponding Coxeter
$2$-cell at $\cp$ (i.e., $\nabla_{\cp}(\cc_{\sigma})$ is the vertex
set of $\sigma$).  Let $(e_1,\ldots,e_{2m})$ be the edge cycle that
starts at $\cp$ and goes once around $\cc_{\sigma}$, and let
$\nabla(e_1)\cdots\nabla(e_{2m})$ be the corresponding word in $V$.
Let $r_{\sigma}$ denote the image of $\nabla(e_1)\cdots\nabla(e_{2m})$
in $F_V$.  Finally, let $A_{\nabla}$ be the quotient of $F_V$ by the
normal subgroup $N$ generated by $\{v\bar{v}\;|\; v\in V\}\cup
\{r_{\sigma}\;|\;\sigma\;\mbox{a $1$-simplex of $L$}\}$.  In other words, $A_{\nabla}$ is the group
with one generator $\alpha_{v}$ for each element of $V$, and relations
$(\alpha_v)^{-1}=\alpha_{\bar{v}}$ as well as relations corresponding to
$r_{\sigma}$ for each $1$-simplex $\sigma$ of $L$.

A different choice of vertex $\cp'$ leads to a relation $r_{\sigma}'$
that is conjugate to $r_{\sigma}$ in $F_V$.  Hence, the choice of
$\cp'$ yields the same normal subgroup $N$ and the same quotient group
$A_{\nabla}$.  Moreover, in the $F_V$-action on $\ccc^{(0)}$, the
subgroup $N$ acts trivially.  Hence, we have a well-defined transitive
action (from the right) of $A_{\nabla}$ on $\ccc^{(0)}$.

\begin{remark}
Something significant has been gained.  Two reduced edge paths $\be$
and $\be'$ with the same endpoints are homotopic rel endpoints if and
only if one can be pushed to the other across $2$-cells.  Therefore,
the element $a_{\be}$ of $A_{\nabla}$ corresponding to $\be$ depends
only on the homotopy class of $\be$.  
\end{remark}

This remark has the following immediate consequence.

\begin{proposition}
Suppose $(\ccc,\nabla)$ is a homogeneously framed Coxeter cell complex
and that $A_{\nabla}$ is the group defined above.  Then the isotropy
subgroup of the $A_{\nabla}$-action on $\ccc^{(0)}$ at
$\cp\in\ccc^{(0)}$ is naturally identified with $\pi_1(\ccc,\cp)$.
\end{proposition}

As a consequence of Lemma~\ref{lem:frame-fv-maps} we get the following.

\begin{proposition}\label{prop:can-iso}
Suppose $(\ccc,\nabla)$ and $(\ccc',\nabla')$ are two homogeneously
framed Coxeter cell complexes with the same framing system and the same set
of gluing isomorphisms $\{j_v\}_{v\in V}$.  Then the groups
$A_{\nabla}$ and $A_{\nabla'}$ are canonically isomorphic.  Moreover,
given vertices $\cp\in{\ccc}^{(0)}$ and $\cp'\in{\ccc'}^{(0)}$ there is
a framed map $f:(\ccc,\nabla)\rightarrow (\ccc',\nabla')$ taking
$\cp$ to $\cp'$ if and only if $\pi_1(\ccc,\cp)$ is a subgroup of
$\pi_1(\ccc',\cp')$.
\end{proposition}

\begin{corollary}\label{cor:symmetric-homogeneous} 
Suppose $(\ccc,\nabla)$ is a homogeneously framed Coxeter cell
complex.  Let $A=A_{\nabla}$ and $\pi=\pi_1(\ccc,\cp)$.  Then the
following statements are true.
\begin{enumerate}
\item[(1)] $\Aut(\ccc,\nabla)\cong N_A(\pi)/\pi$, where $N_A(\pi)$ denotes
the normalizer of $\pi$ in $A$.
\item[(2)] $\nabla$ is symmetric if and only if $\pi$ is normal in
$A$.
\item[(3)] If $\ccc$ is simply-connected, then $\nabla$ is symmetric
and $\Aut(\ccc,\nabla)\cong A$.
\end{enumerate}
\end{corollary}

\begin{remark}
If the framing is homogeneous and if the involution on $V$ is trivial,
then each gluing map $j_v:\star(v)\rightarrow\star(v)$ is an
involution. 
\end{remark}

Suppose $(\ccc,\nabla)$ is homogeneously framed.  Let
\[\cF_{\ccc}=\{v\in V\;|\;\mbox{$j_v$ does not extend to an automorphism of
$(V,M,L)$}\}.\]
Also, let $t_{\ccc}=\card(\cF_{\ccc})$.  The next
lemma, which is a special case of Lemma~\ref{lem:iso-nec-cond}, gives
a necessary condition for two homogeneously framed Coxeter cell
complexes to be isomorphic.  It will be used in
Section~\ref{s:associahedral-tilings} to distinguish among various
associahedral tilings.

\begin{lemma}\label{lem:t-condition}
Suppose that $(\ccc,\nabla)$ and $(\ccc',\nabla')$ are homogeneously
framed with the same framing system $(V,M,L)$.  Suppose further that
there is a covering projection $f:\ccc\rightarrow\ccc'$.  Then there
is an automorphism $\phi\in\Aut(V,M,L)$ that takes $\cF_{\ccc}$ to
$\cF_{\ccc'}$.  In particular, $t_{\ccc}=t_{\ccc'}$.
\end{lemma}

\begin{proof}{}
Let $\phi$ be the automorphism induced by $f$ at a given vertex $\cp$.
Put $J_v=j_{\phi(v)}\circ\phi\circ(j_v)^{-1}$.  By
Lemma~\ref{lem:iso-nec-cond}, $J_v$ extends.  Therefore, $j_v$ extends
if and only if $j_{\phi(v)}$ extends. 
\end{proof}   

\subsection{Maximally symmetric tilings}

\begin{definition}
A proper action of a discrete group $G$ on a space $Y$ is {\em rigid}
if for all $g\in G-\{1\}$, there does not exist a nonempty open subset
$U$ of $Y$ that is fixed pointwise by $g$.
\end{definition}

For example, if $Y$ is a connected manifold, then it follows from
Newman's Theorem (page 153 in \cite{Bre}) that any proper effective
action on $Y$ is rigid.

Suppose $(X,\nabla)$ is a homogeneously framed Coxeter cell complex
with framing system $(V,M,L)$.  Let $\Aut(X)$ denote the automorphism
group of the cell complex $X$, and for $x$ in $X^{(0)}$, let
$\Stab(x)$ be the stabilizer subgroup of $x$.  Assume
further that the action of $\Aut(X)$ on $X$ is rigid.  Let
$\Aut_0(V,M,L)$ be the subgroup of all $\phi\in\Aut(V,M,L)$ that are   
equivariant with respect to the involution $v\mapsto\bar{v}$.  Then
the natural homomorphism $\Phi:\Stab(x)\rightarrow\Aut_0(V,M,L)$ is
injective.  If, in addition, $\Phi$ is surjective we say that $X$ is
{\em maximally symmetric}.  

\begin{proposition}\label{prop:max-symm}
Suppose $X$ is simply connected.  Then $X$ is maximally symmetric if
and only if $\phi\circ j_v\circ\phi^{-1}\circ(j_{\phi(v)})^{-1}$
extends to an automorphism of framing systems for all $v\in V$ and 
$\phi\in\Aut_0(V,M,L)$, 
\end{proposition}

\begin{proof}{}
If $X$ is maximally symmetric, then the condition follows from
Lemma~\ref{lem:iso-nec-cond}.  Let $\nabla'$ be the framing of $X$
defined by $\nabla'=\nabla^{\phi^{-1}}$.  Then $(X,\nabla'$) is
homogeneously framed and the corresponding gluing isomorphisms $j_v'$
are given by $j_v'=\phi^{-1}\circ j_{\phi(v)}\circ\phi$.  If the
condition holds, then $\{j_v'\}_{v\in V}=\{j_v\}_{v\in V}$.  So, by
Proposition~\ref{prop:can-iso}, $A_{\nabla}$ is canonically isomorphic
to $A_{\nabla'}$ (the isomorphism is induced by $v\mapsto\phi(v)$),
and since $X$ is simply connected there is a framed isomorphism
$f:(X,\nabla)\rightarrow(X,\nabla')$ taking $x$ to $x$.  Since
$\Phi(f)=\phi$, $\Phi$ is surjective.
\end{proof}    

\begin{example}\label{ex:max-symm-refl}
If $\Aut(V,M,L)$ is trivial, then $X$ is maximally symmetric and
$\Aut(X)=A_{\nabla}$.  If $X=\Sigma(W,S,L)$ is a tiling of reflection
type, then $X$ is maximally symmetric (each $j_v$ is trivial).  In
fact, in this case, the full symmetry group of $X$ 
is 
\[\Aut(X)=W\rtimes\Aut(V,M,L).\]
\end{example}

\subsection{A sufficient condition for two Coxeter cell complexes to
be isomorphic}\label{ss:iso-reflection}
Suppose that $(\ccc,\nabla)$ is a framed Coxeter cell
complex with framing system $(V,M,L)$.  Let $\Sigma=\Sigma(W,M,L)$ be
the corresponding Coxeter cell complex of reflection type (where
$W=W(M)$).  In this section we give sufficient conditions for finding
a covering projection $p:\Sigma\rightarrow X$ such that $p$ takes the
distinguished vertex $1\in\Sigma^{(0)}$ to a given vertex $x\in
X^{(0)}$ so that the induced map of framing systems $p_x:V\rightarrow
V$ at the vertex $1$ is the identity.  The two conditions we shall
give are necessary by Lemma~\ref{lem:iso-nec-cond} and
Remark~\ref{rem:gluing-composition}.  The first condition is the
following ``Extendability Condition'': 

\begin{enumerate}
\item[(E)] For each $e\in OE(\ccc)$, the map
$j_{e}:\star(\nabla(e))\rightarrow \star(\nabla(\bar{e}))$ 
extends to an automorphism $\tilde{j}_e:V\rightarrow V$ of framing
systems. 
\end{enumerate}

\noindent The second condition (H) is the condition of ``no holonomy around
$2$-cells''.  It is justified by the first paragraph of
Remark~\ref{rem:gluing-composition}.

\begin{enumerate}
\item[(H)] The extended gluing isomorphisms of (E) can be chosen so
that $\tilde{j}_{\bar{e}}=(\tilde{j}_e)^{-1}$ and so that if
$(e_1,\ldots,e_{2m})$ is a cycle around a $2$-cell of $X$, then
$\tilde{j}_{e_{2m}}\circ\cdots\circ\tilde{j}_{e_1}$ is the identity
map.
\end{enumerate}

\begin{remark}
Suppose $(e_1,\ldots,e_{2m})$ is a cycle around a $2$-cell $F\subset
X$, that $y=i(e_1)=t(e_{2m})$, and that $\sigma$ is the $1$-cell of
$L$ spanned by $\nabla_y(e_1)$ and $\nabla_y(\bar{e}_{2m})$.  Then
Remark~\ref{rem:gluing-composition} asserts that
$j_{e_{2m}}\circ\cdots\circ j_{e_1}|_{\star(\sigma)}$ is the identity
map.  The content of Condition (H) is that
$\tilde{j}_{e_{2m}}\circ\cdots\circ\tilde{j}_{e_{1}}$ is an extension
of this to the identity map of $V$.
\end{remark}

\begin{remark} \label{rem:rigid}
If the action of $\Aut(L)$ on $L$ is rigid, then the extension
$\tilde{j}_e$ of Condition (E) is unique (if it exists); moreover,
Condition (H) is then automatic.
\end{remark}

Now suppose that $(X,\nabla)$ satisfies (E) and (H).  Given a word
$(v_1,\ldots,v_n)$ in $V$ we will define an edge path
$\be=(e_1,\ldots,e_n)$ beginning at $x=x_0$.  We will use the notation
$x_k=t(e_k)$.  By definition $e_1$ is the unique edge at $x_0$ such
that $\nabla_{x_0}(e_1)=v_1$.  Suppose by induction that $(e_1,\ldots
e_k)$ has been defined.  Then $e_{k+1}$ is the unique edge at $x_k$
satisfying
$\nabla_{x_k}(e_k)=\tilde{j}_{e_k}\circ\cdots\circ\tilde{j}_{e_1}(v_{k+1})$.
(Compare this to the formula in the second paragraph of
Remark~\ref{rem:gluing-composition}.)

Set $w_k=v_1\cdots v_k$.

\begin{lemma}\label{lem:W-action}
Suppose Conditions (E) and (H) hold.  Then with notation as above, the
vertex $x_n\in X^{(0)}$ depends only on the element $w=w_n\in W$ (and
not on the choice of word $v_1\cdots v_n$ representing $w$).
Moreover, this gives a right action of $W$ on $X^{(0)}$ defined by
$x\cdot w=x_n$.
\end{lemma}

\begin{proof}{}
Any two words for an element $w$ in a Coxeter group differ by a
sequence of moves of one of the following two types: (1) cancelling a
subword of the form $vv$ or (2) replacing a subword which goes half
way around a $2$-cell in $\Sigma$ by the subword which goes half way
around the $2$-cell in the other direction.  A subword of the form
$vv$ in $(v_1,\ldots,v_n)$ corresponds to an edge subpath of the form
$e\bar{e}$ and since $\tilde{j}_e\circ\tilde{j}_{\bar{e}}=\Id$, we can
cancel.  So the issue comes down to the following.  Suppose
$(v_1,\ldots,v_n)$ contains a subword $\bs=(s,t,\ldots)$ which is an
alternating word of length $m=m(s,t)$ in letters $s$ and $t$.  We want
to show that if we replace this subword by the other alternating word
$\bs'=(t,s,\ldots)$ of length $m$ then the initial and final segments
of the new edge path in $X$ remain unchanged.  Suppose
$(e_1,\ldots,e_m)$ is the portion of the edge path corresponding to
$\bs$.  Let $(\bar{e}_{2m},\ldots,\bar{e}_{m+1})$ be the edge path with
the same initial vertex that corresponds to $\bs'$.  Its opposite path
is $(e_{m+1},\ldots,e_{2m})$ and $(e_1,\ldots,e_{2m})$ is a cycle around
a $2$-cell in $X$.  By Condition (H)
\[\tilde{j}_{e_m}\circ\cdots\circ\tilde{j}_{e_1}=
\tilde{j}_{\bar{e}_{m+1}}\circ\cdots\circ\tilde{j}_{\bar{e}_{2m}}.\]
Hence, we can replace the segment $(e_1,\ldots,e_m)$ by
$(\bar{e}_{2m},\ldots,\bar{e}_{m+1})$ without affecting the definition of any
succeeding edges.
\end{proof}

Define $p^{(0)}:\Sigma^{(0)}\rightarrow{\ccc}^{(0)}$ by $p^{(0)}(1\cdot
w)=\cp\cdot w$.  It follows as in Lemma~\ref{lem:frame-fv-maps} that
this extends to a covering projection $p:\Sigma\rightarrow\ccc$.  We
have proved the following.

\begin{theorem}\label{thm:iso-to-reflection}
Suppose that Conditions (E) and (H) hold for a framed Coxeter cell
complex $(\ccc,\nabla)$ and that $\Sigma$ is the associated tiling of 
reflection type.  Then for any vertex $\cp\in{\ccc}^{(0)}$, there is
a covering projection $p:\Sigma\rightarrow\ccc$ taking $1$ to $\cp$.
\end{theorem}

\subsection{The completed complex $\hat{\ccc}$}\label{ss:hat-X}
Suppose $(\ccc,\nabla)$ is a framed Coxeter cell complex with framing
system $(V,M,L)$.  Let $\hat{L}$ denote the simplicial complex
$\cS(V)_{>\emptyset}$.  Our goal in this section is to try to
construct a new Coxeter cell complex $\hat{\ccc}$ such that
\begin{enumerate} 
\item[(a)] $\ccc$ is a subcomplex of $\hat{\ccc}$ and they have the
same $2$-skeleton, and 
\item[(b)] for each vertex $\cp\in\hat{\ccc}^{(0)}$ ($=\ccc^{(0)}$)
its link $\hat{L}_{\cp}$, is isomorphic to $\hat{L}$ and the
framing gives an isomorphism
$\nabla_{\cp}:\hat{L}_{\cp}\rightarrow\hat{L}$. 
\end{enumerate}
If we are successful in this construction then, since
$OE(\hat{\ccc})=OE(\ccc)$, we will have an induced framing on
$\hat{\ccc}$ with framing system $(V,M,\hat{L})$.

There are two obvious conditions for (b) to hold: first, each gluing
isomorphism must extend to an isomorphism of the appropriate stars in
$\hat{L}$ and second, there can be no holonomy around $2$-cells in the
star of the appropriate edge of $\hat{L}$.  For each subset $T$ of
$V$, let $M(T)$ denote the restriction of the Coxeter matrix to $T$.
For each $v\in V$, let $V_v$ denote the vertex set of $\star(v)$.  The
gluing isomorphism
$j_e:\star(\nabla(e))\rightarrow\star(\nabla(\bar{e}))$ restricts to a
map of vertex sets which we also denote by
$j_e:V_{\nabla(e)}\rightarrow V_{\nabla(\bar{e})}$.  We next discuss a
Condition (M) (for ``Matrix Condition'').  It has two parts (M1) and
(M2).  The first part (M1) is a weak version of (E) in
\ref{ss:iso-reflection}, and the second (M2) is a weak version of (H).  

\begin{enumerate}
\item[(M1)] For each $e\in OE(\ccc)$ and for each pair $v_1$, $v_2$ in
$V_{\nabla(e)}-\{\nabla(e)\}$,
$m(j_e(v_1),j_e(v_2))=m(v_1,v_2)$. 
\end{enumerate}
Condition (M1) says that $j_e$ pulls back $M(V_{\nabla(\bar{e})})$ to
$M(V_{\nabla(e)})$.  Hence, if (M1) holds, then each $j_e$
induces a simplicial isomorphism
\[\hat{j}_e:\cS(V)_{\geq\{\nabla(e)\}}\rightarrow
\cS(V)_{\geq\{\nabla(\bar{e})\}}\]
between the appropriate stars in $\hat{L}$.

Note that Condition (M1) is automatic if $\{\nabla(e),v_1,v_2\}$
spans a $2$-simplex in $\star(\nabla(e))$ (since $j_e$ maps this
$2$-simplex to a $2$-simplex in $\star(\nabla(\bar{e}))$).  Thus,
Condition (M1) holds automatically if $L$ has no ``missing $2$-simplices''. 

Suppose $F$ is a $2$-cell of $X$, that $(e_1,\ldots,e_{2m})$ is a
cycle around $F$, and that $x=i(e_1)$.  The second part of Condition
(M) is the following:
\begin{enumerate}
\item[(M2)] For any $2$-cell $F$ of $X$, with notation as above, 
\[\hat{j}_{e_{2m}}\circ\cdots\circ\hat{j}_{e_1}:
\cS_{\geq\nabla_x(F)}\rightarrow\cS_{\geq\nabla_x(F)}\]
is the identity map.
\end{enumerate}

We note that Condition (M2) is automatic if there are no missing
simplices in the star of the edge $\sigma$ of $L$ which corresponds to
$\nabla_x(F)$.  Thus, Condition (M) holds automatically if $L$ has no
missing simplices.

Supposing that (M) holds, the idea for the construction of
$\hat{\ccc}$ is the following.  We want to fill in the ``missing''
cells of $\ccc$ that correspond to the ``missing'' simplicies of
$L$.  The first problem is to describe subcomplexes of $\ccc$ which
will serve as the (partial) boundaries of the missing cells to be
filled in.  This is done essentially by the same method as in
\ref{ss:iso-reflection}.  Condition (M) is exactly the hypothesis
needed to carry out this procedure.  A second problem then arises.
Our putative boundaries of missing cells might not be embedded
subcomplexes of $\ccc$.  Instead, they might be the image of the
actual boundary of a cell under a covering projection.  Thus, the
completion $\hat{\ccc}$ might turn out to be a ``Coxeter orbihedron''
(as defined below).

For each $T\in\cS(V)_{>\emptyset}$ let $L_T$ be the subcomplex of
$L$ spanned by $T$.  By condition (iii) of
Definition~\ref{def:framing-system}, the $1$-skeleton of $L_T$ is
the complete graph on $T$.  Let $W_T$ be the finite Coxeter group
corresponding to $M(T)$, let $\cc$ be the corresponding Coxeter cell,
and let $1\in\cc^{(0)}$ be a distinguished vertex (such that
$OE_1(\cc)=T$).  Let $\cc(L_T)$ be the subcomplex of $\cc$
corresponding to $L_T$, i.e., $Z(L_T)=\Sigma(W_T,T,L_T)$.

Fix a vertex $\cp\in\ccc^{(0)}$ and let $T$ be any spherical subset of
$\nabla_{\cp}(OE_{\cp})$.  We now propose to define a subcomplex
$C(L_T)$ of $\ccc$ containing the vertex $\cp$, and a map
$\phi:\cc(L_T)\rightarrow C(L_T)$ such that $\phi(1)=\cp$.  (If $T$ is
the vertex set of a simplex $\sigma$ in $L$, then $C(L_T)$ will be
$\cc_{\sigma}$, the corresponding Coxeter cell at $\cp$.)

We regard $T$ as a set of fundamental generators for $W_T$.  We
proceed as in \ref{ss:iso-reflection}.  Let $\bu=u_1\cdots u_n$ be a
word in $T$.  For $1\leq k\leq n$, let $w_k$ be the element of $W_T$
represented by $u_k\cdots u_1$.  Assuming that (M) holds, we are going
to define, by induction on $k$, an edge path $\be=(e_1,\ldots,e_n)$
and a sequence of subsets $T_0,\ldots,T_n$ of $V$.  We will
use the notation: $v_k=\nabla(e_k)$, $\cp_k=t(e_k)$, and $\cp_0=\cp$.
We will also show by induction that $T_k$ is a spherical subset of
$\nabla(OE_{\cp_k})$.  Let $T_0=T$ and let $e_1\in OE_{\cp_0}$ be the
unique element satisfying $\nabla(e_1)=u_1$.  Assuming that $e_i$ and
$T_{i-1}$ have been defined for $1\leq i\leq k$, set
$T_k=\hat{j}_{e_k}(T_{k-1})=\hat{j}_{e_k}\circ\cdots\circ\hat{j}_{e_1}(T)$
and let $v_{k+1}\in T_k$ be the element defined by
$v_{k+1}=\hat{j}_{e_k}\circ\cdots\circ\hat{j}_{e_1}(u_{k+1})$.  (Here
we have used induction and Condition (M) to show that
$\hat{j}_{e_{k-1}}\circ\cdots\circ\hat{j}_{e_1}(u_{k+1})$ lies in the
domain of $\hat{j}_{e_k}$.)  Let $e_{k+1}\in OE_{\cp_k}$ be the unique
element satisfying $\nabla(e_{k+1})=v_{k+1}$.  The proof of
Lemma~\ref{lem:W-action} then gives the following.

\begin{lemma}\label{lem:no-monodromy}
Assume (M) holds for $(\ccc,\nabla)$.  With notation as above, the
vertex $\cp_n\in\ccc^{(0)}$, the spherical subset $T_n$, and the
isomorphism $\hat{j}_{e_n}\circ\cdots\circ\hat{j}_{e_1}:T\rightarrow
T_n$ depend only on the element $w=w_n\in W_T$ (and not on the choice
of word $u_n\cdots u_1$ representing $w$).  Hence, we can use the
notation $w\cdot\cp=\cp_n$, $wT=T_n$, and
$w_*=\hat{j}_{e_n}\circ\cdots\circ\hat{j}_{e_1}:T\rightarrow wT$.
\end{lemma}

For each simplex $\sigma$ of $L_T$ and $w\in W_T$, $wW_{\sigma}\cp$ is
the vertex set of a cell in $\ccc$.  The subcomplex $C(L_T)$ is the
union of these cells.  The map $\phi:\cc(L_T)\rightarrow C(L_T)$ is
defined on vertices by $w\mapsto w\cdot\cp$.  The map $\phi$ is a
covering projection and induces an isomorphism
$\cc(L_T)/W_{T,\cp}\cong C(L_T)$, where $W_{T,\cp}$ denotes the
isotropy subgroup of $W_T$ at $\cp$.

A {\em Coxeter orbicell} is the quotient of a Coxeter cell (of type
$W_T$) by a subgroup of $W_T$.  A {\em Coxeter orbihedron} is an
orbihedron that is locally isomorphic to the product of a cell and
some Coxeter orbicell.  Thus, a Coxeter orbihedron is decomposed into
Coxeter orbicells.  The link of a vertex in a Coxeter orbihedron is
defined as before; it is a simplicial cell complex.

The Coxeter orbihedron $\hat{\ccc}$ is constructed by gluing onto
$\ccc$ a copy of $Z/W_{T,\cp}$ for each subcomplex of the form
$C(L_T)$.  In more detail, one constructs a sequence of Coxeter
orbihedra $\ccc=\ccc_2,\ccc_3,\ldots$ where $\ccc_n$ is supposed to
be the union of $\ccc$ and the $n$-skeleton of $\hat{\ccc}$, and
$\ccc_{n+1}$ is obtained from $\ccc_n$ by gluing in the missing
Coxeter orbicells of dimension $n+1$.  We have proved the
following. 

\begin{theorem}\label{thm:orbihedral-completion}
Suppose $(\ccc,\nabla)$ is a framed Coxeter cell complex and that
Condition (M) holds.  Then $\ccc$ is a subcomplex of a Coxeter
orbihedron $\hat{\ccc}$ with the same $2$-skeleton such that for each
$\cp\in\ccc^{(0)}$, $\nabla_{\cp}$ induces an isomorphism from
$\hat{L}_{\cp}$ to $\cS(V)_{>\emptyset}$.  
\end{theorem}

\begin{lemma}\label{lem:condition-M-for-blowups}
As in Section~\ref{s:blow-ups}, suppose $\Sigma_{\#}$ is the
$\cR$-blow-up of $\Sigma(W,S)$.  Then the natural framing on
$\Sigma_{\#}$ (defined in Example~\ref{ex:blow-up-framing}) satisfies
Condition (M).
\end{lemma}

\begin{proof}{}
An empty $k$-simplex of $L_{\#}$, $k\geq 2$, corresponds to a subset
$\{T_1,\ldots,T_{k+1}\}$ of $S_{\#}$ that is not $\cR$-nested, but
such that any proper subset is.  In particular, $T_i\not\subset T_l$
for all $i\neq l$ (if $T_i\subset T_l$ then adding $T_i$ to the
$\cR$-nested set $\{T_1,\ldots,\hat{T_i},\ldots,T_k\}$ would give an
$\cR$-nested set).  Thus, if $T$ is one of these vertices, say 
$T_i$, then it follows from the definition of $j_T$ in
\ref{ss:blow-up-gluing} that $j_T(T_l)=T_l$ for $l\neq i$.  In
other words, $j_T$ acts trivially on any empty $k$-simplex in
$\cN_{\geq\{T\}}$.  This implies that Condition (M) holds.
\end{proof}

\begin{example}\label{ex:n-cube-sharp}
Suppose that $\cc_{\#}$ is the blow-up of the $n$-cube at its center.
Thus, $\cc_{\#}$ is an interval bundle over $\bbR\bbP^{n-1}$
($=\partial\cc_{\#}/a$).  Let $(S_{\#},M_{\#},L_{\#})$ be the framing
system for the natural framing on $\cc_{\#}$ described in
Example~\ref{ex:blow-up-framing}.  If $S=\{s_1,\ldots,s_n\}$, then
$\cR=\{S\}$ and $S_{\#}=\{\{s_1\},\ldots,\{s_n\},S\}$.  The link
$L_{\#}$ is the subdivision of the $(n-1)$-simplex with vertex set
$T=\{\{s_1\},\ldots\{s_n\}\}$ obtained by introducing a barycenter
$S$.  The picture for $n=3$ is given in Figure~\ref{fig:completion}.  

\begin{figure}[ht]
\begin{center}
\psfrag{a}{\scriptsize $\{s_1\}$}
\psfrag{b}{\scriptsize $\{s_2\}$}
\psfrag{c}{\scriptsize $\{s_3\}$}
\psfrag{S}{\scriptsize $S$}
\hspace*{-.2in}\includegraphics[scale = .35]{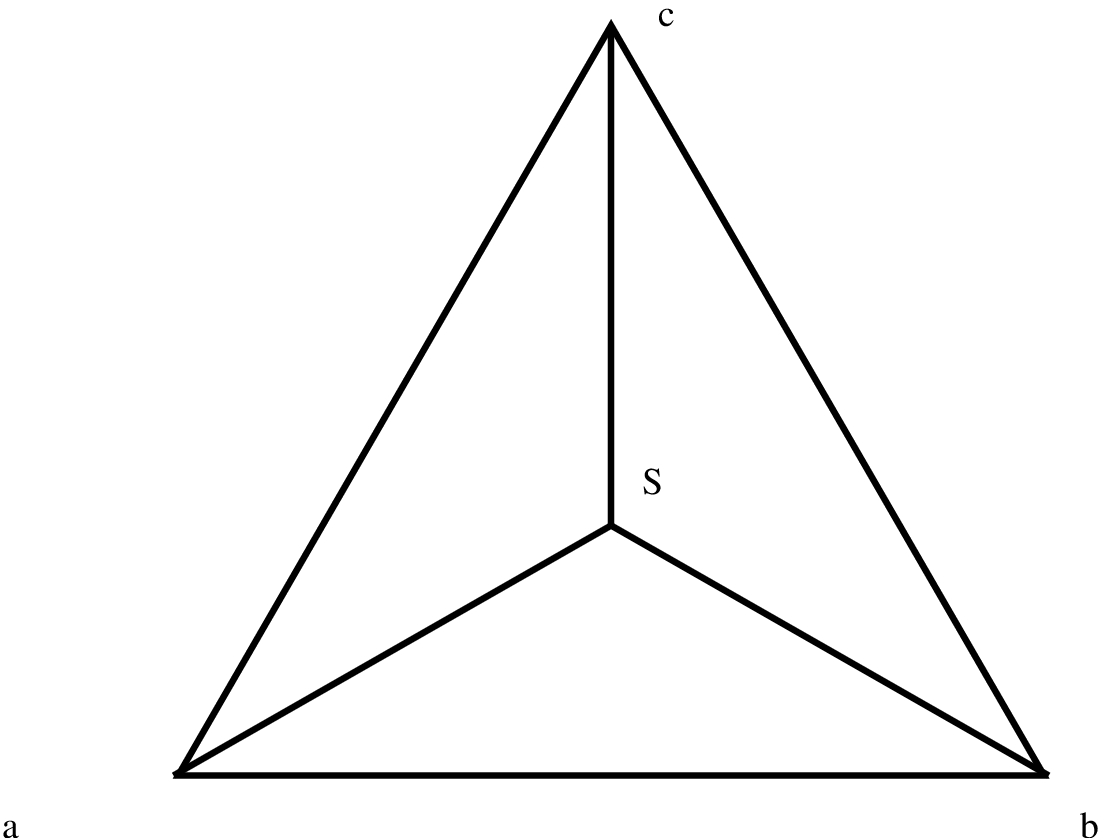}
\caption{\label{fig:completion}}
\end{center}
\end{figure}

The complex $\hat{L}=\cS(S_{\#})_{\geq\emptyset}$ is the full
$n$-simplex on $S_{\#}$.  To describe the orbihedron $\hat{\cc}_{\#}$,
we note that there are two missing simplices in $L_{\#}$.  The first
is the $(n-1)$-simplex spanned by $T$, and the second is the full
$n$-simplex on $S_{\#}$.  Clearly,
$C(L_T)=\partial\cc_{\#}=\partial([-1,1]^n)$.  We fill in the missing
$n$-cell to obtain 
\[\ccc_n=\cc_{\#}\cup_{\partial\cc_{\#}}[-1,1]^n\]
which is isomorphic to $\bbR\bbP^n$ ($=\partial([-1,1]^{n+1})/a$).
The missing $(n+1)$-cell is $[-1,1]^{n+1}$, but its boundary is no
longer embedded (it is the complex $\ccc_n$).  Thus,
\[\hat{\cc}_{\#}=[-1,1]^{n+1}/a.\]
\end{example}

\begin{example}\label{ex:n-cube-sharp-tilde}
Suppose $Y$ is the universal cover of the complex
$\cc_{\#}$ in the previous example.  Thus, $Y$ is isomorphic to
$S^{n-1}\times[-1,1]$ ($=\partial([-1,1]^n)\times[-1,1]$).  This time
there are two missing $n$-cells with boundaries $S^{n-1}\times\{-1\}$
and $S^{n-1}\times\{1\}$.  We fill these in to obtain
$Y_n=\partial([-1,1]^{n+1})$.  Filling in the missing $(n+1)$-cell
we obtain $\hat{Y}=[-1,1]^{n+1}$.
\end{example}

\subsection{Metric flag complexes}
The natural piecewise Euclidean metric on a Coxeter cell complex
induces a piecewise spherical metric on the link of each vertex.
Results of Gromov \cite{G} and Moussong \cite{M} show that the
condition that a Coxeter cell complex be nonpositively curved in the
sense of Aleksandrov and Gromov is equivalent to a condition on the
link of each vertex.  We shall now recall this condition.  (For more
details see \cite{BH} or \cite{DM}.)

A spherical simplex {\em has size $\geq\pi/2$} if the length of each
of its edges is at least $\pi/2$.  Similarly, a piecewise spherical
simplicial cell complex {\em has size $\geq\pi/2$} if each of its
simplices does.

\begin{example}
The link of a vertex in a Coxeter cell of type $W_T$ is a spherical
simplex $\sigma$ with vertex set $T$.  The length of the edge
connecting the vertices $t_1$ and $t_2$ is $\pi-\pi/m(t_1,t_2)$, which
is $\geq\pi/2$.  Hence $\sigma$ has size $\geq\pi/2$.  It follows that
the link of any vertex in a Coxeter cell complex has size
$\geq\pi/2$.  Similarly, if $(V,M,L)$ is any framing system, then the
natural piecewise spherical structure on $L$ (in which the length of
an edge from $v_1$ to $v_2$ is given by $\pi-\pi/m(v_1,v_2)$) has size
$\geq\pi/2$. 
\end{example}

The definition of {\em flag complex} in \ref{ss:flag-complex} can be
generalized as follows.

\begin{definition}
Suppose $L$ is a piecewise spherical simplicial cell complex of size
$\geq\pi/2$.  Then $L$ is a {\em metric flag complex} if it is a
simplicial complex and if given any nonempty finite collection $T$ of
vertices of $L$ such that (a) any two elements of $T$ are connected by
an edge in $L$ and (b) it is possible to find a spherical simplex with
the same set of edge lengths, then $T$ spans a simplex in $L$.
\end{definition} 

\begin{example}
With regard to condition (b), if the length of the edge connecting any
two vertices $t_1,t_2\in T$ is given by a Coxeter matrix (i.e., if it
is $\pi-\pi/m(t_1,t_2)$), then there exists a spherical simplex with
this set of edge lengths if and only if the associated Coxeter group
$W_T$ is finite (see Example~6.7.3 in \cite{DM}).  It follows that if
$(V,M,L)$ is a framing system, then $L$ is a metric flag complex if
and only if $L=\hat{L}$ (where $\hat{L}=\cS(V)_{>\emptyset}$).
\end{example}

\begin{remark}
If $L$ has size $\geq\pi/2$ and if it is a flag complex, then it is a
metric flag complex.  Conversely, when the length of each edge is
$\pi/2$, then $L$ is a metric flag complex only if it is a flag
complex. 
\end{remark}

Moussong generalized Gromov's criterion for nonpositive curvature of a
cubical complex (stated in \ref{ss:flag-complex}) by showing that a
Coxeter cell complex is nonpositively curved if and only if the link
of each vertex is a metric flag complex.  (Proofs can be found in
\cite{M} or in Section~6.7 of \cite{DM}.)  A corollary of this is the
following. 

\begin{theorem}\label{thm:moussong}
(Moussong) Suppose $\ccc$ is a framed Coxeter cell complex with
framing system $(V,M,L)$.  Then $\ccc$ is nonpositively curved if and
only if $L=\hat{L}$.
\end{theorem}

\subsection{$\hat{\ccc}$ is nonpositively curved}\label{ss:cat-0}
In this subsection $(\ccc,\nabla)$ is a framed Coxeter cell complex
satisfying Condition (M) of \ref{ss:hat-X} and $\hat{\ccc}$ is the
Coxeter orbihedron constructed in
Theorem~\ref{thm:orbihedral-completion}.  It follows from Moussong's
result (Thoerem~\ref{thm:moussong}) that $\hat{\ccc}$ is nonpositively
curved. 

\begin{theorem}
If $\ccc$ is simply connected, then $\hat{\ccc}$ is actually a Coxeter
cell complex.
\end{theorem}

Some evidence for this theorem is provided by
Examples~\ref{ex:n-cube-sharp} and \ref{ex:n-cube-sharp-tilde}.

\begin{proof}{}
Nonpositively curved orbihedra are ``developable'' (see page 562 in
\cite{BH}).  This means that the local isotropy group at any point in
the universal orbihedral cover is trivial.  Thus, the universal
orbihedral cover of any nonpositively curved Coxeter orbihedron is a
Coxeter cell complex.  The $2$-skeleton of $\hat{\ccc}$ (as an
orbicell complex) is the same as that of $\ccc$.  So, if $\ccc$ is
simply connected, then so is $\hat{\ccc}$ and hence, $\hat{\ccc}$ is
its own orbihedral cover.
\end{proof}

If a complete metric space is nonpositively curved and
simply connected, then it is a {\em CAT($0$)-space} (meaning that Gromov's
CAT($0$)-inequalities hold for all triangles in the space).  This
implies, for instance, that the space is contractible.  Let us say
that a group $G$ is a {\em CAT($0$)-group} if it admits a representation
as a discrete, cocompact group of isometries on a finite dimensional
CAT($0$)-space.  In the following theorem we list some well-known
properties of CAT($0$)-groups.  (Proofs can be found in \cite{BH}.)

\begin{theorem}\label{thm:CAT(0)-groups}
Suppose $G$ is a discrete, cocompact group of isometries of a CAT($0$)
polyhedron $K$.  Then the following statements are true.
\begin{enumerate}
\item[(1)] $G$ is finitely presented.
\item[(2)] The word problem and the conjugacy problem are solvable for
$G$.
\item[(3)] $G$ has only finitely many congugacy classes of finite
subgroups. 
\item[(4)] If $H$ is any solvable subgroup of $G$, then $H$ is
virtually abelian.
\item[(5)] $H_*(G;\bbQ)$ is finite dimensional and $\operatorname{cd}_{\bbQ}(G)$ is
no greater than $\dim K$.
\item[(6)] If $\pi$ is any torsion-free subgroup of $G$, then $\pi$
acts freely on $K$ and $K/\pi$ is a $K(\pi,1)$-complex for $\pi$.
\end{enumerate}
\end{theorem}

If $G$ is any discrete, cocompact group of automorphisms of the
Coxeter cell complex $\ccc$, then it extends to a cocompact group of
automorphisms of $\hat{\ccc}$.  If, in addition, $\ccc$ is simply
connected, then $\hat{\ccc}$ is CAT($0$) and these automorphisms are
isometries.  Thus, $G$ will be a CAT($0$)-group.  So, as a corollary
to Theorem~\ref{thm:CAT(0)-groups} we have the following.

\begin{theorem}\label{thm:CAT(0)-groups-from-framings}
Suppose that $\ccc$ is a simply-connected, framed Coxeter cell complex
satisfying Condition (M) of \ref{ss:hat-X}.  Suppose further that
$\ccc$ admits a discrete, cocompact group of automorphisms $G$.  Then
the following statements are true. 
\begin{enumerate}
\item[(1)] $G$ is finitely presented.
\item[(2)] The word problem and the conjugacy problem are solvable for
$G$.
\item[(3)] $G$ has only finitely many congugacy classes of finite
subgroups. 
\item[(4)] If $H$ is any solvable subgroup of $G$, then $H$ is
virtually abelian.
\item[(5)] $H_*(G;\bbQ)$ is finite dimensional and
$\operatorname{cd}_{\bbQ}(G)$ is no greater than $\dim\hat{\ccc}$.
\item[(6)] If $\pi$ is any torsion-free subgroup of $G$, then $\pi$
acts freely on $\ccc$.  If, in addition, $\pi$ has finite index in
$G$, then we can add finitely many cells to $\ccc/\pi$ to obtain a
$K(\pi,1)$-complex.
\end{enumerate}
\end{theorem}

\begin{proof}{}
All statements but (6) follow immediately from the previous theorem
applied to the $G$-action on $\hat{\ccc}$.  For (6), since the
$\pi$-action is free on $\hat{\ccc}$, it is also free on $\ccc$.  The
$K(\pi,1)$-complex is then $\hat{\ccc}/\pi$. 
\end{proof} 

From Theorem~\ref{thm:CAT(0)-groups-from-framings} and
Lemma~\ref{lem:condition-M-for-blowups} we immediately get the
following.

\begin{corollary} \label{cor:complete-blow-up}
Suppose that $\tilde{\Sigma}_{\#}$ is
the universal cover of an $\cR$-blow-up of $\Sigma(W,S)$ and that $A$
is its symmetry group as defined in \ref{ss:group-A}.  Then $A$ is a
CAT($0$)-group (and, hence, has all the properties listed in
Theorem~\ref{thm:CAT(0)-groups-from-framings}).
\end{corollary}

\subsection{Construction of the group and the complex from the gluing
isomorphisms}\label{ss:gluing2group}

Suppose we are given a Coxeter matrix $M$ ($=m(u,v)$) on a finite set
$V$ and an involution on $V$ denoted by $v\mapsto\bar{v}$.  Call two
distinct elements $u$ and $v$ of $V$ {\em adjacent} if
$m(u,v)\neq\infty$.  Let $L$ denote the $1$-skeleton of
$\cS(V)_{>\emptyset}$ (i.e., $L$ is the graph with vertex set $V$ and
with an edge connecting any two adjacent vertices).  As before, $V_{v}$
denotes the star of $v$ in $L$ (i.e., $V_v=\{u\in V\;|\;
m(u,v)\neq\infty\}$).  Also suppose that we are given as ``gluing
data'' a set $\{j_v\}_{v\in V}$ of bijections $j_v:V_v\rightarrow
V_{\bar{v}}$ such that $j_v(v)=\bar{v}$.  The question we address here
is: when can we find a group $A$ and a simply connected,
$2$-dimensional, Coxeter cell complex $X$ such that (a) the link of
each vertex in $X$ is $L$, (b) $A$ acts simply transitively on
$X^{(0)}$, and (c) $\{j_v\}_{v\in V}$ is the corresponding set of
gluing isomorphisms?

There are four obvious conditions that should be imposed on the $j_v$.
First, 
\begin{enumerate}
\item[(1)] $j_{\bar{v}}=(j_v)^{-1}$ for all $v\in V$.
\end{enumerate}
The second condition is just a rewording of Condition (M1) from
\ref{ss:hat-X}:
\begin{enumerate}
\item[(2)] For each pair $v_1,v_2$ in $V_v-\{v\}$,
$m(j_v(v_1),j_v(v_2))=m(v_1,v_2)$.
\end{enumerate}
Before stating the third condition we need to develop some more
notation.  Given an ordered pair $(u,v)$ of adjacent elements in $V$,
define a sequence of elements $v_0,v_1,v_2,\ldots$ by the formulas:
\[v_0=\bar{u},\;\; v_1=v,\;\;\mbox{and}\;\;
v_k=j_{v_{k-1}}(\bar{v}_{k-2})=j_{v_{k-1}}j_{v_{k-2}}(v_{k-2})\;\;
\mbox{for $k\geq 2$}\hspace{.5in}(\ast)\]
(Just as in \ref{ss:hat-X}, it follows by induction that
$\bar{v}_{k-2}\in V_{v_{k-1}}$ and hence, that the above formula makes
sense.)  We note that $v_{2l}=j_{v_{2l-1}}\cdots j_{v_1}(u)$ and
$v_{2l+1}=j_{v_{2l}}\cdots j_{v_1}(v)$.  Our third condition is the
following:
\begin{enumerate} 
\item[(3)] Given any ordered pair $(u,v)$ of adjacent elements in $V$,
let $v_0,v_1,\ldots$ be the sequence defined above, and let
$m=m(u,v)$.  Then $v_{2m}=v_0=\bar{u}$ and $v_{2m+1}=v_1=v$.
\end{enumerate}
We note that (3) implies that $v_{2m+k}=v_k$ for all $k\geq 0$.

These first three conditions are enough for us to be able to define
the group $A$ by the same procedure as in
\ref{ss:homogeneous-framings}.  For each ordered pair $(u,v)$ of
adjacent elements in $V$, let $r(u,v)$ be the word in $V$ defined by 
\[r(u,v)=v_1v_2\cdots v_{2m}.\]
Let $A$ be the quotient of $F_V$ (the free group on $V$) by the normal
subgroup generated by $\{v\bar{v}\;|\; v\in
V\}\cup\{r(u,v)\;|\;\mbox{$u$ and $v$ are adjacent}\}$.  Let
$\alpha_v$ denote the image of $v$ in $A$ and let
$\cA=\{\alpha_v\}_{v\in V}$ be the corresponding set of generators.
Finally, for each ordered pair $(u,v)$ of adjacent elements we define
a sequence $a_1(u,v),a_2(u,v),\ldots$ of elements of $A$ by the
formula:
\[\hspace*{2in}a_k(u,v)=\alpha_{v_1}\alpha_{v_2}\cdots\alpha_{v_k}\hspace{1.5in}
(\ast\ast)\] 
where $v_1,v_2,\ldots$ is the sequence defined by ($\ast$).  Our last
condition is the following:
\begin{enumerate}
\item[(4)](i) $\alpha_v\neq 1$, for each $v\in V$.\\
(ii) $\alpha_u\neq\alpha_v$, if $u\neq v$.\\
(iii) For each ordered pair $(u,v)$ of adjacent elements, let
$(a_k)=(a_k(u,v))$ be the above sequence and let $m=m(u,v)$.  Then the
elements $1,a_1,a_2,\ldots,a_{2m-1}$ are distinct.
\end{enumerate}
   
\begin{theorem}\label{thm:gluings-to-group}
Suppose $V$, $M$, and $L$ are as above and that $\{j_v\}_{v\in V}$ is
a set of gluing isomorphisms satisfying conditions (1), (2), (3), and
(4), above.  Then there is a simply connected, $2$-dimensional Coxeter
cell complex $X$ and a group $A$ of automorphisms of $X$ such that (a)
the link of each vertex is $L$, (b) $A$ is simply transitive on
$X^{(0)}$, and (c) $\{j_v\}_{v\in V}$ is the corresponding set of
gluing isomorphisms.
\end{theorem}

\begin{proof}{}
The Cayley $2$-complex associated to the presentation of $A$ is such
an $X$.  (Condition (4) is needed to check that $X$ is a Coxeter cell
complex; for example, (4)(iii) implies that the $2$-cell associated to
$r(u,v)$ is a $2m(u,v)$-gon.)
\end{proof}

\begin{remark}\label{rem:gluing2cat0group}
If, in addition, $X$ satisfies Condition (M2) of \ref{ss:hat-X}, then 
applying Theorems~\ref{thm:orbihedral-completion} and
\ref{thm:CAT(0)-groups} we see that $X$ is actually the $2$-skeleton
of a CAT($0$) Coxeter cell complex $\hat{X}$ on which $A$ acts as a
group of isometries.
\end{remark}

\begin{example}
Here we present an example of a finite $2$-dimensional Coxeter cell
complex $X$ which violates the conclusion of the previous remark.
Although $X$ will be homogeneous and simply connected, it cannot be
completed to a CAT($0$) complex $\hat{X}$.  The reason is that
Condition (M2) of \ref{ss:hat-X} does not hold: there is nontrivial
holonomy around each $2$-cell.

Suppose that $V$ consists of four elements $\{a,b,c,d\}$, that $L$ is
the complete graph on $V$ (i.e., $L$ is the $1$-skeleton of a
$3$-simplex), and that $M$ is the right-angled Coxeter matrix
associated to $L$ (i.e., $m(u,v)=2$ for any two distinct elements
$u,v\in V$).  Since the automorphism group of $L$ is $S_4$ (the
symmetric group on $4$ letters), we can specify candidates for gluing
isomorphisms by giving $4$ permutations.  We choose the following
involutions: 
\[\begin{array}{l}
j_a=(b\,d)\\
j_b=(c\,d)\\
j_c=(a\,d)\\
j_d=\Id.\end{array}\]
(Thus, $j_a$, $j_b$, and $j_c$ are transpositions.)  Conditions (1)
and (2) clearly hold.  Next we need to check Condition (3) for each
pair $(u,v)$ of distinct elements in $V$.  Since reversing the order
of $(u,v)$ only reverses the order of the $v_k$, it suffices to check
each of the $6$ unordered pairs.  In each case, we get a relation of
length four:
\[\begin{array}{ll}
r(a,b)=bada,\hspace{1in} & r(a,d)=daba\\
r(b,c)=cbdb,\hspace{1in} & r(b,d)=dbcb\\
r(c,d)=dcac,\hspace{1in} & r(a,d)=acdc.\end{array}\]
Notice that the relations in the same row differ by a cyclic
permutation, so we really have only three independent relations.  

This gives a presentation for a group
\[A=\langle a,b,c,d\,|\,
a^2,b^2,c^2,d^2,r(a,b),r(b,c),r(c,d)\rangle.\]
It can be checked that $A$ is isomorphic to the dihedral group of
order $14$, and it follows from this that condition (4) holds.  As
above, let $X$ be the Cayley $2$-complex of the presentation.  It is
simply connected and $A$ acts simply transitively on $X^{(0)}$.  Each
$2$-cell is a square and the link of each vertex is $L$.  $X$
obviously cannot be completed to a CAT($0$) complex $\hat{X}$.  Indeed,
the Coxeter cell corresponding to $L$ is the $4$-cube, so $\hat{X}$
would have to be the $4$-cube.  But the $4$-cube has $16$ vertices
while $X$ has only $14$.

The problem is that Condition (M2) does not hold.  (It is not
automatic since $L$ is not rigid.)  For example, if we compute the
holonomy around a $2$-cell labeled $adab$, we get
$j_b\circ j_a\circ j_d\circ j_a=(c\,d)$.
\end{example}
 
\section{A linear representation}\label{s:linearity}

\subsection{The partial order}
Let $M$ be a Coxeter matrix on the set $V$, and let $\{j_v\}_{v\in V}$
be a set of gluing maps as defined in \ref{ss:gluing2group}.
In addition, we assume that the involution $v\mapsto\bar{v}$ on $V$ is
trivial, thus each $j_v$ is an automorphism of $V_v$.  In this section
we describe conditions that allow us to define a linear representation
of the resulting group $A$.  We will show that when $M$, $V$, and 
$\{j_v\}$ arise from an $\cR$-blow-up $\Sigma_{\#}$, these conditions
are satisfied, giving a linear representation of the group $A$ acting
on the universal cover $\tilde{\Sigma}_{\#}$ (as described in
\ref{ss:group-A}).  This representation generalizes the standard
geometric representation of a Coxeter group.     

Our first assumption is that there is a partial order
(denoted by $<$) on $V$ satisfying the following condition:
\begin{enumerate}
\item[(P)] (i) If $u$ and $v$ are comparable (i.e., $u < v$ or $v <
u$), then they are adjacent (i.e., $m(u,v)<\infty$).\\
(ii) If $3\leq m(u,v)<\infty$, then $u$ and $v$ are minimal.\\
(iii) If $u'\leq u$, $v'\leq v$, and $u$ and $v$ are noncomparable and
adjacent, then $u'$ and $v'$ are noncomparable and adjacent.\\
(iv) Let $\bar{M}$ denote the Coxeter matrix $M$ with all $\infty$ entries
replaced by $2$'s.  Then for each $v$ in $V$, the subset $V_{<v}$
($=\{u\in V\;|\; u<v\}$) is spherical with respect to $\bar{M}$.
\end{enumerate}

Next we consider the set of gluing automorphisms $\{j_v:V_v\rightarrow
V_v\}_{v\in V}$ with $j_v(v)=v$ as in \ref{ss:gluing2group}.  In the
presence of the partial order, we can formulate the following
condition on the $j_v$'s, which is much simpler than are the
conditions in \ref{ss:gluing2group}.

\begin{enumerate}
\item[(C)](i) $j_v$ is an involution.\\
(ii) If $v_1,v_2\in V_v-\{v\}$, then
$m(j_v(v_1),j_v(v_2))=m(v_1,v_2)$.\\
(iii) $j_v$ preserves the partial order on $V_v$.\\
(iv) If $j_v(u)\neq u$, then $u<v$.\\
(v) If $u<v$ and $u'=j_v(u)$, then $j_uj_{v}j_{u'}j_v(y)=y$
for all $y<u$.
\end{enumerate}

\begin{lemma}
If $\{j_v\}_{v\in V}$ satisfies Conditions (C)(i)-(iv), then Conditions (1),
(2), and (3) of \ref{ss:gluing2group} hold. 
\end{lemma}

\begin{proof}{}
(C)(i) implies Condition (1), and (C)(ii) is the same as Condition
(2).  Consider Condition (3) of \ref{ss:gluing2group}.  If $3\leq
m(u,v)<\infty$, then by (P)(ii) $u$ and $v$ are minimal and hence, by
(C)(iv), $j_v(u)=u$ and $j_u(v)=v$ and the sequence
$v_0,v_1,v_2,v_3,\ldots$ is $u,v,u,v,\ldots$.  The same conclusion
holds if $m(u,v)=2$ and $u$ and $v$ are noncomparable.  Hence,
Condition (3) holds in these cases.  Finally, if $u<v$, then
$m(u,v)=2$ and the sequence $v_0,v_1,v_2,v_3,\ldots$ begins as
$u,v,j_v(u),v,\ldots$ and has the same period, $4$, so Condition (3)
again holds.  
\end{proof}

\begin{remark}\label{rem:holonomy} 
Condition (C)(v) is related to the holonomy condition (M2) of 
\ref{ss:hat-X}.  Given any pair $(u,v)$ with $m=m(u,v)<\infty$, the
sequence $v_0=u,v_1=v,v_2=j_{v_1}(v_0),\ldots$ gives rise to the
holonomy automorphism $j_{v_{2m}}\cdots j_{v_1}:V_{uv}\rightarrow
V_{uv}$ where $V_{uv}=\{y\in V\;|\; \{u,v,y\}\; \mbox{is
spherical}\}$.  For most pairs, this holonomy will be trivial as a
consequence of the properties listed in (P) and (C)(i)-(iv).  For
example, if $m(u,v)\geq 3$, then $u$ and $v$ are minimal; thus $j_u$
and $j_v$ are identity maps so the holonomy automorphism $j_uj_v\cdots
j_uj_v$ is trivial.  However, to ensure trivial holonomy in
the case $u<v$ the additional condition (C)(v) is needed.
\end{remark}

Henceforth, we assume that (C) holds.  It follows that we can define
the group $A$ as in \ref{ss:gluing2group}.  Condition (4) of
\ref{ss:gluing2group} then becomes the following:
\begin{enumerate}
\item[(4$^{\prime}$)](i) $\alpha_v\neq 1$.\\
(ii) $\alpha_u\neq\alpha_v$ if $u\neq v$.\\
(iii) If $u$ and $v$ are adjacent and noncomparable, then the order of
$\alpha_u\alpha_v$ is $m(u,v)$.  If $u<v$ and $u'=j_v(u)$, then the
elements $1,\alpha_u,\alpha_u\alpha_v,\alpha_u\alpha_v\alpha_{u'}$ are
distinct. 
\end{enumerate} 
In the next subsection, we shall define a representation for the group
$A$ and use it to show that (4$^{\prime}$) always holds (cf.,
Corollary~\ref{cor:poset2cat0gp}, below).   

\begin{lemma}
Let $\Sigma_{\#}$ denote the $\cR$-blow-up of $\Sigma(W,S,L)$, and let
$(V,M)=(S_{\#},M_{\#})$.  Then the partial order on $S_{\#}$ defined by
inclusion satisfies Condition (P) and the natural set of gluing
involutions $\{j_T\}_{T\in S_{\#}}$ satisfies Condition (C).    
\end{lemma}

\begin{proof}{}
If $T$ and $T'$ are comparable, then $m_{\#}(T,T')=2$, so $T$ and $T'$
are adjacent.  Thus, (P)(i) holds.  If $3\leq m_{\#}(T,T')<\infty$,
then $T$ and $T'$ are singleton subsets of $S$, hence they are
minimal, so (P)(ii) holds.  For (P)(iii), we just note that two
subsets $T$ and $U$ are adjacent and noncomparable if and only if
they are both minimal or they are completely disjoint.  Thus, if $T$
and $U$ are noncomparable and adjacent, $T'\subset T$, and $U'\subset
U$, then $T'$ and $U'$ must be noncomparable and adjacent.  

To show that (P)(iv) holds, suppose $T\in S_{\#}$, and let $\bar{M}_{\#T}$
denote the restriction of $\bar{M}_{\#}$ to $(S_{\#})_{<T}$.  We have to
show that $\bar{M}_{\#T}$ is the matrix for a finite Coxeter group.  Let
$\Gamma_T$ denote the Coxeter diagram for the (finite) special
subgroup $W_T$.  Then the Coxeter diagram for $\bar{M}_{\#T}$ is obtained
from $\Gamma_T$ in the following two steps.  First, for each subset
$T'\in (S_{\#})_{<T}$ of the form $T'=\{s,s'\}$, we delete the edge
joining $s$ and $s'$.  (This corresponds to replacing the
$\infty$-entry $m_{\#}(\{s\},\{s'\})$ of $M_{\#}$ with a $2$.)  Since
the resulting diagram represents a product of special subgroups of
$W_T$, its Coxeter group is finite.  Second, for each
nonsingleton element $T'\in (S_{\#})_{<T}$, we add a new disjoint
node.  (This new node is not connected to any other node 
since $m_{\#}(T',T'')$ is either $2$ or $\infty$ for any nonminimal
$T'\in S_{\#}$.)  Since adding a disjoint node to a Coxeter diagram
corresponds to adding a $\bbZ_2$ factor to its Coxeter group, the
resulting diagram represents a finite Coxeter group. 

Condition (C) follows directly from the definition of the
$j_T$'s given in \ref{ss:blow-up-gluing}.
\end{proof}

\subsection{The representation}
Let $M$ be a Coxeter matrix on a set $V$, and let $E$ denote the
vector space $\bbR^V$ with standard basis $\{\epsilon_v\}_{v\in V}$.
For each $t\in\bbR$, we define a symmetric bilinear form $B_t$
($=B_t(M)$) on $E$ by    
\[B_t(\epsilon_u,\epsilon_v)=\left\{\begin{array}{ll}
-\cos(\pi/m(u,v))&\mbox{if $m(u,v)<\infty$}\\  
-t &\mbox{if $m(u,v)=\infty$.}\end{array}\right.
\] 
(Note that when $t=1$, $B_t(M)$ is the canonical bilinear form
associated to the Coxeter matrix $M$.)

\begin{lemma}
Assume that (P) holds for $(V,M)$ and that $\{j_v\}_{v\in V}$ is a set
of gluing involutions satisfying Condition (C).  For each $v\in V$,
let $E_v\subset E$ denote the subspace defined by 
\[E_v=\Span\{\epsilon_{u}-\epsilon_{u'}\;|\;\mbox{$u\in V_v$ and
$u'=j_v(u)$}\}.\]  
Then for $t$ sufficiently large, the restriction of $B_t$ to $E_v$ is
nondegenerate for all $v\in V$. 
\end{lemma}

\begin{proof}{}
Let $v\in V$, and let $U_v\subset E$ be the subspace
$U_v=\Span\{\epsilon_u\;|\;u\in V_{<v}\}$.  
We let $\bar{B}$ denote the canonical bilinear form associated to the
Coxeter matrix $\bar{M}$.  Then when $t=0$, $B_t(M)$ coincides with
$\bar{B}$, and the restriction $B_t|_{U_v}$ coincides with the
restriction $\bar{B}|_{U_v}$.  By (P)(iv), the latter is positive
definite, hence $\det(\bar{B}|_{U_v})\neq 0$.  It follows that
$\det(B_t|_{U_v})$ is a nonzero polynomial in $t$, so for $t$
sufficiently large, $B_t|_{U_v}$ is nondegenerate.  To complete the
proof, we note that $E_v\subset U_v$ (by Condition (C)(iv)) and that
$U_v$ splits as an orthogonal direct sum   
\[U_v=E_v\oplus\Span\{\epsilon_u+\epsilon_{u'}\;|\;\mbox{$u\in V_{<v}$ 
and $u'=j_v(u)$}\}\]
(by (C)(ii)).  Thus, for $t$ sufficiently large, $B_t|_{E_v}$ is
nondegenerate.  
\end{proof}
   
Let $B_t$ be one of the bilinear forms that is nondegenerate on the
subspace $E_v$ for all $v\in V$.  Let $v\in V$.  The definition of the 
gluing involution $j_v$ implies that $\epsilon_v$ is orthogonal to the
subspace $E_v$. Moreover, since $B_t(\epsilon_v,\epsilon_v)=1$ and $B_t$ is
nondegenerate on $E_v$, we know that $B_t$ is nondegenerate on $\bbR
\epsilon_v\oplus E_v$.  Letting $F_v$ denote the orthogonal complement of
$\bbR \epsilon_v\oplus E_v$, we then have an orthogonal decomposition
\[E=\bbR \epsilon_v\oplus E_v\oplus F_v.\]
With respect to this decomposition, we define an involution
$\rho_v:E\rightarrow E$ by the formula 
\[\rho_v=-\Id|_{\bbR \epsilon_v}\oplus -\Id|_{E_v}\oplus\Id|_{F_v}.\]
It is clear that $\rho_v$ preserves the bilinear form $B_t$, that
$\rho_v(\epsilon_v)=-\epsilon_v$ and that $\rho_v(E_v)=E_v$ for all $v\in V$. 

\begin{lemma}\hspace{.1in}\label{lem:adjacent-pairs}
\begin{enumerate}
\item[(1)] If $v$ is minimal, then $E_v=\{0\}$ and
$\rho_v$ is the orthogonal reflection across the hyperplane
$F_v=\epsilon_v^{\perp}$.  In other words,
$\rho_v(\lambda)=\lambda-2B_t(\lambda,\epsilon_v)\epsilon_v$ for all $\lambda\in E$.
\item[(2)] If $u<v$, then $\rho_v(\epsilon_u)=\epsilon_{u'}$ and
$\rho_v(E_u)=E_{u'}$ where $u'=j_v(u)$.
\item[(3)] If $u$ and $v$ are noncomparable and $m(u,v)=2$, then
$\rho_v(\epsilon_u)=\epsilon_u$ and $\rho_v(E_u)=E_{u}$.  
\end{enumerate}
\end{lemma}
  
\begin{proof}{}
For (1), if $v$ is minimal, the involution $j_v$ is trivial.  This means
$E_v=\{0\}$, so $\rho_v$ is the orthogonal reflection with the
given formula.  For (2), we note that both $u$ and $u'$ are adjacent to
$v$, and $v$ is nonminimal; hence, by (P)(ii), $m(u,v)=m(u',v)=2$.  It
follows that the vector $\epsilon_u+\epsilon_{u'}$ is orthogonal to
$\epsilon_v$.  By (C)(ii), this vector is also orthogonal to the
subspace $E_v$.  So $\epsilon_{u}+\epsilon_{u'}$ is in $F_v$ and,
thus, fixed by $\rho_v$.  We then have 
\[\rho_v(2\epsilon_u)=\rho_v(\epsilon_{u}-\epsilon_{u'})+\rho_v(\epsilon_{u}+\epsilon_{u'})=
-(\epsilon_{u}-\epsilon_{u'})+(\epsilon_{u}+\epsilon_{u'})=2\epsilon_{u'}.\]
Hence, $\rho_v(\epsilon_u)=\epsilon_{u'}$.  To see that
$\rho_v(E_u)=E_{u'}$, suppose $y < u$.  By (C)(v), we have
$j_uj_vj_{u'}j_v(y)=y$ or, in other words, $j_vj_u(y)=j_{u'}j_v(y)$.
Applying $\rho_v$ to $\epsilon_y-\epsilon_{j_u(y)}\in E_u$, and using
the fact that $j_u(y)<u<v$, we obtain
\[\rho_v(\epsilon_y-\epsilon_{j_u(y)})=
\rho_v(\epsilon_{y})-\rho_v(\epsilon_{j_u(y)})=
\epsilon_{j_v(y)}-\epsilon_{j_{v}j_u(y)}= 
\epsilon_{j_v(y)}-\epsilon_{j_{u'}j_v(y)}.\]
Since this last term is in $E_{u'}$, it follows that
$\rho_v(E_{u})\subset E_{u'}$.  The same argument shows 
$\rho_v(E_{u'})\subset E_{u}$; hence, $\rho_v(E_u)=E_{u'}$.  The proof
of (3) is similar. 
\end{proof}
  
\begin{theorem}\label{thm:A-representation}
The map $\alpha_v\mapsto\rho_v$ extends to a homomorphism
$\rho:A\rightarrow GL(E)$.  (In fact the image of $\rho$ lies
in the orthogonal subgroup $O(B_t)\subset GL(E)$.) 
\end{theorem} 

\begin{proof}{}
$A$ is the group with one generator $\alpha_v$ for each $v\in V$ and
relations of the form 
\begin{enumerate}
\item[(a)] $(\alpha_v)^2=1$ for all $v\in V$,
\item[(b)] $(\alpha_u\alpha_v)^{m}=1$ if $u$ and $v$ are
minimal and $m=m(u,v)$,
\item[(c)] $\alpha_u\alpha_v\alpha_{u'}\alpha_v=1$ if $u<v$ and
$u'=j_v(u)$, and 
\item[(d)] $(\alpha_u\alpha_v)^2=1$ if $u$ and $v$ are noncomparable
and $m(u,v)=2$.
\end{enumerate}
Since each $\rho_v$ is an involution (preserving $B_t$), it suffices
to show that relations (b), (c), and (d) hold under the
substitution $\alpha_v\mapsto\rho_v$.  In case (b), $j_u$ and
$j_v$ are trivial; hence, the involutions $\rho_u$ and $\rho_{v}$ are the
usual orthogonal reflections through the hyperplanes $F_u$ and
$F_v$, respectively.  Since $\rho_u$ and $\rho_v$ fix the
codimension-two subspace $F_u\cap F_v$ pointwise, it suffices to show
that $\rho_u\rho_v$ has order $m$ when restricted to
$\Span\{\epsilon_u,\epsilon_v\}$.  A simple 
calculation shows that $B_t$ is positive definite on
$\Span\{\epsilon_u,\epsilon_v\}$ and that $\rho_u\rho_v$ is a
rotation through an angle of $2\pi/m$.  Thus, $(\rho_u\rho_v)^m=\Id$.  

In case (c), Lemma~\ref{lem:adjacent-pairs} (part 2) implies that the
following diagram commutes:
\[\xymatrix{(\bbR \epsilon_{u}\oplus E_{u})\oplus F_{u} 
\ar[d]^{\rho_{u}=-\Id\oplus\Id}
\ar[r]^{\rho_v} &
(\bbR \epsilon_{u'}\oplus E_{u'})\oplus F_{u'}
\ar[d]^{\rho_{u'}=-\Id\oplus\Id}\\
(\bbR \epsilon_{u}\oplus E_{u})\oplus F_{u} 
\ar[r]^{\rho_v} &
(\bbR \epsilon_{u'}\oplus E_{u'})\oplus F_{u'}}.\]
Since all of these maps are involutions, we have 
$\rho_u\rho_v\rho_{u'}\rho_{v}=\Id$.

In case (d), Lemma~\ref{lem:adjacent-pairs}
(part 3) implies that we have the same commutative diagram as in
case (c) but with $u'=u$.  Thus, $\rho_u\rho_v\rho_u\rho_v=\Id$.
\end{proof}

\begin{remark}
It seems likely that all of the mock reflection groups considered in
this paper are linear.  Indeed, the representation $\rho : A \rightarrow
GL(E)$, constructed above, is probably always faithful; however, we do
not have a proof of this and the linearity of $A$ is an open
question. 
\end{remark}

\begin{example}\label{ex:a3-rep}
Let $S$ be the set $\{a,b,c\}$, and let $(W,S)$ be the Coxeter system
corresponding to the diagram $\bA_3$.  The corresponding
Coxeter cell $Z(W,S)$ is the $3$-dimensional permutohedron.
The collection $\cR=\{\{a,b\},\{b,c\}\}$ is admissible with
respect to $\cP=\cS-\{\{a,b,c\}\}$, and the corresponding $\cR$-blow-up
is a blow-up of $\partial Z$.  (The Coxeter tile is a pentagon, the
$2$-dimensional associahedron $\Delta_{\#}$ as in
Example~\ref{ex:minbu-assoc}.)  The corresponding group
$A$ has a generator for each element of
\[S_{\#}=\{\{a\},\{b\},\{c\},\{a,b\},\{b,c\}\}.\]  
With respect to the standard basis 
$\epsilon_a,\epsilon_b,\epsilon_c,\epsilon_{ab},\epsilon_{bc}$ for
$E=\bbR^5$, the family of forms $B_t$ is given by 
\[[B_t]=\left[\begin{array}{ccccc}
1&-t&0& 0 & -t \\
-t & 1 & -t & 0 & 0\\
0&-t&1&-t&0\\
0&0&-t&1&-t\\
-t&0&0&-t&1 \end{array}\right],\]
and the representation $\rho:A\rightarrow GL(\bbR^5)$ is defined by
the involutions 
\[\rho_a=\left[\begin{array}{ccccc}
-1&2t&0&0&2t\\
0&1&0&0&0\\
0&0&1&0&0\\
0&0&0&1&0\\
0&0&0&0&1\end{array}\right]\hspace{.17in}
\rho_b=\left[\begin{array}{ccccc}
1&0&0&0&0\\
2t&-1&2t&0&0\\
0&0&1&0&0\\
0&0&0&1&0\\
0&0&0&0&1
\end{array}\right]\hspace{.17in}
\rho_c=\left[\begin{array}{ccccc}
1&0&0&0&0\\
0&1&0&0&0\\
0&2t&-1&2t&0\\
0&0&0&1&0\\
0&0&0&0&1\end{array}\right]\]
\[\rho_{ab}=\left[\begin{array}{ccccc}
0&1&\frac{-t}{1+t}&0&\frac{t}{1+t}\\
1&0&\frac{t}{1+t}&0&\frac{-t}{1+t}\\
0&0&1&0&0\\
0&0&2t&-1&2t\\
0&0&0&0&1\end{array}\right]\hspace{.3in}
\rho_{bc}=\left[\begin{array}{ccccc}
1&0&0&0&0\\
\frac{t}{1+t}&0&1&\frac{-t}{1+t}&0\\
\frac{-t}{1+t}&1&0&\frac{t}{1+t}&0\\
0&0&0&1&0\\
2t&0&0&2t&-1\end{array}\right]\]
\end{example}

\begin{corollary}\label{cor:poset2cat0gp}
Suppose $(V,M)$ satisfies Condition (P), and $\{j_v\}_{v\in V}$ is a
set of involutions satisfying Condition (C).  Then Condition
(4$^{\prime}$) holds.  In particular, the triple $(V,M,\{j_v\})$
yields a group $A$ and an action of $A$ on a CAT($0$) Coxeter cell
complex. 
\end{corollary}   

\begin{proof}{}
(4$^{\prime}$)(i) follows since $\rho_v$ is a nontrivial involution for every
$v\in V$.  Similarly, since $u\neq v$ implies $\rho_u$ and $\rho_v$
have different $-1$-eigenspaces, we have $\rho_u\neq\rho_v$. Thus,
$\alpha_u\neq\alpha_v$, so (4$^{\prime}$)(ii) holds.  Cases (b) and (d) of the
proof of Theorem~\ref{thm:A-representation} show that if $u$ and $v$
are adjacent and noncomparable, then $\rho_u\rho_v$ has order
$m(u,v)$.  This means that the order of $\alpha_u\alpha_v$ is at least 
$m(u,v)$, and therefore exactly $m(u,v)$.  Similarly, if
$u<v$, then case (c) of Theorem~\ref{thm:A-representation} shows that
$\Id$, $\rho_u$, $\rho_u\rho_v$, and $\rho_u\rho_v\rho_{u'}$ are all
distinct.  Thus,
$1,\alpha_u,\alpha_u\alpha_v,\alpha_u\alpha_v\alpha_{u'}$ must be
distinct elements of $A$.  By Theorem~\ref{thm:gluings-to-group}, $A$
acts on a $2$-dimensional Coxeter cell complex, and since condition
(M2) holds (Remark~\ref{rem:holonomy}), this complex can be completed
to a CAT($0$) Coxeter cell complex 
(Remark~\ref{rem:gluing2cat0group}). 
\end{proof}

\section{Permutohedral tilings}\label{s:permuto}

Suppose that we are given a tiling of a Coxeter cell complex with
fundamental tile an $n$-simplex.  Its maximal blow-up (as defined 
in \cite{DJS}, Section 4.1) will then be a cubical $n$-manifold tiled
by permutohedra.  This situation can arise from a Coxeter system
$(W,S)$ of rank $n+1$ in two ways.  The first is where $W$ is finite
and we take the maximal blow-up of $\partial Z(W,S)$ (the boundary
of the Coxeter cell).  The second way occurs when $(W,S)$ is a
``simplicial'' Coxeter system (defined in subsection \ref{ss:simplicial-cs},
below), and we take the maximal blow-up of $\Sigma(W,S)$ (the complete
reflection tiling of type $(W,S)$).  As it turns out, it follows from
Theorem~\ref{thm:iso-to-reflection} that the universal covers of all
such examples 
yield the same right-angled tiling of $\bbR^n$ by permutohedra.  Hence,
the fundamental groups of any two closed $n$-manifolds that arise from
such constructions are commensurable.  (Essentially, the same result
was asserted in Section 4.2 of \cite{DJS}; however, as we shall
explain in subsection \ref{ss:bu-simp-arrangement}, below, there was a
mistake in the proof.) 

\subsection{The permutohedron}
There are three equivalent definitions of the $n$-dimensional
permutohedron $P$.  First, it can be defined as the convex polytope
obtained by truncating all of the faces of an $n$-simplex $\Delta^n$
of codimension $\geq 2$.  A second definition is that it is the
polytope whose boundary complex is dual to the barycentric subdivision
of $\partial\Delta^n$.  The third definition is that $P^n$ is the
Coxeter cell associated to the symmetric group $S_{n+1}$  ($S_{n+1}$
is the Coxeter group with Coxeter graph ${\bf A}_n$.)

Let us fix a set $S$ of cardinality $n+1$ and suppose that the
elements of $S$ index the codimension-one faces of $\Delta^n$.  Regard
$S_{n+1}$ as the symmetric group on $S$.  Let $V(P^n)$ denote the set
of all proper nonempty subsets of $S$.  Thus, $V(P^n)$ naturally
indexes the set of codimension-one faces of $P^n$.  Let $\cN(P^n)$
denote the poset of all subsets of $V(P^n)$ that are chains.  (If
$\cR=V(P^n)$, then $\cN(P^n)$ is the poset of all $\cR$-nested subsets
of $V(P^n)$ as in Definition~\ref{def:R-nested}.)  Thus,
$\cN(P^n)_{>\emptyset}$ is the barycentric subdivision of
$\partial\Delta^n$.  We shall also denote this simplicial complex by
$L(P^n)$.  

The action of $S_{n+1}$ on $S$ induces an action on $P^n$ as a group
of combinatorial automorphisms.  A fundamental domain for this action
is the associated Coxeter block $B^n$.  As is the case for any Coxeter
cell, the orbit space $P^n/S_{n+1}$ can be identified with this
Coxeter block.  The permutohedron has one further symmetry: the
antipodal map.  In fact, it is easy to see that its full group of
combinatorial symmetries, $\Aut(P^n)$, is just $S_{n+1}\times\bbZ_2$. 

\subsection{The reflection tiling}\label{ss:permuto-reflection-tiling}
Let $M'$ be the right-angled Coxeter matrix associated to the flag
complex $L(P^n)$ and let $W'$ be the associated Coxeter group.  That
is to say, $W'$ has a generator for each codimension-one face of $P^n$
(i.e., for each element of 
$V(P^n)$), and two such generators commute if and only if the
corresponding faces intersect.  Let $\Sigma_{P^n}$ be the complete
reflection tiling of type $(W',V(P_n))$.  The natural framing is
symmetric, and the group of frame-preserving automorphisms is
$A(\Sigma_{P^n})=W'$.    

The full symmetry group of $\Sigma_{P^n}$ is
$G(\Sigma_{P^n})=A(\Sigma_{P^n})\rtimes\Aut(P^n)$.  This has a 
subgroup of index two, $G_0(\Sigma_{P^n})$ defined by 
\[G_0(\Sigma_{P^n})=A(\Sigma_{P^n})\rtimes S_{n+1}.\]
In fact, $G_0(\Sigma_{P^n})$ is a also a Coxeter group and its action on $\Sigma_{P^n}$
is as a group generated by reflections.  A fundamental chamber for
this action is the Coxeter block $B^n$.  Thus,
$\Sigma_{P^n}/G_0(P^n)\cong P^n/S_{n+1}\cong B^n$.  (The Coxeter diagram for
$G_0(P^n)$ is given in Figure 7, page 536 of \cite{DJS}.)
 
\subsection{Simplicial Coxeter systems}\label{ss:simplicial-cs}

\begin{definition}
A Coxeter system $(W,S)$ is {\em simplicial} if $W$ is infinite and
each proper subset of $S$ is spherical.
\end{definition}

Suppose that $(W,S)$ is simplicial and that $\card(S)=n+1$.  Then the
fundamental chamber for the $W$-action on $\Sigma(W,S)$ is an
$n$-simplex.  In 1950 in \cite{L}, Lanner showed that each such
$\Sigma(W,S)$ can be identified with either hyperbolic $n$-space
$\bbH^n$ or Euclidean $n$-space $\bbE^n$ so that the $W$-action is as
a classical group of isometries generated by reflections across the
faces of a hyperbolic or Euclidean $n$-simplex.  He also listed
possible Coxeter diagrams of simplicial Coxeter systems.

In the Euclidean case, there are four families in each dimension
$n\geq 4$.  Their Coxeter diagrams are denoted $\tilde{\bf A}_n$,
$\tilde{\bf B}_n$ ($n\geq 2$), $\tilde{\bf C}_n$ ($n\geq 3$), and
$\tilde{\bf D}_n$ ($n\geq 4$).  There are also five exceptional
Euclidean simplicial Coxeter systems: $\tilde{\bf G}_2$, $\tilde{\bf
F}_4$, $\tilde{\bf E}_6$, $\tilde{\bf E}_7$, $\tilde{\bf E}_8$.  In
the hyperbolic case, in dimension two, there are the $(p,q,r)$
triangle groups (where $p^{-1}+q^{-1}+r^{-1}<1$).  In dimension three,
there are the nine tetrahedral hyperbolic Coxeter systems, and in
dimension four there are five more hyperbolic examples.  Finally,
there are no hyperbolic simplicial Coxeter systems in dimensions $>4$.
(The Coxeter diagrams of all these groups can be found, for example,
in \cite{Bo} or \cite{L}). 

Because $\Sigma(W,S)$ is either $\bbE^n$ or $\bbH^n$, its quotient by
any torsion free subgroup $\Gamma\subset W$ will be a manifold $M$.  If
$\Sigma_{\#}$ is any $\cR$-blow-up, then the quotient
$M_{\#}=\Sigma_{\#}/\Gamma$ is a blow-up of $M$.  

\begin{remark}
The group $W$ corresponding to the $\tilde{\bf A}_{n+1}$ Coxeter
diagram is the semidirect product of $\bbZ^n$ with the symmetric group
$S_{n+1}$.  In this case, the fundamental tile of the minimal blow-up
$\Sigma_{\#}$ is a polytope called the ``cyclohedron''.  The quotient
of this tiling by $\bbZ^n$ is discussed in \cite{De} in relation to
compactifications of configuration spaces.
\end{remark}  

\subsection{Maximal blow-ups of simplicial tilings}
Suppose $(W,S)$ is a Coxeter system with $\card(S)=n+1$ and with $W$
finite.  Let $\tilde{Z}_{\#}$ be the universal cover of the maximal
blow-up of the Coxeter cell $Z(W,S)$, and let $\ccc$ be the
universal cover of $(\partial Z)_{\#}$.  Since $Z_{\#}$ is an interval
bundle over $(\partial Z)_{\#}/a_{\#}$ and $a_{\#}$ is a free
involution, $\tilde{Z}_{\#}$ is an interval bundle over $\ccc$.
It is easy to see that the natural framing on $\ccc$ obtained by
restricting the framing on $\cc_{\#}$ is symmetric and the group of
frame-preserving symmetries is the subgroup $A$ ($=A(\ccc)$) of
$A(W,S)$ (see \ref{ss:group-A}) generated by all $\alpha_T$ where $T$
is a proper nonempty subset of $S$.   

\begin{theorem}\label{thm:max-boundary}
Suppose $(W,S)$ is a Coxeter system with $\card(S)=n+1$ and with $W$
finite.  Let $\ccc$ denote the universal cover of the maximal
blow-up of $\partial Z(W,S)$, and let $A$ ($=A(\ccc)$) be the
group described above.  Then 
\begin{enumerate}
\item[(1)] $\ccc$ is isomorphic to $\Sigma_{P^n}$, and 
\item[(2)] $A$ is isomorphic to a subgroup of index $n!$ in
$G_0(\Sigma_{P^n})$.
\end{enumerate}
\end{theorem}

\begin{theorem}\label{thm:max-simplicial}
Suppose $(W,S)$ is a simplicial Coxeter system, with $\card(S)=n+1$.
Let $\tilde{\Sigma}_{\#}$ denote the universal cover of the maximal
blow-up of $\Sigma(W,S)$, and let $A$ ($=A(\tilde{\Sigma}_{\#})$)
denote its frame-preserving symmetry group from \ref{ss:group-A}.
Then 
\begin{enumerate}
\item[(1)] $\tilde{\Sigma}_{\#}$ is isomorphic to $\Sigma_{P^n}$, and 
\item[(2)] $A$ is isomorphic to a subgroup of index $n!$ in
$G_{0}(\Sigma_{P^n})$. 
\end{enumerate}
\end{theorem}

Both of these theorems follow from
Theorem~\ref{thm:iso-to-reflection}.  We must first verify that  
Conditions (E) and (H) hold.  For Condition (E), we must show that for each
proper, nonempty subset $T$ of $S$, the involution 
\[j_T:\cN(P^n)_{\geq\{T\}}\rightarrow\cN(P^n)_{\geq\{T\}}\]
extends to an involution of $\cN(P^n)$.  The map $j_T$ is induced by
the permutation of $T$ defined by $a_Tw_T$ (we continue to denote this
permutation by $j_T$).  Extend $j_T$ to a permutation $\tilde{j}_T$ of
$S$ by letting $\tilde{j}_T|_{S-T}$ be the identity permutation of
$S-T$.  Thus, $\tilde{j}_T\in S_{n+1}$ and hence can naturally be
regarded as an element of $\Aut(\cN(P^n))$ ($=\Aut(P^n)$).  

Condition (H) follows from Remark~\ref{rem:rigid} and the fact
that the link of a vertex in $X$ (or $\tilde{\Sigma}_{\#}$) is the boundary
complex of the $n$-dimensional octahedron.  

\begin{example}
As in Theorem~\ref{thm:max-simplicial} suppose that $(W,S)$ is a
simplicial Coxeter system of rank $n+1$, that $\Sigma=\Sigma(W,S)$,
and that $\Sigma_{\#}$ is the maximal blow-up.  Let $H$ be a
torsion-free subgroup of finite index in $W$.  Then $\Sigma/H$ and
$\Sigma_{\#}/H$ are nonpositively curved, closed $n$-manifolds.  If
$\phi:A(\tilde{\Sigma}_{\#})\rightarrow W$ denotes the natural
projection, then $\pi_1(\Sigma_{\#}/H)\cong\phi^{-1}(H)$.  So,
Theorem~\ref{thm:max-simplicial} implies that $\pi_1(\Sigma_{\#}/H)$
is also a subgroup of finite index in $G_0(\Sigma_{P^n})$.
\end{example}

\subsection{Maximal blow-ups of simplicial
arrangements}\label{ss:bu-simp-arrangement}
Suppose we are given a simplicial hyperplane arrangement in
$\bbR^{n+1}$.  Let $Z$ be the corresponding zonotope, and let $K$ be
the triangulation of $S^n$ induced by the arrangement.  In
Corollary~4.2.8 of \cite{DJS}, we asserted that the universal cover of
the maximal blow-up $(\partial Z)_{\#}$ of $\partial Z$ could be
identified with $\Sigma_{P^n}$.  Although this is correct, the proof given in
\cite{DJS} is not.  We take this opportunity to correct it.  

The ``proof'' of \cite{DJS} had three steps.
\begin{enumerate}
\setlength{\itemindent}{.2in}
\item[\em Step 1:] There is a simplicial projection (or ``folding map'')
$p:K\rightarrow\Delta^n$ (Lemma~4.2.6 in \cite{DJS}).
\item[\em Step 2:] The map $p$ induces $p_{\#}:(\partial
Z)_{\#}\rightarrow\Delta^n_{\#}=P^n$ (Corollary~4.2.7 in \cite{DJS}). 
\item[\em Step 3:] The map $p_{\#}$ is the projection map of an orbifold
covering.
\end{enumerate}

Step 2 is incorrect--the map $p_{\#}$ is not well-defined.  Of course,
the problem is caused by the fact that $p$ need not be compatible with
the antipodal map $a$: if $\sigma^n$ is an $n$-simplex in $K$ then, in
general, $p|_{\sigma^n}\neq p\circ a|_{-\sigma^n}$.  (This phenomenon
causes a problem when we are considering the normal arrangement to a
subspace.)  However, if we divide out by $S_{n+1}$, the symmetry group
of $\Delta^n$, then we do get a well-defined map $(\partial
Z)_{\#}\rightarrow\Delta^n_{\#}/S_{n+1}$.  To see this, let $\Delta'$
be an $n$-simplex in the barycentric subdivision of $\Delta^n$.  Then
$\Delta^n/S_{n+1}=\Delta'$, so we have a folding map
$q:\Delta^n\rightarrow\Delta'$.  Let $p'=q\circ
p:K\rightarrow\Delta'$.  Then Step 2 can be replaced by the following:
\begin{enumerate}
\setlength{\itemindent}{.25in}
\item[\em Step 2\/$'$:] $p'$ induces a map $p'_{\#}:(\partial
Z)_{\#}\rightarrow\Delta_{\#}/S_{n+1}$.  
\end{enumerate}

As for the last step, we have
$\Delta_{\#}/S_{n+1}=P^n/S_{n+1}=\Sigma_{P^n}/G_{0}(\Sigma_{P^n})$.  Furthermore, the
map $p'_{\#}:(\partial Z)_{\#}\rightarrow \Sigma_{P^n}/G_0(\Sigma_{P^n})$ is an orbifold
covering.  (This just means that the map is locally isomorphic to
$\bbR^n\rightarrow\bbR^n/H$ where the finite group $H$ is either a
subgroup of $(\bbZ_2)^n$ or $S_{n+1}$.)  Therefore, we have proved the
following result.

\begin{theorem}\label{thm:correction}
Suppose that $Z$ is an $(n+1)$-dimensional simple zonotope (i.e., it
is associated to a simplicial hyperplane arrangement in $\bbR^{n+1}$)
and let $(\partial Z)_{\#}$ denote the maximal blow-up of its
boundary.  Then 
\begin{enumerate}
\item[(1)] The universal cover of $(\partial Z)_{\#}$ is isomorphic to
$\Sigma_{P^n}$, and 
\item[(2)] $\pi_1((\partial Z)_{\#})$ is a subgroup of finite index in
$G_0(\Sigma_{P^n})$.
\end{enumerate}
\end{theorem}
 
\begin{remark}
Note that this proof of Theorem~\ref{thm:correction} gives another
proof of Theorem~\ref{thm:max-simplicial}.
\end{remark}

\section{Associahedral tilings}\label{s:associahedral-tilings}
Manifolds tiled by associahedra arise as minimal blow-ups of the
boundaries of certain Coxeter cells and as the minimal blow-ups of
$\Sigma(W,S)$ for certain simplicial Coxeter systems.  In contrast to
permutohedral tilings, the universal covers of these examples tend not
to be isomorphic (although they all give tilings of $\bbR^n$).  The
reason is that the associahedron is less symmetric than the
permutohedron.  More precisely, two adjacent associahedral tiles are glued
together by an involution $j_T$ of a codimension-one face, and these
gluing involutions tend not to all extend to symmetries of the full
associahedron.  

\subsection{The associahedron}\label{ss:associahedron}
Following \cite{Lee} we give two equivalent descriptions of a
simplicial complex $\cN_{>\emptyset}$ that is dual to the boundary
complex of the $n$-dimensional associahedron $K^n$.  The first
description is the one given in \ref{ss:finite-W-tiles}: it shows how
$K^n$ arises as a truncation of the $n$-simplex.  The second
description is in terms of diagonals in an $(n+3)$-gon: it is more
convenient for describing the group of combinatorial symmetries of
$K^n$.

Let $S$ be a set with $n+1$ elements and suppose that $\Gamma$ is a
graph with vertex set $S$ such that $\Gamma$ is homeomorphic to an
interval.  We might as well assume that $S=\{1,2,\ldots,n+1\}$ and
that $\Gamma$ is the interval $[1,n+1]$.  Let $V$ be the set of proper
nonempty subsets $T$ of $S$ such that the full subgraph $\Gamma_T$
spanned by $T$ is connected.  In other words, $V$ can be identified
with sets of consecutive integers of the form $[k,l]$ where $1\leq
k\leq l\leq n+1$ and $l-k<n$.  A {\em decomposition} of a subset $T$
of $S$ is a collection $\{T_1,\ldots,T_k\}$ of disjoint subsets of $T$
such that $T=T_1\cup\cdots\cup T_k$, each $T_i\in V$, and
$\Gamma_{T_1},\ldots,\Gamma_{T_k}$ are the connected components of
$\Gamma_T$.  A subset $\cT$ of $V$ is {\em nested} if either
$\cT=\emptyset$ or if the maximal elements $T_1,\ldots,T_k$ in $\cT$
give a decomposition of $T_1\cup\cdots\cup T_k$.  $\cN$ is defined as
the poset of all such nested subsets of $V$.  Then $\cN_{>\emptyset}$
is a simplicial complex of dimension $n-1$ which can be identified
with a certain simplicial subdivision of $\partial\Delta^n$ (see
\cite{Lee} or \cite{DJS}).  Moreover, it is proved in \cite{Lee} that
this subdivision of $\partial\Delta^n$ can be idenitifed with the
boundary complex of a convex simplicial polytope in $\bbR^n$.  The
dual polytope $K^n$ is the {\em $n$-dimensional associahedron}.  $K^0$
is a point, $K^1$ is an interval, and $K^2$ is a pentagon.  We shall
also denote the simplicial complex $\cN_{>\emptyset}$ by $L(K^n)$. 

The second description of this complex is more illuminating.  Let
$P_{n+3}$ be a regular $(n+3)$-gon.  A {\em diagonal} $d$ in $P_{n+3}$
is a line segment connecting two nonadjacent vertices of $P_{n+3}$.
Two diagonals $d$ and $d'$ are {\em noncrossing} if their interiors
are disjoint.  Let $V'$ denote the set of all diagonals of $P_{n+3}$.
We next define a simplicial complex $\cN'_{>\emptyset}$ with vertex
set $V'$.  A $k$-simplex om $\cN'_{>\emptyset}$ is a collection
$\sigma=\{d_0,\ldots,d_k\}$ of pairwise noncrossing diagonals.  So, a
maximal simplex in $\cN'_{>\emptyset}$ corresponds to a triangulation
of $P_{n+3}$ with no additional vertices.  The dimension of such a
maximal simplex is easily seen to be $n-1$.  

A bijection $V\rightarrow V'$ is defined as follows.  Number the
vertices of $P_{n+3}$ cyclically around the boundary by
$0,1,\ldots,n+2$.  Let $[k,l]\in V$.  The corresponding diagonal
$d_{[k,l]}$ is defined to be the diagonal connecting the vertices
$k-1$ and $l+1$ of $P_{n+3}$.  Clearly, this is a bijection.
Furthermore, it is easy to see (c.f., \cite{Lee}) that it induces an
isomorphism $\cN_{>\emptyset}\cong\cN'_{>\emptyset}$.  Henceforth, we
identify $V$ with $V'$ and $\cN_{>\emptyset}$ with
$\cN'_{>\emptyset}$.

Each diagonal $d\in V$ corresponds to a codimension-one face, $F(d)$,
of $K^n$.  Next, we investigate the combinatorial type of $F(d)$.

\begin{figure}[ht]
\begin{center}
\psfrag{d}{$d$}
\psfrag{Q}{$Q$}
\psfrag{Q'}{$Q'$}
\includegraphics[scale=.4]{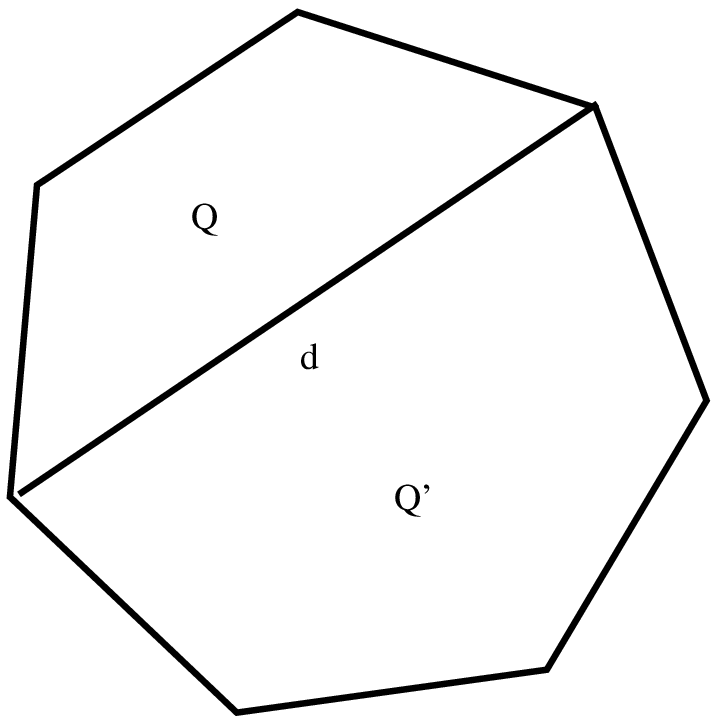}
\end{center}
\caption{\label{fig:polysides}}
\end{figure}

The diagonal $d$ divides $P_{n+3}$ into two polygons $Q$ and $Q'$ as
indicated in Figure~\ref{fig:polysides}.
Let $m(Q)+3$ and $m(Q')+3$ be the number of vertices of $Q$ and $Q'$,
respectively.  One checks easily that $m(Q)+m(Q')=n-1$.  Without loss
of generality, we may suppose that $m(Q)\leq m(Q')$.  Set
$m(d)=m(Q)$. 

It is completely straightforward to check that the link of $d$ in
$\cN_{>\emptyset}$ is isomorphic to the join of the corresponding
complexes of diagonals for $Q$ and $Q'$.  This implies the following
result.

\begin{proposition}\label{prop:faces}
Suppose the diagonal $d$ divides $P_{n+3}$ into two polygons $Q$ and
$Q'$ as above.  Then 
\[F(d)\cong K^{m(Q)}\times K^{m(Q')}.\]
\end{proposition}

We note that if $m(d)=0$, then $K^{m(Q)}$ is a point and hence $F(d)$
is an $(n-1)$-dimensional associahedron.  We will need to use the
following lemma in the next subsection.  

\begin{lemma}\label{lem:assoc-faces}
$K^n$ is not combinatorially isomorphic to a product of the form
$K^i\times K^{n-i}$, with $0<i<n$.
\end{lemma}

\begin{proof}{}
Let $v(n)={n+3\choose 2}-(n+3)$ be the number of diagonals in
$P_{n+3}$.  A computation shows that $v(n)\geq v(i)+v(n-i)$, with
equality if and only if $i=0$ or $i=n$.  In other words, for $0<i<n$,
$K^i\times K^{n-i}$ has fewer codimension-one faces than does $K^n$.
\end{proof}

\subsection{Symmetries of the associahedron}
Let $\Aut(K^n)$ ($=\Aut(\cN)$) denote the group of combinatorial
symmetries of $K^n$. 

The symmetry group of $P_{n+3}$ is $D_{n+3}$, the dihedral group of
order $2(n+3)$.  An element of $D_{n+3}$ takes diagonals to diagonals
and collections of noncrossing diagonals to collections of noncrossing
diagonals.  Hence, it gives an automorphism of $K^n$.  This defines a
homomorphism $\phi:D_{n+3}\rightarrow\Aut(K^n)$.

\begin{lemma}
For $n\geq 2$, $\phi:D_{n+3}\rightarrow\Aut(K^n)$ is an isomorphism.
\end{lemma}

\begin{proof}{}
For each $i=0,\ldots,n+2$, we let $d_i$ denote the diagonal with
$m(d_i)=0$ that cuts off the vertex $i$ from $P_{n+3}$.  We then let
$F_i$ be the corresponding codimension-one face $F(d_i)$.  By
Proposition~\ref{prop:faces} each $F_i$ is isomorphic to $K^{n-1}$ and
by Lemma~\ref{lem:assoc-faces} these are the only faces isomorphic to
$K^{n-1}$.   Thus, the collection $\{F_0,\ldots,F_{n+2}\}$ is
preserved by any combinatorial automorphism of $K^n$.  In fact, 
the relative positions of the vertices $0,1,\ldots,n+2$ on the circle
are uniquely determined by the poset $\cN_{>\emptyset}$.  To see this,
we just note that $d_i$ and $d_j$ are crossing diagonals (i.e.,
$\{d_i,d_j\}\not\in\cN_{>\emptyset}$) if and only if $j=i+1$ or
$j=i-1$ (modulo $n+3$).  It follows that any combinatorial
automorphism of $K^n$ must preserve the relative positions of
$0,\ldots,n+2$, hence $\phi$ is surjective. 

To see that $\phi$ is injective we simply note that the face $F_0$ is
stabilized by the two-element subgroup generated by the reflection
$r\in D_{n+3}$ across the line perpendicular to the diagonal $d_0$.
The reflection $r$ does not act trivially on $K^n$ since it exchanges
the diagonals $d_1$ and $d_{n+2}$ (and, hence, the faces $F_1$ and
$F_{n+2}$).  On the other hand, any nontrivial element of $D_{n+3}$
other than $r$ moves the diagonal $d_0$, hence does not act trivially
on $K^n$.  It follows that $\phi$ is injective.
\end{proof}

\begin{remark}
For $n=1$, $\phi:D_4\rightarrow\Aut(K^1)$ is not injective.  There are
two types of reflections in $D_4$.  A line of symmetry of the square
$P_4$ can connect either two opposite vertices or the midpoints of two
opposite edges:  we say that the coresponding reflection is of {\em
vertex type} or {\em edge type}, respectively.  Clearly, $\phi$ takes
each vertex type reflection to the identity element of $\Aut(K^1)$ and
each edge-type reflection to the nontrivial element.
\end{remark}

A nontrivial involution in $D_{n+3}$ is either the antipodal map (when
$n+3$ is even) or a reflection.  Let us say that an involution
$f:K^n\rightarrow K^n$ is of {\em R-type} if $f=\phi(r)$ for some
reflection $r\in D_{n+3}$. 

\begin{figure}[ht]
\begin{center}
\psfrag{d}{$d$}
\psfrag{L}{$L(d)$}
\includegraphics[scale=.5]{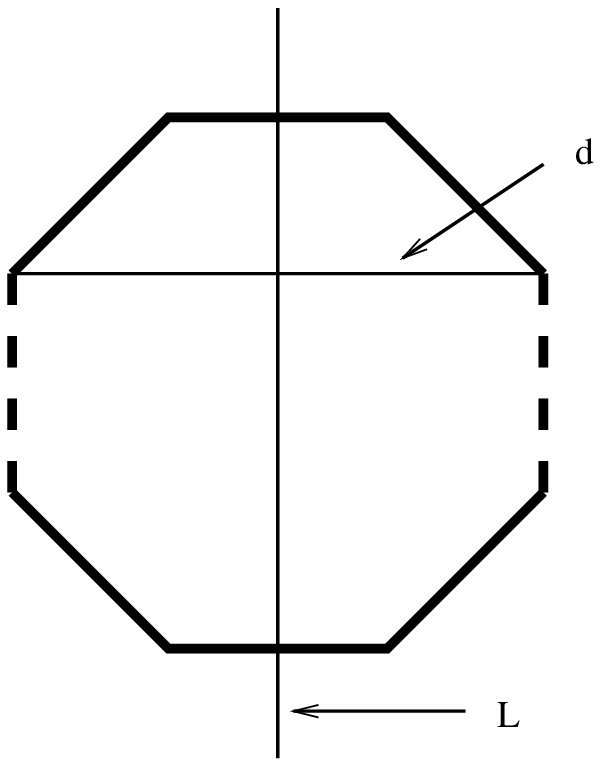}
\end{center}
\caption{\label{fig:stab-action}}
\end{figure}

Suppose $d$ is a diagonal of $P_{n+3}$ subdividing it into polygons
$Q$ and $Q'$.  Let $p_d$ be the midpoint of $d$ and let $L(d)$ be the
line perpendicular to $d$ at the point $p_d$.  Thus, $L(d)$ is a line
of symmetry of $P_{n+3}$.  The corresponding reflection in $D_{n+3}$ is
denoted by $r_{L(d)}$.  (See Figure~\ref{fig:stab-action}, below.)

A diagonal of $P_{n+3}$ is a {\em main diagonal} if $n+3$ is even and
if it connects opposite vertices of $P_{n+3}$.

The next lemma is geometrically clear.

\begin{lemma}\label{lem:d-stab}
Let $d$ be a diagonal $P_{n+3}$.
\begin{enumerate}
\item[(1)] If $d$ is not a main diagonal, then its stabilizer in
$D_{n+3}$ is the cyclic group of order $2$ generated by the reflection
$r_{L(d)}$.
\item[(2)] If $d$ is a main diagonal, then its stabilizer in $D_{n+3}$
is isomorphic to $\bbZ_2\times\bbZ_2$ and the generators can be taken
to be $r_{L(d)}$ and the reflection $r_d$ across $d$.
\end{enumerate}
\end{lemma}

\begin{corollary}\label{cor:d-stab-action}
Suppose $d$ is a diagonal of $P_{n+3}$, $n\geq 2$, and that $L=L(d)$
is the corresponding line of symmetry.
\begin{enumerate}
\item[(1)] If $m(d)=0$, then $F(d)\cong K^{n-1}$ and $\phi(r_L)$ acts
on $F(d)$ as a nontrivial R-type involution.
\item[(2)] If $m(d)>0$, then $F(d)\cong K^{m(d)}\times K^{n-m(d)-1}$
and the restriction of $\phi(r_L)$ to $F(d)$ can be written as
$\phi(r_L)=r_1\times r_2$ where $r_1:K^{m(d)}\rightarrow K^{m(d)}$ and  
$r_2:K^{n-m(d)-1}\rightarrow K^{n-m(d)-1}$ are both (nontrivial)
R-type involutions.
\item[(3)] Suppose $d$ is a main diagonal.  Then $F(d)\cong K^k\times
K^k$ where $k=(n-1)/2$ and $\phi(r_d)$ acts on $F(d)$ by switching the
factors.  Furthermore, the restrictions of $\phi(r_L)$ and $\phi(r_d)$
are not equal.
\end{enumerate}
\end{corollary}

\subsection{Gluing involutions}
Manifolds tiled by associahedra arise from minimal blow-ups in cases
where the Coxeter diagram is an interval.  The associahedra are glued
together via involutions $j_{T}$, $T\in V$, defined on the
codimension-one faces of $K^n$.  If $T$ corresponds to a subinterval
of the Coxeter diagram, then $j_T$ is induced by an involution of the
subinterval.  This involution might be trivial (e.g., if the subgraph
is of type $\bB_m$) or it might be the nontrivial involution that
flips the subinterval (e.g., if it is of type $\bA_m$, $m\geq 2$).  Here
we are concerned with the question of when $j_T$ extends to an
automorphism of $K^n$.  Since the identity always extends, we shall
analyze the case $j_T=i_T$, where $i_T$ is induced by the nontrivial
involution of the subinterval.

As in \ref{ss:associahedron}, let $S$ be the set of integers in
$[1,n+1]$.  $T$ will denote a proper nonempty subset of $S$ consisting
of consecutive integers.  By an abuse of notation, we will write
$T=[k,l]$ to mean $T=\{k,\ldots,l\}$.  Let $i_T$ be the order reversing
involution of $T$.  Then $i_T$ induces an involution of
$\cN_{\geq\{T\}}$, also denoted by $i_T$.  (If $\{T',T\}$ is a vertex of
$\cN_{>\{T\}}$ such that $T\subset T'$ or such that the subintervals
corresponding to $T$ and $T'$ are disjoint, then $i_T(T')=T'$.)  Its
geometric realization is again denoted $i_T$.  Let $d(T)$ be the
diagonal of $P_{n+3}$ corresponding to $T$.  

We note that $i_T:F(d(T))\rightarrow F(d(T))$ extends to an
automorphism of $K^n$ if and only if it coincides with the action of
an element of the stabilizer of $d(T)$ in $D_{n+3}$ on this face.
The main result of this subsection is the following key lemma, which
determines precisely when this happens.

\begin{lemma}\label{lem:mock-or-not}
The involution $i_T:F(d(T))\rightarrow F(d(T))$ lies in the stabilizer
of $d(T)$ in $D_{n+3}$ if and only if $m(d(T))=0$.
\end{lemma}

\begin{proof}{}
Suppose $T=[k,l]$.  So, $d(T)$ connects the vertices of $P_{n+3}$
numbered $k-1$ and $l+1$.  It divides $P_{n+3}$ into two polygons
$Q_1$ and $Q_2$ where $Q_2$ contains the vertices numbered
$k-1,\ldots,l+1$.  Let $L=L(d(T))$ be the line of symmetry for
$P_{n+3}$ that stabilizes $d(T)$.  Then $L$ is also a line of symmetry
for $Q_1$ and $Q_2$.  Let $r_2$ denote the restriction of $r_L$ to
$Q_2$.  

The codimension-one face $F(d(T))$ decomposes as
$F(d(T))=K^{m(Q_1)}\times K^{m(Q_2)}$.  It follows from the
definitions that 
\[i_T=\Id\times\phi(r_2):K^{m(Q_1)}\times
K^{m(Q_2)}\rightarrow K^{m(Q_1)}\times K^{m(Q_2)}.\]
Comparing this with Lemma~\ref{lem:d-stab} and
Corollary~\ref{cor:d-stab-action}, we see that $i_T$ does not belong
to the 
stabilizer of $d(T)$ unless $m(Q_2)=0$ (in which case $i_T$ is the
identity) or $m(Q_1)=0$ (in which case $i_T=\phi(r_2)$).
\end{proof}

\begin{remark}
There are two ways in which it can happen that $m(d(T))=0$.  The first
is that $T$ is a singleton.  In this case $i_T$ is the identity and
$F(d(T))$ is a reflection-type face.  The second way is that $T$ is a
maximal proper subinterval of $[1,n+1]$, i.e., $T=[1,n]$ or
$T=[2,n+1]$.  In both these cases, the involution $i_T$ is nontrivial,
but it is the restriction of a symmetry of the full associahedron.
\end{remark}

\subsection{Schl\"{a}fli symbols}  
As before, suppose that
$S=\{1,\ldots,n+1\}$ and that there is a corresponding Coxeter diagram
$\Gamma$ with underlying graph the interval $[1,n+1]$.  The labels on
the edges of $\Gamma$ are then given by an $n$-tuple
$(m_1,\ldots,m_n)$ of integers, each $\geq 3$, where $m_i$ is the
label on $[i,i+1]$.  Classically, this $n$-tuple is called the
{\em Schl\"{a}fli symbol} of $\Gamma$.  For example, the Schl\"{a}fli symbol
for ${\bf A}_{n+1}$ is $(3,\ldots,3)$, while for ${\bf B}_{n+1}$ it is
$(4,3,\ldots,3)$. 

By allowing these integers to be $2$, this notation can be extended to
cover certain reducible Coxeter diagrams, namely, those with
underlying graph a disjoint union of subintervals $[1,n+1]$.  For
example, the Schl\"{a}fli symbol $(2,\ldots,2)$ should be understood
as representing the diagram consisting of $n+1$ vertices and no edges,
i.e., $\bA_1\times\cdots\times\bA_1$ ($\cong(\bbZ_2)^{n+1}$).

The Schl\"{a}fli symbols that we will be interested in correspond to
Coxeter systems that are either simplicial (c.f.,
\ref{ss:simplicial-cs}) or spherical.  Moreover, in the simplicial
case, the diagram is necessarily irreducible.  The Schl\"{a}fli
symbols corresponding to irreducible spherical or simplicial Coxeter
systems are listed in Table~\ref{tab:schlafli}, below.

\begin{table}[ht]
\caption{\label{tab:schlafli}Schl\"{a}fli symbols of irreducible spherical and
simplicial Coxeter systems}
\begin{tabular}{||c||l|c||l|c||l|c||}
\hline
\mbox{dimension}&\multicolumn{2}{c||}{\mbox{spherical}} &
\multicolumn{2}{c||}{\mbox{Euclidean}}
&\multicolumn{2}{c||}{\mbox{hyperbolic}}\\ \cline{2-7}
&\mbox{symbol}&t&\mbox{symbol}&t&\mbox{symbol}&t\\ \hline
2&(3,3)&0&(4,4)&0&(p,q) \mbox{with}&0\\ 
& (4,3)&0&(6,3)&0& $p^{-1}+q^{-1}+2^{-1}<1$&\\
& (5,3)&0& &&&\\ \hline
3&(3,3,3)&3&(4,3,4)&1&(3,5,3)&3\\
& (4,3,3)&2&&        &(5,3,4)&2\\
& (5,3,3)&3&&        &(5,3,5)&3\\
& (3,4,3)&2&&        &&\\ \hline
4&(3,3,3,3)&7&(4,3,3,4)&3&(5,3,3,3)&6\\
& (4,3,3,3)&5&(3,4,3,3)&4&(5,3,3,4)&4\\
&          & &         & &(5,3,3,5)&5\\ \hline
$n\geq 5$&(3,3,\ldots,3)&$a_n$&(4,3,\ldots,3,4)&$c_n$&&\\
&(4,3,\ldots,3)&$b_n$&&&&\\ \hline
\end{tabular}\vspace{.2in}
\[\mbox{where}\;\; a_n=\frac{n(n+1)-6}{2},\;\; b_n=\frac{n(n-1)-2}{2},\;\;
c_n=\frac{(n-2)(n-1)}{2}\]
\end{table}

\subsection{Examples of symmetric associahedral tilings}
\label{ss:symm-assoc-tilings} 

As previously, we let $S=\{1,\ldots,n+1\}$ and $\cR$ be the set of all
subsets of $S$ that correspond to proper subintervals of $[1,n+1]$ of
nonzero length.  In what follows we only consider spherical or
simplicial Coxeter systems that can be described by Schl\"{a}fli
symbols as in the previous subsection.

\begin{example}\label{ex:min-sph}
{\em (Minimal blow-ups of boundaries of Coxeter cells.)}
Suppose $(W,S)$ is spherical and irreducible with Schl\"{a}fli symbol
$(m_1,\ldots,m_n)$ (see Table~\ref{tab:schlafli}) and
that $Z=Z(m_1,\ldots,m_n)$ is the corresponding Coxeter cell.  Then
$\cR$ is the set for the minimal blow-up $(\partial Z)_{\#}$ of
$\partial Z$.  Hence, $(\partial Z)_{\#}$ is tiled by associahedra (as
is its universal cover).
\end{example}

\begin{example}\label{ex:nonmin-sph}
{\em (Nonminimal blow-ups of boundaries of Coxeter cells.)}
Suppose that $(m_1,\ldots,m_n)$ is the Schl\"{a}fli symbol for a finite
Coxeter group $W$ (not necessarily irreducible).  For $n=1$ or $2$,
the set $\cR$ is always admissible (cf.,
Definition~\ref{def:admissible}); however, for $n\geq 3$, it might not
be.  In fact, for $n\geq 3$, $\cR$ is admissible if and only if any
time a $2$ occurs in $(m_1,\ldots,m_n)$, the numbers before and after
it are both even.  (This can be checked by using
Remark~\ref{rem:a-is-longest-word}.)  For example, if
$(m_1,\ldots,m_n)$ is $(2,\ldots,2)$ or
$(3,\ldots,3,4,2,4,3,\ldots,3)$, then $\cR$ is admissible.  For any
Schl\"{a}fli symbol such that $W$ is finite and $\cR$ is admissible,
the $\cR$-blow-up of $\partial\cc(m_1,\ldots,m_n)$ will be tiled by
associahedra.
\end{example}

\begin{example}\label{ex:min-simp}
{\em (Minimal blow-ups of simplicial Coxeter systems.)}  Let
$(m_1,\ldots,m_n)$ be a Schl\"{a}fli symbol of a simplicial Coxeter
system (see Table~\ref{tab:schlafli}).  Then $\cR$
is the set for the minimal blow-up of $\Sigma$
($=\Sigma(m_1,\ldots,m_n)$).  Hence $\Sigma_{\#}$ is tiled by
associahedra.  
\end{example}

\begin{notation}
Given a Schl\"{a}fli symbol $(m_1,\ldots,m_n)$ as above, let
$\ccc(m_1,\ldots,m_n)$ denote the universal cover of the $\cR$-blow-up
described in either Example~\ref{ex:min-sph}, \ref{ex:nonmin-sph}, or
\ref{ex:min-simp}.  Also, let $A(m_1,\ldots,m_n)$ denote the symmetry
group of the natural framing on $\ccc(m_1,\ldots,m_n)$.
\end{notation}

\begin{remark}
In Example~\ref{ex:nonmin-sph}, when the Schl\"{a}fli symbol is
$(2,\ldots,2)$, $\cR$ is admissible.  In this case, each gluing
involution $j_T$ is trivial.  Hence, the universal cover
$\ccc(2,\ldots,2)$ is the usual reflection type tiling corresponding
to $(V,M,L(K^n))$.
\end{remark} 

In the next four theorems we classify some of these examples
$X(m_1,\ldots,m_n)$ up to isomorphism.  In fact, in the first three
theorems we classify all such examples for $n\leq 4$.  In the last
theorem, we show that in each dimension $\geq 5$, the three
irreducible examples are distinct.

The basic method for showing that two such tilings are not isomorphic
is to use Lemma~\ref{lem:t-condition}.  In fact, in most cases, the
number $t_{\ccc}$ of mirrors (i.e., codimension-one faces) of the
associahedron that have nonextendable gluing involutions is sufficient
to distinguish among the examples.  The number $t_{\ccc}$ is easily
computable.  It is the number of subsets $T$ of $\{1,\ldots,n+1\}$ such
that $\Gamma_T$ is connected,  $1<\card(T)<n$ and such that $j_T$ is not 
the antipodal map (cf. Remark~\ref{rem:a-is-longest-word}).  (The
numbers $t_{\ccc}$ are also given in Table~\ref{tab:schlafli}.)
Conversely, the basic method for 
showing that two such tilings are isomorphic is to use
Proposition~\ref{prop:can-iso}. 

\begin{theorem}\label{thm:2d-assoc-class}
{\em (Dimension $2$.)}  Let $(m_1,m_2)$ be a pair of integers $\geq
2$.  Then $\ccc(m_1,m_2)$ is isomorphic to the reflection tiling
$\ccc(2,2)$.  Thus, all $2$-dimensional examples are isomorphic to the
tiling of the hyperbolic plane by right-angled pentagons.
\end{theorem}

\begin{proof}{}
In dimension $2$, all gluing involutions are extendable.
\end{proof}

\begin{theorem}\label{thm:3d-assoc-class}
{\em (Dimension $3$.)}  For $n=3$, the universal covers of the
examples in \ref{ex:min-sph}, \ref{ex:nonmin-sph}, and
\ref{ex:min-simp} fall into four isomorphism classes as indicated
below.  (The corresponding values of $t=t_{\ccc}$ are indicated in
parentheses.  Also, $p$ and $q$ denote arbitrary even integers $\geq
2$.)
\begin{enumerate}
\item[(i)] ($t=0$): all three integers are even, i.e., $\ccc(p,2,q)$
or $\ccc(2,p,2)$.
\item[(ii)] ($t=1$): $\ccc(4,3,4)\cong\ccc(2,4,3)$.
\item[(iii)] ($t=2$): $\ccc(4,3,3)\cong\ccc(4,3,5)\cong\ccc(3,4,3)$. 
\item[(iv)] ($t=3$): $\ccc(3,3,3)\cong\ccc(5,3,3)\cong\ccc(3,5,3)\cong
\ccc(5,3,5)$.
\end{enumerate}
\end{theorem}

\begin{proof}{}
The four isomorphism types can be distinguished by the indicated value
of $t$.  The existence of the indicated isomorphisms follows from
Proposition~\ref{prop:can-iso} (possibly after changing one of the
framings by an automorphism of framing systems).
\end{proof}

\begin{theorem}\label{thm:4d-assoc-class}
{\em (Dimension $4$.)}  For $n=4$, the universal covers of the
examples in \ref{ex:min-sph}, \ref{ex:nonmin-sph}, and
\ref{ex:min-simp} fall into nine isomorphism classes as indicated
below.  (The corresponding values of $t=t_{\ccc}$ are indicated in
parentheses.  Also, $p$ and $q$ denote arbitrary even integers $\geq
2$.)
\begin{enumerate}
\item[(i)] ($t=0$): all four integers are even, i.e., $\ccc(2,p,2,q)$
($\cong\ccc(p,2,2,q)$).
\item[(ii)] ($t=1$): $\ccc(3,4,2,p)$.
\item[(iii)] ($t=3$): $\ccc(4,3,3,4)\cong\ccc(2,4,3,3)$.
\item[(iv)] ($t=4$): $\ccc(3,3,4,3)$.
\item[(v)] ($t=4$): $\ccc(5,3,3,4)$.
\item[(vi)] ($t=5$): $\ccc(4,3,3,3)$.
\item[(vii)] ($t=5$): $\ccc(5,3,3,5)$.
\item[(viii)] ($t=6$): $\ccc(5,3,3,3)$.
\item[(ix)] ($t=7$): $\ccc(3,3,3,3)$.
\end{enumerate}
\end{theorem}

\begin{proof}{}
The only question is to distinguish the example in (iv) from (v) and
to distinguish (vi) from (vii).  In $\ccc(3,3,4,3)$ there are three
mirrors of type $K^1\times K^2$ with extendable gluing involutions (in
fact the identity maps), namely, the mirrors corresponding to $(3,4)$,
$(4,3)$ and $(4)$.  The mirror corresponding to $(4)$ intersects the
other two.  In $\ccc(5,3,3,4)$ there are also three such mirrors
corresponding to $(5,3)$, $(3,4)$, and $(4)$.  The mirror
corresponding to $(5,3)$ is disjoint from the other two.  Hence, there
is no automorphism of $K^4$ that takes the first set of mirrors into
the second.  So by Lemma~\ref{lem:t-condition},
$\ccc(3,3,4,3)\not\cong\ccc(5,3,3,4)$.  For a similar reason,
$\ccc(4,3,3,3)\not\cong\ccc(5,3,3,5)$.
\end{proof}

In dimensions $>4$, we shall only consider the examples coming from
irreducible Coxeter systems.

\begin{theorem}\label{thm:>4d-assoc-class}
{\em (Dimension $n>4$.)}  For $n>4$, the universal covers of the
minimal blow-ups in \ref{ex:min-sph} and \ref{ex:min-simp}, namely
$\ccc(3,\ldots,3)$, $\ccc(4,3,\ldots,3)$, and $\ccc(4,3,\ldots,3,4)$
are mutually non-isomorphic.
\end{theorem}

\begin{proof}{}
The corresponding values of $t$ given as $a_n$, $b_n$, and $c_n$ in
Table~\ref{tab:schlafli} are distinct for each $n\geq 5$.
\end{proof}

\begin{remark}
If $\ccc(m_1,\ldots,m_n)\cong\ccc(m_1',\ldots,m_n')$, then the
corresponding groups $A(m_1,\ldots,m_n)$ and $A(m_1',\ldots,m_n')$ are
commensurable.
\end{remark}

More generally, we pose the following question.

\begin{question}
When do the different associahedral tilings of
Theorems~\ref{thm:3d-assoc-class}, \ref{thm:4d-assoc-class}, and
\ref{thm:>4d-assoc-class} give commensurable mock reflection groups?
When do they give quasi-isometric mock reflection groups?  (We answer
the $3$-dimensional quasi-isometry question below in
\ref{ss:3-d-qiso}.) 
\end{question} 

\subsection{Maximally symmetric associahedral tilings}

Let $X$ be one of the associahedral tilings discussed above.  Recall
that $d\mapsto F(d)$ defines a bijection between the set of diagonals
in $P_{n+3}$ and the set of codimension-one faces of $K^n$.  Let $\sD$
denote the set of diagonals, and let $j_d$ denote the gluing
involution on the face $F(d)$.  Then by
Proposition~\ref{prop:max-symm}, we know that $X$ will be maximally 
symmetric if and only if for all $\phi\in\Aut(K^n)$ ($\cong D_{n+3}$)
and $d\in \sD$, the composition $\phi\circ 
j_d\circ\phi^{-1}\circ(j_{\phi(d)})^{-1}$ is the restriction of an
element of $\Aut(K^n)$.  For all of the tilings in
\ref{ss:symm-assoc-tilings} except $X(2,2,\ldots,2)$ and
$X(3,3,\ldots,3)$ there exists a symmetry $\phi\in\Aut(K^n)$ that
conjugates a nonextendable gluing involution to an extendable one,
hence these cannot be maximally symmetric.  In the case of
$X(2,2,\ldots,2)$ the tiling is of reflection type, so we already know
it is maximally symmetric (Example~\ref{ex:max-symm-refl}). Moreover,
its symmetry group is  
\[\Aut(X)=A\rtimes D_{n+3}.\]

In the case of $X(3,3,\ldots,3)$, each gluing involution $j_d$
($=j_{d(T)}$) is the involution $i_T$ described in the proof of
Lemma~\ref{lem:mock-or-not}.  Recall the face  
$F(d')$ is adjacent to $F(d)$ if and only if the diagonals $d'$ and $d$ do not
cross.   Letting $\sD_d$ denote the set of diagonals that do not cross
$d$, we see that $j_d:F(d)\rightarrow F(d)$ induces an involution
(which we also denote by $j_d$) on $\sD_d$.  Let $e$ denote the edge 
$\{0,n+2\}$ of $P_{n+3}$.  Then the involution $j_d:\sD_d\rightarrow\sD_d$
is given by 
\[j_d(d')=\left\{\begin{array}{ll}
r_{L(d)}(d') & \mbox{if $d'$ and $e$ are on opposite sides of $d$}\\
d' &\mbox{if $d'$ and $e$ are on the same side of
$d$}\end{array}\right.\]
(where $L(d)$ is as in Figure~\ref{fig:stab-action}).  Now suppose
$\phi\in D_{n+3}$.  If 
$d'$ and $\phi(d)$ do not cross, then  
\[(\phi\circ j_d\circ\phi^{-1})(d')=\left\{\begin{array}{ll} 
j_{\phi(d)}(d') & \mbox{if $\phi(e)$ and $e$ are on the same side of
$d$}\\
(r_{L(\phi(d))}\circ j_{\phi(d)})(d') & \mbox{if $\phi(e)$ and $e$ are
on opposite sides of $d$.}\end{array}\right.\]
Thus, $\phi\circ j_d\circ\phi^{-1}\circ(j_{\phi(d)})^{-1}$ extends to
an automorphism of $K^n$ (it is either the restriction of $\Id$ or the
restriction of $r_{L(\phi(d))}$), so by
Proposition~\ref{prop:max-symm}, $X(3,3,\ldots,3)$ is maximally symmetric.

In the remainder of this subsection, we will discuss the symmetry
group of the tiling $X=X(3,3,\ldots,3)$.  By
Theorem~\ref{thm:mock-presentation} the subgroup $A$ of $\Aut(X)$ is 
generated by involutions $\alpha_T$ where 
$T$ is a proper subinterval of $[1,n+1]$.  (In what follows we shall
denote this generator by $\alpha_d$ where $d\in \sD$ is the diagonal
corresponding to $T$.)  Since $X$ is maximally symmetric, we also
know that for any element $\phi$ of $D_{n+3}$ there is a unique lift
$\tilde{\phi}$ to $\Aut(X)$ that stabilizes the fundamental tile
$K^n$.  We let $\rho_d:X\rightarrow X$ denote the lift of the reflection 
$r_{L(d)}$ in $D_{n+3}$.  Then the group $\Aut(X)$ is generated by
the $\alpha_d$ and $\rho_d$, $d\in \sD$.  This generating set, however, is
not symmetric with respect to $\Aut(K^n)$.  A more symmetric
generating set arises from the following observation.

\begin{lemma}
The involutions $\rho_d$ and $\alpha_d$ commute.
\end{lemma}

\begin{proof}{}
Let $Q_1$ and $Q_2$ be the two subpolygons of $P$ with diagonal $d$,
and assume the 
edge $e$ (with vertex labels $0$ and $n+2$) is contained in $Q_1$.  Let
$r_1$ and $r_2$ denote the restriction of $r_{L(d)}$ to $Q_1$ and $Q_2$,
respectively.  If $Q_1$ is an $(m_1+3)$-gon and $Q_2$ is an $(m_2+3)$-gon, then
the face $F(d)$ is isomorphic to the product $K^{m_1}\times K^{m_2}$, and
as in the proof of Lemma~\ref{lem:mock-or-not}, the restriction of
$\alpha_d$ to $F(d)$ is $\Id\times\phi(r_2)$.  By
Corollary~\ref{cor:d-stab-action},
the restriction of $\rho_d$ to $F(d)$ is $\phi(r_1)\times\phi(r_2)$.
It follows that the automorphism $(\alpha_d\rho_d)^2$ takes the
fundamental tile $K^n$ to itself and fixes the face $F(d)$ pointwise.
By rigidity, it must be the trivial automorphism.
\end{proof} 

Let $\sS$ denote the set of all subpolygons of $P_{n+3}$, and let
$\delta:\sS\rightarrow \sD$ be the $2$-to-$1$ map that takes each
subpolygon to its corresponding diagonal.  Letting $Q$ be an element of
$\sS$ and $d=\delta(Q)$, we define an involution $\beta_Q:X\rightarrow
X$ by    
\[\beta_Q=\left\{\begin{array}{ll}\alpha_d & \mbox{if
$e\not\subset Q$}\\
\rho_d\alpha_d & \mbox{if
$e\subset Q$.}\end{array}\right.\]
It follows that if $Q_1$ and $Q_2$ are the two subpolygons sharing the
diagonal $d$, then the two involutions $\beta_{Q_1}$ and
$\beta_{Q_2}$ both take the fundamental tile $K^n$ to the adjacent
tile across the face $F(d)$.  As in \ref{ss:group-A}, we obtain
relations among these involutions by considering local pictures around
codimension-$2$ faces  of $K^n$ (or, dually, by considering
$2$-dimensional cells in the dual 
cellulation of $X$ by Coxeter cells--in this case cubes).  Thus,
around any codimension-$2$ face, there are {\em a priori\/} $2^4$
automorphisms of the form
$\beta_{Q_1}\beta_{Q_2}\beta_{Q_3}\beta_{Q_4}$ that take $K^n$ to 
itself, and since the stabilizer of $K^n$ is $D_{n+3}$, there is an
element $\phi$ ($=\phi(Q_1,Q_2,Q_3,Q_4)$) in $D_{n+3}$ such that 
\[\beta_{Q_1}\beta_{Q_2}\beta_{Q_3}\beta_{Q_4}=\tilde{\phi}.\]
We work out these relations explicitly, below.

\begin{figure}[ht]
\begin{center}
\psfrag{Q0}{\scriptsize$Q_0$}
\psfrag{Q1}{\scriptsize$Q_1$}
\psfrag{Q2}{\scriptsize$Q_2$}
\psfrag{Q3}{\scriptsize$Q_3$}
\psfrag{Q4}{\scriptsize$Q_4$}
\psfrag{Q4=f(Q0)}{\scriptsize$Q_4=\phi(Q_0)$}
\psfrag{Q5=f(Q1)}{\scriptsize$Q_5=\phi(Q_1)$}
\psfrag{r1}{\scriptsize$r_1$}
\psfrag{r2}{\scriptsize$r_2$}
\psfrag{r3}{\scriptsize$r_3$}
\psfrag{r4}{\scriptsize$r_4$}
\psfrag{b}{\scriptsize$b$}
\psfrag{c}{\scriptsize$c$}
\psfrag{F1}{}
\psfrag{F2}{}
\psfrag{f}{\scriptsize$\phi$}
\psfrag{1}{$1$}
\psfrag{b1}{$\beta_{Q_1}$}
\psfrag{b1b2}{$\beta_{Q_1}\beta_{Q_2}$}
\psfrag{b1b2b3}{$\beta_{Q_1}\beta_{Q_2}\beta_{Q_3}$}
\psfrag{b1b2b3b4}{$\beta_{Q_1}\beta_{Q_2}\beta_{Q_3}\beta_{Q_4}$}
\includegraphics[scale = .5]{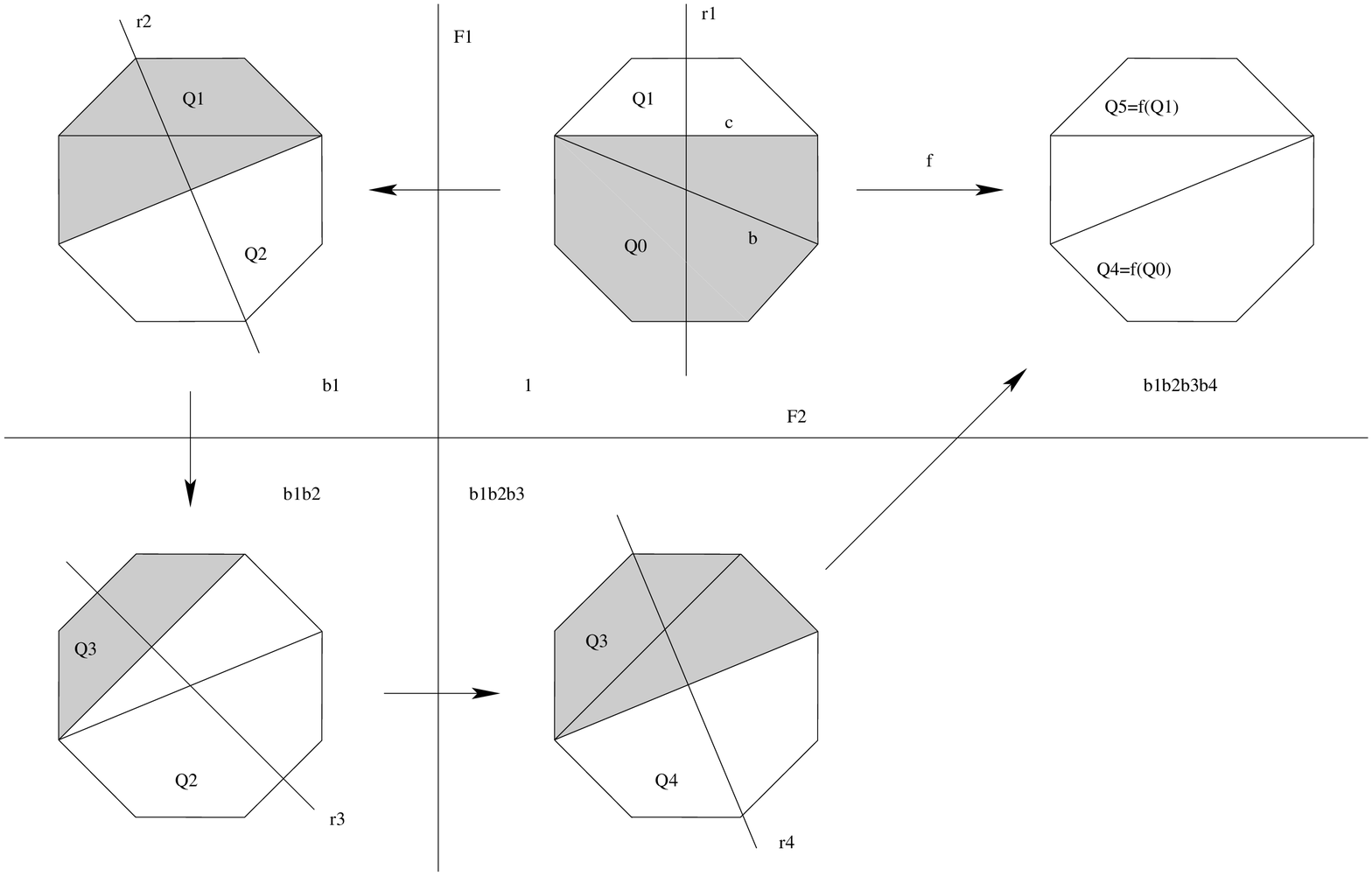}
\caption{\label{fig:symmpres}}
\end{center}
\end{figure}

Suppose $b,c\in \sD$ is a pair of noncrossing diagonals, and
$I=(i_1,\ldots,i_4)$ is a $4$-tuple in $(\bbZ_2)^4$.  Let $B$
denote the subpolygon of $P_{n+3}$ such that $\delta(B)=b$ and $c\not\subset
B$, and let $C$ denote the subpolygon such that $\delta(C)=c$ and
$b\not\subset C$.  We define a sequence of subpolygons inductively by 
\[\begin{array}{l}Q_0=B\\ Q_1=C\\ Q_2=(r_1)^{i_1}(Q_0)\\ Q_3=(r_2)^{i_2}(Q_1)\\
Q_4=(r_3)^{i_3}(Q_2)\\Q_5=(r_4)^{i_4}(Q_3)\end{array}\]
where $r_i$ denotes the reflection of $P_{n+3}$ that takes the
subpolygon $Q_i$ to itself.  The product
$\beta_{Q_1}\beta_{Q_2}\beta_{Q_3}\beta_{Q_4}$ will take the
fundamental tile of $X$ to itself, and the resulting automorphism of
$K^n$ corresponds to the element $\phi\in D_{n+3}$ defined by
$\phi(Q_0)=Q_4$ and $\phi(Q_1)=Q_5$.  In other words,
\[\tilde{\phi}=(\rho_c)^{i_1}(\rho_b)^{i_2}
(\rho_c)^{i_3}(\rho_b)^{i_4}.\]    
Figure~\ref{fig:symmpres} shows an example with $n=5$ and
$I=(1,1,0,1)$.  The fundamental tile $K^n$ is labeled $1$, and the tile
$\beta K^n$ is labeled $\beta$. ($\beta$ is just one of the possible
labels on the tile $\beta K^n$, since any label of the form
$\beta\tilde{\phi}$ where $\phi\in D_{n+3}$ describes the same tile).
The shading indicates which side of the fixed diagonal is affected by
the next gluing involution. The shaded side contains the other
diagonal if and only if the corresponding element of
$(i_1,i_2,i_3,i_4)$ is $1$.  

Let $R_I(b,c)$ denote the word
\[R_I(b,c)=\beta_{Q_1}\beta_{Q_2}\beta_{Q_3}\beta_{Q_4}(\rho_b)^{i_4}
(\rho_c)^{i_3}(\rho_b)^{i_2}(\rho_c)^{i_1}.\]
Letting $m(b,c)$ denote the order of the rotation $r_{b}r_{c}$
in $D_{n+3}$, we then have the following.

\begin{theorem}\label{thm:symmpres}
The group $\Aut(X)$ has a presentation with generators $\beta_Q$,
$Q\in\sS$, and $\rho_d$,$d\in\sD$, and relations: 
\[\begin{array}{ll}
(\beta_Q)^2 & \mbox{for all $Q\in\sS$}\\
(\rho_d)^2 & \mbox{for all $d\in \sD$}\\
(\beta_{Q_1}\beta_{Q_2})^2 & \mbox{whenever
$\delta(Q_1)=\delta(Q_2)$}\\
(\beta_{Q}\rho_d)^2 & \mbox{whenever $\delta(Q)=d$}\\
(\rho_b\rho_c)^{m(b,c)} &\mbox{for all $b,c\in \sD$}\\
R_I(b,c) & \mbox{for all $I\in(\bbZ_2)^4$ and noncrossing diagonals
$b,c$}\end{array}\]
\end{theorem}

\begin{remark}
The generators $\rho_d$ can be eliminated from the presentation, since  
$\rho_d=\beta_{Q_1}\beta_{Q_2}$ if $Q_1$ and $Q_2$ are the two
subpolygons that share the diagonal $d$. 
\end{remark}  
       
Let $S_{n+3}$ be the group of permutations on the set
$\{0,1,\ldots,n+2\}$ (i.e., the set of vertex labels for $P_{n+3}$).
Then for any $d\in\sD$, the reflection $r_{L(d)}$ induces an
involution $\bar{\rho}_d\in S_{n+3}$.  Similarly, for any $Q\in\sS$,
we obtain an involution $\bar{\beta}_Q\in S_{n+3}$ as follows.  Let
$a_0,a_1,\ldots,a_{k+1}$ be labels on the vertices of $Q$ ordered
sequentially with $a_0$ and $a_{k+1}$ being the vertices of the
diagonal $\delta(Q)$.  Then $\bar{\beta}_Q$ is the involution that
reverses the order of the sequence $a_1,a_2,\ldots,a_k$.  

\begin{proposition}
There is a surjective homomorphism $\psi:\Aut(X)\rightarrow S_{n+3}$
defined by $\rho_d\mapsto\bar{\rho}_d$,
$\beta_Q\mapsto\bar{\beta}_Q$. 
\end{proposition}

\begin{proof}{}
All of the relations in Theorem~\ref{thm:symmpres} hold for
$\bar{\rho}_d$ and $\bar{\beta}_Q$.
\end{proof}

\begin{remark}
Let $M^n$ denote the minimal blow-up of the projectivized braid
arrangement in $\bbR\bbP^n$.  (That is, $M^n=(\partial Z_{\#})/a_{\#}$
where $Z$ is the Coxeter cell of type $\bA_{n+1}$.)  Then $M^n$ can be
identified with the real points of the moduli space
$\overline{\cM}_{0,n+3}$ (see \cite{Ka1,Ka2}), and the
$S_{n+3}$-action on $\overline{\cM}_{0,n+3}(\bbR)$ respects the
Coxeter cell decompositon of $M^n$.  The homomorphism
$\psi:\Aut(X)\rightarrow S_{n+3}$ arises when one lifts the
$S_{n+3}$-action to the universal cover $X=\tilde{M}^n$.  In
particular, $\pi_1(M^n)=\ker(\psi)$.
\end{remark}

\subsection{The $3$-dimensional examples}\label{ss:3-d-qiso}
In this subsection we discuss the question of when two $3$-dimensional
tilings are quasi-isometric.  Nowadays any such discussion should be
within the context of Thurston's Geometrization Conjecture.  A closed
orientable irreducible $3$-manifold with infinite fundamental group
has a canonical ``JSJ-decomposition'' into ``simple pieces'' and
Seifert fibered pieces.  Each such piece is a compact $3$-manifold
with boundary, and each boundary component is a torus.  Thurston's
Conjecture is that each simple piece is hyperbolic.  By definition, a
compact $3$-manifold $M^3$ is a {\em hyperbolic piece} if its interior
is homeomorphic to a complete hyperbolic $3$-manifold of finite
volume.  Each boundary component then has a collared neighborhood $C$
such that each component of the inverse image of $C$ in $\bbH^3$ is a
horoball.  Identifying $M^3$ with the complement of a collared
neighborhood of the boundary, we obtain an identification of its
universal cover $\tilde{M}^3$ and $\bbH^3$ with all these horoballs
chopped off.  Such an $\tilde{M}^3$ is called a {\em neutered
hyperbolic space}.

In the case of $3$-manifolds that are tiled by associahedra or
permutohedra, it turns out that (1) each piece in the
JSJ-decomposition is hyperbolic, and (2) the neutered hyperbolic
spaces that arise as universal covers of hyperbolic pieces are all
identical.  The question of whether the universal covers of two such
tilings are quasi-isometric then comes down to the question of whether
or not the lifts of certain gluing involutions extend to
quasi-isometries of the neutered hyperbolic space.  It turns out,
somewhat surprisingly, that the lift of such a gluing involution
extends to an an isometry of $\bbH^3$ that commensurates the lattice
associated to the neutered hyperbolic space.  (Hence, the gluing
involution extends to a quasi-isometry of the neutered hyperbolic
space.)  

\begin{theorem}\label{thm:q-iso}
Suppose $A_1$ and $A_2$ are mock reflection groups associated either
the $3$-dimensional permutohedral tiling $\Sigma_{P^3}$ (see
Theorems~\ref{thm:max-boundary} and \ref{thm:max-simplicial}) or to
one of the $3$-dimensional associahedral tilings of
Theorem~\ref{thm:3d-assoc-class}.  Then $A_1$ and $A_2$ are 
quasi-isometric. 
\end{theorem}

Before proving this theorem we need to develop some notation.  Let
$\bbR^{3,1}$ denote Minkowski space, that is, it is a $4$-dimensional
real vector space with coordinates $x=(x_1,x_2,x_3,x_4)$, equipped
with the indefinite bilinear form defined by 
\[\langle x,y\rangle=x_1y_1+x_2y_2+x_3y_3-x_4y_4.\]
Hyperbolic $3$-space $\bbH^3$ can be defined as one sheet of the
hyperboloid $\langle x,x\rangle=-1$, defined by $x_4>0$.  Let $O(3,1)$
denote the isometry group of the bilinear form and let $O_+(3,1)$ be
the index-two subgroup that preserves the sheets of the hyperboloid.
Then $O_+(3,1)$ is the isometry group of the Riemannian manifold
$\bbH^3$. 

Given a spacelike vector $v\in\bbR^{3,1}$ (i.e., a vector $v$ with
$\langle v,v\rangle>0$), define a reflection $r_v\in O_+(3,1)$ by the
formula 
\[r_v(x)=x-2\frac{\langle x,v\rangle}{\langle v,v\rangle}v.\]

\begin{remark}
Let $O_+(3,1;\bbZ)$ denote the subgroup of $O_+(3,1)$ that preserves
the standard integer lattice $\bbZ^4\subset\bbR^{3,1}$
(i.e., $O_+(3,1;\bbZ)=O_+(3,1)\cap GL_4(\bbZ)$).  If $v\in\bbZ^4$ and
if $\langle v,v\rangle =1$ or $2$, then $r_v$ preserves $\bbZ^4$,
i.e., $r_v\in O_+(3,1;\bbZ)$.
\end{remark}

Suppose $P$ is the permutohedron or the associahedron.  It turns out
that after collapsing each rectangular face of $P$ to a vertex, one ends up
with a polytope $Q$ that can be realized as a right-angled convex
polytope in $\bbH^3$ of finite volume.   Moreover, the vertices
of $Q$ corresponding to the collapsed faces of $P$ will be ideal
vertices of the realization.  This is a special case of a well-known
theorem of Andreev, but we shall verify it directly.  If $P$ is a
permutohedron, then $Q$ is an octahedron.  If $P$ is an associahedron,
then $Q$ is a double pyramid on a triangular base.  (In other words,
$Q$ is the suspension of a triangle.)  When we write down specific
realizations of these polytopes, we find an interesting surprise: the
normal vectors to their faces are integral vectors $v$ satisfying
$\langle v,v\rangle=1$ or $2$.  In fact, consider the following four
polytopes in $\bbH^3$.  

\begin{enumerate}
\item[$\bullet$] {\em The fundamental $3$-simplex} $Q_0$.  Normal
vectors to the faces are $u_1=(1,1,1,1)$, $v_1=(1,0,0,0)$,
$w_1=(1,-1,0,0)$ and $t_1=(0,1,-1,0)$.  $Q_0$ has one ideal vertex,
where the faces normal to $v_1$, $w_1$, and $t_1$ meet.  The Coxeter
group generated by the reflections across its faces has Coxeter
diagram 
\begin{figure}[ht]
\begin{center}
\psfrag{4}{$4$}
\includegraphics[scale=.6]{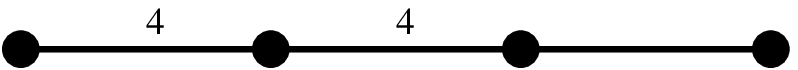}
\end{center}
\end{figure}

\item[$\bullet$] {\em The pyramid} $Q_1$.  Normal vectors to the faces
are $u_1=(1,1,1,1)$, $v_1=(1,0,0,0)$, $v_2=(0,1,0,0)$, and
$v_3=(0,0,1,0)$.  It is again a $3$-simplex, this time with $3$ ideal
vertices.  The corresponding Coxeter diagram is 
\begin{figure}[ht]
\begin{center}
\psfrag{4}{$4$}
\includegraphics[scale=.6]{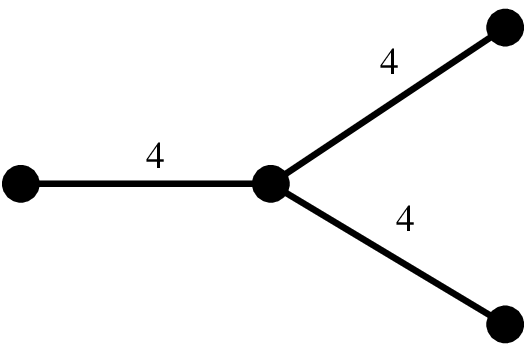}
\end{center}
\end{figure} 

\noindent The {\em base} of the pyramid is the face normal to $u_1$.  It is an
ideal triangle.  The other three faces meet at the vertex
$(0,0,0,1)\in\bbH^3$. 

\item[$\bullet$] {\em The double pyramid} $Q_2$.  The reflection
$r_{u_1}$ carries the vectors $v_1$, $v_2$, and $v_3$ into
$v_1'=(0,-1,-1,-1)$, $v_2'=(-1,0,-1,-1)$, and $v_3'=(-1,-1,0,-1)$,
respectively.  The normal vectors to $Q_2$ are then
$v_1,v_2,v_3,v_1',v_2',v_3'$.

\item[$\bullet$] {\em The regular ideal octahedron} $Q_3$.  The normal
vectors to the faces are the eight vectors $(\pm 1,\pm 1,\pm 1,1)$.
\end{enumerate}

\begin{observations}\hspace{.1in}
\begin{enumerate}
\item[(1)] $Q_0\subset Q_1$, $Q_1\subset Q_2$, and $Q_1\subset Q_3$.
\item[(2)] The Coxeter group generated by reflections across the faces
of $Q_0$ is $O_+(3,1;\bbZ)$.
\item[(3)] The symmetry group of the associahedron (i.e., the dihedral
group $D_6$ of order $12$) acts on $Q_2$, and $Q_0$ is a fundamental
domain. 
\item[(4)] The symmetry group of the permutohedron (i.e.,
$S_4\times\bbZ_2$) acts on the octahedron $Q_3$ as its full symmetry
group.  Again $Q_0$ is a fundamental domain.  ($Q_0$ is a simplex in
the barycentric subdivision of $Q_3$). 
\end{enumerate}
\end{observations}

We are now in a position to prove Theorem~\ref{thm:q-iso}.  Let $X_1$
and $X_2$ be associahedral tilings corresponding to groups $A_1$ and
$A_2$.  Consider a nonextendable gluing involution $i$ defined on a
face of the 
$3$-dimensional associahedron.  By Lemma~\ref{lem:mock-or-not}, the
face is rectangular.  By the proof of Lemma~\ref{lem:mock-or-not}, $i$
is a reflection of the rectangular face about a line of symmetry
connecting the midpoints of two opposite edges.  Any such rectangular
face corresponds to an ideal vertex of the double pyramid $Q_2$.  For
the sake of definiteness, let us fix this vertex to be the one where
the faces normal to $v_1$ and $v_2$ intersect the base triangle
(normal to $u_1$).  The reflection $r_w$ defined by the vector
$w=(0,2,1,1)$ then has the desired effect -- its restriction to the
corresponding rectangular face in the horosphere is $i$.  Since
$\langle w,w\rangle=4$, $r_w$ is not represented by an integral
matrix, rather its entries are rational numbers with denominators at
most $2$.  It follows that $r_w$ commensurates $O_+(3,1;\bbZ)$.  (In
fact, conjugation by $r_w$ maps the congruence subgroup, consisting of
all matrices that are congruent to the identity mod $2$, into itself.) 

Let $\Omega$ denote the neutered hyperbolic space for
$O_+(3,1;\bbZ)$.  The isometry $r_w$ does not quite map $\Omega$ into
itself ($r_w$ does map the set of lifts of all cusps into itself,
but it might not preserve their horoball neighborhoods).  However, it
can be modified to a homeomorphism preserving $\Omega$ that
extends the gluing involution $i$ and that is a bounded distance from
$r_w$.  The hyperbolic pieces of $X_i$, $i=1,2$, give a partition into
copies of $\Omega$, and these copies are glued together via
quasi-isometries.  Hence, $X_1$ is quasi-isometric to $X_2$. 
Similarly, $X_1$ and $X_2$ are both quasi-isometric to the
simply-connected, symmetric permutohedral tiling.  This completes the
proof of Theorem~\ref{thm:q-iso}.

\end{document}